\documentclass[12pt]{amsart}
\usepackage{amsmath,amsthm,amsfonts,amssymb}
\usepackage{graphicx}
\usepackage{mathrsfs}
\usepackage{enumitem}
\usepackage{etoolbox}
\usepackage{xcolor}
\apptocmd{\sloppy}{\hbadness 10000\relax}{}{}
\apptocmd{\sloppy}{\vbadness 10000\relax}{}{}
\usepackage[pdfpagelabels]{hyperref}
\usepackage[letterpaper,margin=1.1in]{geometry}
\usepackage{stmaryrd}

\numberwithin{equation}{section}
\theoremstyle{plain}

\newtheorem{theorem}{Theorem}[section]

\newtheorem{porism}[theorem]{Porism}
\newtheorem{corollary}[theorem]{Corollary}

\newtheorem{lemma}[theorem]{Lemma}

\theoremstyle{definition}
\newtheorem{remark}[theorem]{Remark}
\newtheorem{definition}[theorem]{Definition}

\newcommand{\norm}[1]{\left\lVert#1\right\rVert}

\def\Xint#1{\mathchoice
{\XXint\displaystyle\textstyle{#1}}%
{\XXint\textstyle\scriptstyle{#1}}%
{\XXint\scriptstyle\scriptscriptstyle{#1}}%
{\XXint\scriptscriptstyle\scriptscriptstyle{#1}}%
\!\int}
\def\XXint#1#2#3{{\setbox0=\hbox{$#1{#2#3}{\int}$ }
\vcenter{\hbox{$#2#3$ }}\kern-.6\wd0}}

\def\dashint{\Xint-}

\newcommand{\dist}{\mathop\mathrm{dist}\nolimits}

\newcommand{\spt}{\text{spt}}

\newcommand{\res}{\hbox{ {\vrule height .22cm}{\leaders\hrule\hskip.2cm} }}

\title{Critical Allard regularity: pointwise tilt-excess estimates}

\author{Sean McCurdy}

\keywords{Allard regularity, tilt-excess, harmonic approximation, varifolds, generalized mean curvature}

\address{Department of Mathematics\\ National Taiwan Normal University\\ Taipei, Taiwan}
\email{smccurdy@ntnu.edu.tw}

\begin{document}

\begin{abstract}
The main results of this paper provide VMO-type estimates for the quadratic tilt-excess on varifolds with critical generalized mean curvature.  These estimates apply to varifolds with ``almost-integral" density which are close to a multiplicity one $m$-disc in a ball in the usual senses.  The class of almost-integral varifolds allows for varifolds with non-perpendicular mean curvature.  Moreover, the estimates hold \emph{uniformly for every point} in a relatively open set in $\text{spt}\norm{V}$ and naturally imply a Reifenberg-type parametrization.  The proof relies upon generalizing the $Q$-valued Lipschitz approximation and Sobolev-Poincar\'e estimates of \cite{Menne10_Sobolev} to almost-integral rectifiable varifolds.
\end{abstract}

\maketitle

\tableofcontents

\section{Introduction}

In this paper we continue the study of varifolds with critical generalized mean curvature. The main results show that rectifiable varifolds with almost-integral density (see Definition \ref{d:almost integral densities}, below) satisfy VMO-type qualitative estimates on the quadratic tilt-excess under the hypotheses of the critical Allard tilt-excess decay regime.  These VMO-type estimates are shown to hold \textit{uniformly and point-wise everywhere} in a relatively open set in $\text{spt}\norm{V}$.  This uniform qualitative behavior allows us to appeal to recent results by \cite{BiZhou2022_bilipschitz} on general chord-arc varifolds to obtain a Reifenberg-type parametrization in the critical case.  Thus, this paper generalizes the results of \cite{Zhou_WillmoreEnergy22} in a natural way to dimensions $m\not=2$ and provides the natural analog for Allard's regularity theorem \cite{Allard72} in the critical case (see Corollary \ref{t: VMO main theorem rect} and Corollary \ref{c:chord arc}).  Along the way, we generalize the $Q$-valued Lipschitz approximation and Poincar\'e estimates of \cite{Menne10_Sobolev} to rectifiable varifolds with almost-integral density.  Though we only apply them with $Q=1$ in this work, these results may be of independent interest.

More precisely, we study the following situation. See Section \ref{s:preliminaries} for definitions. 

\subsection{The situation to be studied} (Varifolds with generalized mean curvature)\label{s: the set up} \textit{Let $1\le m< n$ be integers, $1\le p\le \infty$, and $B_1^n(0) \subset \mathbb{R}^n$.  Let $V$ be an $m$-dimensional rectifiable varifold in $B_1^n(0)$ with weight measure $\norm{V}$ that satisfies the following properties.
\begin{enumerate}
    \item $0 \in \emph{spt}\norm{V}$.
    \item The density $\Theta^m(\norm{V}, x) \ge 1$ for $\norm{V}$-a.e. $x \in \emph{spt}\norm{V}$.
    \item $V$ has locally bounded first variation, and there exists a $\boldsymbol{h}(V, \cdot) \in L^p_{loc}(\norm{V}, \mathbb{R}^{n})$ such that for all $\theta \in C^{\infty}_c(B_{1}^n(0), \mathbb{R}^n)$
\begin{align*}
    \int D\theta(x) \cdot S_{\#} dV(x, S) = -\int \theta(x) \cdot \boldsymbol{h}(V, x)d\norm{V}(x).
\end{align*}
\end{enumerate}}

\begin{definition}\label{d:regularity regime}(The Allard Tilt-excess Decay Regime)
    Let $m,n,p,V$ be as in \ref{s: the set up}, and let $0<\tilde{\gamma}< 1$, $0< \delta_0 \le 1/4$, and $B_{\rho}^n(x) \subseteq B_1^n(0)$. We shall say that $V$ \textit{satisfies the hypotheses of the $(\delta_0, \tilde{\gamma})$-Allard tilt-excess decay regime} in $B_{\rho}^n(x)$ if the following conditions hold.
\begin{enumerate}
\item (Density Bounds) \begin{align*}
    \delta_0 \omega_m \rho^m \le \norm{V}(B_\rho^n(x)) \le (2-\delta_0)\omega_m \rho^m.
\end{align*}
\item (Tilt-Excess Bound) There exists a $T \in G(n, m)$ such that
\begin{align*}
E(x, \rho, T) := \left(\rho^{-m}\int_{B_{\rho}^n(x)\times G(n, m)}|S_\# - T_{\#}|^2dV(\xi,S)\right)^{1/2} & \le \tilde{\gamma}.\end{align*}
\item (Mean Curvature Bound) If $p\ge 2$, $\gamma_{\hbar, p}(x, \rho):= \norm{\boldsymbol{h}(V, \cdot)}_{L^p(\norm{V} \res B_{\rho}^n(x))} \le \tilde{\gamma}^2$. If $1\le p<2$, then $\gamma_{\hbar, p}(x, \rho):= \norm{\boldsymbol{h}(V, \cdot)}_{L^p(\norm{V} \res B_{\rho}^n(x))} \le \tilde{\gamma}^4$.\end{enumerate}
Considering the behavior of $\gamma_{\hbar, p}$ under homothetic rescalings of the varifold $V$, we are led to the following terminology.  When $1\le p<m$ we call this the \emph{sub-critical case}.  When $m=p$ we shall call this the \emph{critical case} or \emph{critical Allard tilt-excess decay regime}.  When $m<p$ we shall call this the \emph{super-critical case}.
\end{definition}

This collection of hypotheses is named after Allard, who first studied the \textit{super-critical} case for rectifiable varifolds when $p \ge 2$ \cite{Allard72}. Among other things, Allard proved that if $V$ satisfies the hypotheses of the super-critical tilt-excess decay regime in $B_1^n(0)$ for $0<\tilde{\gamma}$ sufficiently small and $\Theta^m(\norm{V}, x) \in [1, 1+\tilde{\gamma}]$ for $\norm{V}$-a.e. $x \in \text{spt}\norm{V}$, then $\inf_{T \in G(n,m)}E(x, r, T) \le C r^{1-m/p}$ uniformly in $x$ in a neighborhood of $0$.  This uniform decay upon the tilt-excess allowed him to prove that $\text{spt}\norm{V}$ is the graph of a $C^{1, 1-m/p}$ function in a small neighborhood of $0$. Tilt-excess decay is the key ingredient for varifold regularity in the super-critical case and has been a powerful tool in the field ever since. For example, further studies of the super-critical case may be found in \cite{Duggan_varifold86}\cite{GruterJost86}\cite{Schatzle04}\cite{KolasinskiStrzeleckivonderMosel13}\cite{Bourni16}\cite{BourniVolkmann16}.

In contrast to the super-critical case where good parametrizations are obtainable, the sub-critical case does not admit such parametrization.  Indeed, for $m\ge2$ and any $0<\tilde{\gamma}$ it is easy to construct examples of integral varifolds that satisfy the hypotheses of the sub-critical $(\delta_0, \tilde{\gamma})$-Allard tilt-excess decay regime which have dense (or highly disconnected) support in an open neighborhood of $0$.  Nonetheless, decay rates for the quadratic tilt-excess provide a weak notion of regularity.  Studies of the relationship between the integrability of the generalized mean curvature, $\norm{V}$-almost everywhere decay rates for the quadratic tilt-excess, and $C^2$-rectifiability have been further carried out for \textit{integral} varifolds by Brakke \cite{Brakke78}, Sh\"atzle \cite{Schatzle01}\cite{Schatzle09}, Menne \cite{Menne10_Sobolev}\cite{Menne12_decay}, and Kolasinski and Menne \cite{KolasinskiMenne17_decay}.  The results of these papers provide a beautiful $\norm{V}$-almost everywhere \textit{infinitesimal} varifold regularity theory. In particular, this infinitesimal regularity theory applies to points of higher multiplicity.

The study of the critical case goes back to \cite{Allard72} (see Lemma \ref{t:Allard Lemma 8.4}).  However, more recent approaches to the critical case have focused upon $m=p=2$ and proceeded under a set of hypotheses which is closely related to the Allard tilt-excess decay regime.

\begin{definition}(The Allard Regularity Regime) Let $m,n,p,V$ be as in \ref{s: the set up}, and let $0<\tilde{\gamma}$.  For $B_{\rho}^n(0) \subset B_1^n(0)$ we shall say that $V$ \textit{satisfies the hypotheses of the $(m,p,\tilde{\gamma})$-Allard regularity regime} in $B_{\rho}^n(0)$ if the following conditions hold.
\begin{itemize}
    \item (Mass Bounds)
    $\norm{V}(B_{\rho}^n(0))\le (1+\tilde{\gamma})\omega_m \rho^m.$
    \item (Mean Curvature Bounds) $\gamma_{\hbar, p}(x, \rho) = \norm{\boldsymbol{h}(V, \cdot)}_{L^p(\norm{V} \res B_\rho^n(x))} \le \tilde{\gamma}.$
\end{itemize}
As above, if $1\le p<m$, we call this the \textit{sub-critical} case.  If $p=m$, we call this the \textit{critical case} or the \textit{critical Allard regularity regime}.  If $m<p$, we shall call this the \textit{super-critical case}.
\end{definition}
We note that the $(m,p,\tilde{\gamma})$-Allard regularity regime implies the $(\delta_0, \tilde{\gamma}')$-Allard tilt-excess decay regime for $\tilde{\gamma}$ sufficiently small.  See Remark \ref{r:reg regime implies tilt decay regime}.

For $2$-dimensional surfaces, $\gamma_{\hbar, 2}$ corresponds to the Willmore energy, for which there is a large literature (see for example the recent paper \cite{Pozzetta_Willmore23} and the references therein).  While studying topological constraints for surfaces with finite Willmore energy, Jie Zhou proved the following theorem.

\begin{theorem}\label{t:zhou22}\emph{(\cite{Zhou_WillmoreEnergy22})}
    Let $m=p=2$ and $n,V$ be as in \ref{s: the set up}. Let $0<\tilde{\gamma}$, $0<\alpha<1$, and suppose that $V$ satisfies the hypotheses of the $(2,2,\tilde{\gamma})$-Allard regularity regime in $B_1^{n}(0)$. If $\boldsymbol{h}(V, x)$ is perpendicular to $\emph{ap Tan}_x \emph{spt}\norm{V}$ for $\norm{V}$-a.e. $x \in \emph{spt}\norm{V}$ then there exists a $0<\gamma(n,\alpha)$ such that if $0<\tilde{\gamma}\le \gamma(n, \alpha)$, then there exists an $f: \mathbb{R}^m \rightarrow \mathbb{R}^n$ such that $f, f^{-1} \in C^{0,\alpha}$ and
    \begin{align*}
    \emph{graph}_{B_{\rho/2}^m(0)}f \subset \emph{spt}\norm{V} \cap B_{\rho}^m(0) \subset \emph{graph}_{B_{2\rho}^m(0)}f,
    \end{align*}
    where $\rho = \rho(n) \in (0,1/2)$.
\end{theorem}

We note that the hypothesis that the generalized mean curvature is perpendicular is naturally satisfied by integral varifolds by the work of \cite{Brakke78}, but holds for a strictly larger class than disjoint sums of constant multiples of integral varifolds.

More recently, Yuchen Bi and Jie Zhou \cite{BiZhou2022_bilipschitz} have proven that in the case $V$ is an \emph{integral} varifold and satisfies the hypotheses of the $(2,2,\tilde{\gamma})$-Allard regularity regime in $B_1^{n}(0)$ for sufficiently small $0<\tilde{\gamma}$ there exists a bi-Lipschitz parametrization of $\text{spt}\norm{V}$ in a neighborhood of $0$.  

Unfortunately, the method of proof in both \cite{Zhou_WillmoreEnergy22}, \cite{BiZhou2022_bilipschitz} is restricted to the two-dimensional case.  Theorem \ref{t:zhou22} uses a monotonicity formula from \cite{Simon_willmoreminimizers93} which only holds for $2$-dimensional varifolds with perpendicular generalized mean curvature.  The proof of the bi-Lipschitz parametrization in \cite{BiZhou2022_bilipschitz} uses conformal parametrization (isothermal coordinates) and is therefore not available for higher dimensions.  

\subsection{The results of the present paper.}  In this paper, we continue to study varifold regularity under the hypotheses of the critical $(\delta_0, \tilde{\gamma})$-Allard tilt-excess decay regime for $m\ge 1$.  In particulae, we generalize the density assumption from \cite{Allard72} below in order to obtain results upon certain classes of rectifiable varifolds.

\begin{definition}\label{d:almost integral densities}(Almost-integral densities)
    Let $0<\gamma_{\Theta}<1$. For $V$ as in \ref{s: the set up} we shall say that $V$ has $\gamma_{\Theta}$-\textit{almost integral density} if   
\begin{align*}
    \Theta^m(\norm{V}, x) \in \cup_{n \in \mathbb{N}}[n, n(1+\gamma_{\Theta})] \quad \text{for $\norm{V}$-a.e. } x \in \text{spt}\norm{V}.
\end{align*}
In practice, we shall often refer to such $V$ as having \textit{almost-integral density}, or being \textit{almost-integral} rectifiable varifolds. That is, the quantifier will be suppressed unless needed for the statement of a lemma.
\end{definition}

Now, we state our main results.

\begin{theorem}\label{t:BMO main theorem rect}\emph{(The BMO estimate)}
Let $0<\delta_0\le 1/4$, $0<r \le 1$.  There exist constants
\begin{align*}
    0< \eta_0(m, \delta_0), \quad 0<\gamma(m, n, \delta_0), \quad 0< \gamma_{\Theta}(m, n, \delta_0) \le 1/2 
\end{align*} such that if $V$ has $\gamma_{\Theta}$-almost integral density and satisfies the hypotheses of the critical $(\delta_0, \tilde{\gamma})$-Allard tilt-excess decay regime in $B_r^n(0)$ with constant $\tilde{\gamma} \le \gamma(n,m,\delta_0)$ then \emph{for all} $x \in B_{\eta_0 r}^n(0) \cap \emph{spt}\norm{V}$ and all $0< \rho \le (1-\eta_0)r/3$ 
\begin{align*} 
\inf_{R \in G(n, m)}E(x, \rho, R) \le C(n,m,\delta_0)\tilde{\gamma},
\end{align*}
\end{theorem}

\begin{remark}
    The constant $\eta_0$ is given in Definition \ref{d:eta and eta 0}. The requirement that $\gamma_{\Theta}\le 1/2$ comes from the assumption that $\delta_0 \le 1/4$.  Taking $\delta_0$ smaller allows this crude upper bound to be larger.  However, the existence of $\gamma_{\Theta}(n,m,\delta_0)$ is given implicitly by compactness results and is expected to be very, very small.  
\end{remark}

In the critical case, rescaling $V$ does not imply decay on $\norm{\boldsymbol{h}(V, \cdot)}_{L^m(\norm{V})}$, so no decay rate is obtainable as in the super-critical case.  Instead, we simply have that $\lim_{\rho \rightarrow 0} \gamma_{\hbar, m}(x, \rho) = 0$.  This is gives the following corollary.

\begin{corollary}\label{t: VMO main theorem rect}\emph{(The uniform VMO estimate)}
Under the assumptions of Definition \ref{d:regularity regime} with $0<\eta_0(m, n, \delta_0)$, $0< \gamma(m, n, \delta_0)$, $0< \gamma_{\Theta}(m, n, \delta_0)$ as in Theorem \ref{t:BMO main theorem rect}.  If $V$ has $\gamma_{\Theta}$-almost integral density and satisfies the hypotheses of the critical $(\delta_0, \tilde{\gamma})$-Allard tilt-excess decay regime in $B_r^n(0)$ with constant $\tilde{\gamma} \le \gamma(n,m,\delta_0)$ then
\begin{align*} 
\lim_{\rho \rightarrow 0} \sup_{x \in B_{\eta_0 r}^n(0) \cap \emph{spt}\norm{V}}\inf_{R \in G(n,m)} E(x, \rho, R) = 0.
\end{align*}
\end{corollary}

By looking at Lebesgue points of $\boldsymbol{h}(V, \cdot),\Theta^m(\norm{V}, \cdot)$ one may trivially obtain decay rates for $\norm{V}$-a.e. $x \in B_{\eta_0 r}^n(0) \cap \text{spt}\norm{V}$.

\begin{porism}\label{t:decay rates}\emph{($\norm{V}$-almost everywhere decay rates)}
Let $m\ge 2$ and $0\le \alpha <1$ be given. Let $0<\eta_0(m, n, \delta_0)$, $0<\gamma(n,m,\delta_0)$, and $0< \gamma_{\Theta}(m, n, \delta_0)$ be as in Theorem \ref{t:BMO main theorem rect}. If $V$ has $\gamma_{\Theta}$-almost integral density and satisfies the hypotheses of the critical $(\delta_0, \tilde{\gamma})$-Allard tilt-excess decay regime in $B_r^n(0)$ with constant $\tilde{\gamma} \le \gamma(n,m,\delta_0)$, then for $\norm{V}$-a.e. $x \in \emph{spt}\norm{V} \cap B^n_{\eta_0 r}(0)$
\begin{align*}    \lim_{\rho \rightarrow 0}\inf_{R \in G(n,m)} \rho^{-\alpha}E(x, \rho, R) =0\end{align*}

If $m=1$, we may take $\alpha \in [0,1]$.
\end{porism}

If we restrict to \textit{integral} varifolds, Porism \ref{t:decay rates} partially recovers the results of \cite{Menne12_decay}, which contain sharp decay rates for the quadratic tilt excess (see \cite{Menne12_decay} Theorem 10.2).  However, these decay estimate appear to be new for almost-integral varifolds.  In analogy with \cite{Menne12_decay} Theorem 10.2, we do not expect this to be optimal.

\begin{remark}(Sharpness of the uniform VMO estimate) 
The conclusion of Corollary \ref{t: VMO main theorem rect} does not imply that all points $x \in \text{spt}\norm{V} \cap B_{\eta_0 r}^n(0)$ have unique tangents.  Indeed, the example from \cite{HutchinsonMeier1986} shows that the comparison planes $R \in G(n, m)$ may rotate or oscillate uncontrollably as the scale $\rho \rightarrow 0$ at some points.

On the other hand, the example in \cite{HutchinsonMeier1986} shows that Corollary \ref{t: VMO main theorem rect} is sharp in the following sense.  It is possible for an integral varifold $V$ to satisfy \eqref{s: the set up} and the hypotheses of the critical $(\delta_0, \tilde{\gamma})$-Allard tilt-excess decay regime in $B_{1}^n(0)$ and still have $\inf_{R \in G(n,m)}E(0, \rho, R) \rightarrow 0$ slower than any power.
\end{remark}

In analogy with the parametrization results of \cite{Allard72} and \cite{BiZhou2022_bilipschitz}, Corollary \ref{t: VMO main theorem rect} implies that $B_{\eta_0 r/2}^n(0) \cap \text{spt}\norm{V}$ may be parametrized by bi-$W^{1, k}$ a map in a neighborhood of the origin.  To state the result, we recall the definition of chord-arc varifolds from \cite{BiZhou2022_bilipschitz}.  The term ``chord-arc" follows the terminology of Semmes \cite{Semmes91a}\cite{Semmes91b}\cite{Semmes91c} and Blatt \cite{Blatt09}.

\begin{definition}\label{d:chord arc varifold}
 (Chord-arc varifolds). Let $1\le m \le n$ be integers, $U \subset \mathbb{R}^n$, and $V \in \textbf{RV}_m(U)$.  Let $B_r(x) \subset U$.  We say that $V$ is a chord-arc $m$-varifold in $B_r^n(x)$ with constant $0<\tau$ if the following conditions hold for any $B_\rho^n(z) \subset B_r^n(x)$ with $z \in \text{spt}\norm{V}$.
 \begin{enumerate}
     \item (Ahlfors regularity) \begin{align*}
         1-\tau \le \frac{\norm{V}(B_{\rho}^n(z))}{\omega_m r^m} \le 1+\tau.
     \end{align*}
     \item (Reifenberg Condition) There exists a $T_{z, \rho} \in G(n,m)$ such that 
     \begin{align*}
         \frac{1}{\rho}\dist_{\mathcal{H}}(B_{\rho}^n(z) \cap \text{spt} \norm{V}, B_{\rho}^n(z) \cap (T_{z, \rho} + z)) \le \tau.
     \end{align*}
     \item (Tilt-excess bound) For the $T_{z, \rho} \in G(n,m)$ in the Reifenberg condition
     \begin{align*}
         \rho^{-m}\int_{B_\rho^n(z)}|(T_{z, \rho})_{\#} - S_{\#}|^2 dV(x, S) \le \tau^2.
     \end{align*}
 \end{enumerate}
\end{definition}

Recently, \cite{BiZhou2022_bilipschitz} proved a Reifenberg-type parametrization result for chord-arc varifolds.

\begin{theorem}\label{t:bi_zhou_chord arc}\emph{(\cite{BiZhou2022_bilipschitz}, Theorem 1.5)}
If $V$ is an $m$-dimensional chord-arc varifold in $B_{1}(0)$ with constant $0<\tau$ for $0<\tau$ sufficiently small, $T_{0, 1} \in G(n,m)$ is the $m$-plane guaranteed by Definition \ref{d:chord arc varifold}, and $D_1 \subset T_{0,1}$ the disc of radius $1$ then there exists a map $f: D_1 \rightarrow \emph{spt}\norm{V}$ such that
    \begin{enumerate}
        \item For any $x \in D_1$,
\begin{align*}
|f(x)-x|\le C\tau^{1/4};
\end{align*}
\item For any $x, y \in D_1$,
\begin{align*}
(1-C\tau^{1/4})|x-y|^{1+o(\tau)} \le |f(x)-f(y)|\le (1+C\tau^{1/4})|x-y|^{1-o(\tau)};
\end{align*}
\item $B_{1-o(\tau)}^n(0) \cap \emph{spt}\norm{V} \subset f(D_1)$;
\item For functions $f^*, f_*$ defined by 
\begin{align*}
f^*(x) : = \sup_{y \in D_1}\frac{|f(x) - f(y)|}{|x-y|}, \qquad 
f_*(x) : = \inf_{y \in D_1}\frac{|f(x) - f(y)|}{|x-y|}\\
\end{align*}
and for $k<k(\tau) = \frac{\log(4)}{\log(1+C\tau^{1/4})} \rightarrow +\infty$ as $\tau \rightarrow 0$, there holds
\begin{align*}
\int_{D_1}(f^*)^k + (f_*)^{-k}dx \le C_k, \qquad 
\int_{D_1}|(f-i)^*|^k dx \le C_k \tau^{k/4},
\end{align*}
where $i: D_1 \rightarrow \mathbb{R}^n$ is the inclusion map, $o(\tau) \rightarrow 0$ as $\tau \rightarrow 0$, and the constants $C, C_k$ are independent of $0<\tau$.
\end{enumerate}
\end{theorem}

For $\gamma, \gamma_{\Theta}$ small we may apply Theorem \ref{t:bi_zhou_chord arc} to show that $\text{spt}\norm{V}$ is a topological manifold in a neighborhood of the origin.

\begin{corollary}\label{c:chord arc}\emph{(Parametrization under the critical Allard tilt-excess decay regime)} For any $0<\tau$, there are constants 
\begin{align*}
0<\gamma(m,n,\delta_0, \tau), \quad  0<\gamma_{\Theta}(m,n,\delta_0, \tau)
\end{align*} such that if $V$ has $\gamma_{\Theta}$-almost integral density and satisfies the hypotheses of the critical $(\delta_0, \tilde{\gamma})$-Allard tilt-excess decay regime in $B_1^n(0)$ with constant $\tilde{\gamma}\le \gamma$ then $V$ is an $m$-dimensional chord-arc varifold in $B_{\eta_0 r}^n(0)$ with constant $\tau$.  Thus, for $0<\tau$ sufficiently small we may apply Theorem \ref{t:bi_zhou_chord arc} in $B_{\eta_0 r}^n(0)$ to obtain a bi-$W^{1,k(\tau)}(B_{\eta_0 r}^m(0), \mathbb{R}^n)$ function such that
\begin{align*}
    \emph{graph}_{B_{(1-2 \cdot o(\tau))\eta_0 r}^n(0)} f \subset \emph{spt}\norm{V} \cap B_{(1-o(\tau))\eta_0 r}^n(0) \subset \emph{graph}_{B_{\eta_0 r}^n(0)} f.
\end{align*}
Note that $\gamma(m,n,\delta_0, \tau), \gamma_{\Theta}(m,n,\delta_0, \tau) \rightarrow 0$ as $\tau \rightarrow 0$.
\end{corollary}

Corollary \ref{c:chord arc} provides a direct analog of Allard's results in the critical case.  It also can be viewed as, in some sense, generalizing Theorem \ref{t:zhou22} to $m \not = 2$ in the special case of $\gamma_{\Theta}$-almost integral varifolds. In particular, almost-integral varifolds allow for the generalized mean curvature to be non-perpendicular.  However, the hypotheses of Theorem \ref{t:zhou22} do not assume that the possible densities have a large gap as in the case of $\gamma_{\Theta}$-almost integral densities.

\begin{remark}(On the importance of the almost-integral density assumption)
The assumption that $V$ is $\gamma_{\Theta}$-almost integral for $\gamma_{\Theta}<1$ is essential to the conclusion of Theorem \ref{t:BMO main theorem rect} and the bi-H\"older parametrization following from Corollary \ref{c:chord arc}. This emphasizes the fundamental difference in the super-critical and critical cases.  In \cite{Duggan_varifold86} it is shown that one may remove the density hypotheses $\Theta^m(\norm{V}, x) \in [1, 1+\tilde{\gamma}]$ from Allard's proof (see \cite{Simon83_book} 23.2 for more conditions under which this density hypotheses may be removed).  However, \cite{Menne16} Example 6.11 shows that no such result is possible in the critical case for general rectifiable varifolds with $\Theta^m(\norm{V}, x)\ge 1$.   In particular, for any $0<\tilde{\gamma}$ Menne shows that there exists a rectifiable varifold for which the density is allowed to range in $[1,3)$ for $\norm{V}$-a.e. point and $V$ satisfies the hypotheses of the critical $(\delta_0, \gamma)$-Allard tilt-excess decay regime in $B_1^n(0)$, but $\text{spt}\norm{V}$ is not a topological disc in any neighborhood of the origin.  Similar, though less pathological, examples may be easily produced for rectifiable varifolds for which the density is allowed to range in $[1,2]$ for $\norm{V}$-a.e. point for which the parametrization results of Corollary \ref{c:chord arc} fail.

It would be a very interesting question to see just how large $\gamma_{\Theta}$ may be taken.  Unfortunately, the proof in this paper is not able to shed any light on estimating $\gamma_{\Theta}$.
\end{remark}

\subsection{Discussion of the proofs and outline of the paper.}

Section \ref{s:preliminaries} below gives the basic definitions necessary for the statement and proof of the main results.  Section \ref{s:prelim results} recalls some well-known preliminary results for compactness of varifolds, harmonic functions, and almost weakly harmonic functions.  Additionally, in Remark \ref{r:reg regime implies tilt decay regime}, it is proven that for any $0<\tilde{\gamma}$ sufficiently small, there is a $0< \tilde{\gamma}'$ such that varifolds which satisfy the $\tilde{\gamma}'$-Allard regularity regime also satisfies the hypotheses of the $(\delta_0, \tilde{\gamma})$-Allard tilt-excess decay regime.

The proof of Theorem \ref{t:BMO main theorem rect}, Corollary \ref{t: VMO main theorem rect}, and Porism \ref{t:decay rates} follows the outlines of the proof of Allard's tilt-excess decay theorem (\cite{Allard72} Theorem 8.16) in the super-critical case as rewritten in \cite{Simon83_book}.  Roughly speaking, the proof involves two steps:
\begin{enumerate}
    \item First, under the hypotheses of the tilt-excess decay regime it is shown that at some initial scale a large portion of the varifold $V$ may be captured in the graph of a Lipschitz function.  Importantly, the ``symmetric difference" between the graph and the support of the varifold may be controlled by the tilt-excess.
    \item Then, because the varifold $V$ is close to a minimal surface and minimal surfaces are close to the graphs of harmonic functions, it is shown that the Lipschitz approximation is  very close to a harmonic function.  This closeness is then pushed back onto the support of the varifold to show that $V$ satisfies the conditions of the tilt-excess decay regime at some smaller scale.  Iterating gives the conclusion.
\end{enumerate}
We note that this argument uses ideas from De Giorgi \cite{DeGiorigi61} in the context of codimension 1 area minimizing surfaces and introduced for varifolds by Almgren \cite{Almgren68}. Additionally, Allard-type regularity results have been obtained for minimal surfaces using viscosity type methods which completely avoid the need for the Lipschitz approximation (see \cite{Savin07}).  However, in the interest of greater applicability and providing a proper analog to Allard's results this paper shall parallel his approach.  

From this perspective, the key ingredient is the Lipschitz approximation. In Allard's original Lipschitz approximation, he uses some estimates on the mean curvature of the form
\begin{align*}
    R \norm{\delta V}(B_r^n(x)) \le \tilde{\gamma}\norm{V}(B_r^n(x))
\end{align*}
for $0<r<R$ (see \cite{Allard72} Lemma 6.2).  However, by assuming that $\boldsymbol{h}(V, \cdot) \in L^m(\norm{V})$ we naturally only have $\norm{\delta V}(B_r^n(x)) \le \norm{\boldsymbol{h}(V, \cdot)}_{L^m(B_r^n(x))}\norm{V}(B_r^n(x))^{1-1/m}$.  In \cite{Brakke78} (Lemma 5.4), Brakke obtains a Lipschitz approximation under assumptions of the form $r \norm{\delta V}(B_r^n(x)) \le \tilde{\gamma}\norm{V}(B_r^n(x))$, which is satisfied in the critical case under lower density bounds. More recently, the full $Q$-valued Lipschitz approximation for integral varifolds in the critical case were obtained in \cite{Menne10_Sobolev} Lemma 3.15.  The main technical contribution of this paper is to generalize the $Q$-valued Lipschitz approximation of \cite{Menne10_Sobolev} to almost-integral rectifiable varifolds. This is done in Section \ref{s:lipschitz approx}.
 
Having obtained a $Q$-valued Lipschitz approximation lemma, we apply it twice.  First, we use it to obtain Sobolev-Poincar\'e-type estimates on rectifiable varifolds with almost-integral density (Section \ref{s:sobolev poincare}).  These estimates generalize \cite{Menne10_Sobolev} Theorem 4.4 to almost-integral varifolds and may be of independent interest.  However, we shall only use the Poincar\'e estimate with $Q=1$ to obtain local flatness results in Section \ref{s:height bound}. Second, we use the Lipschitz approximation in step (2) in the tilt-excess decay argument, above.

With these modifications, we follow the tilt-excess decay argument as in \cite{Simon83_book} Theorem 22.5 closely. This is carried out in detail in Section \ref{s:tilt-excess decay} for completeness. Then, the proofs of Theorem \ref{t:BMO main theorem rect}, Corollary \ref{t: VMO main theorem rect}, Porism \ref{t:decay rates}, and Corollary \ref{c:chord arc} follow easily in Section \ref{s:proof of main results} by examining the proof of Theorem \ref{t:BMO main theorem rect}.

\section{Notation}\label{s:preliminaries}

For any $n \in \mathbb{N}$, $A \subset \mathbb{R}^n$, and $0<r<\infty$ we let $B^n_r(A) = \{x \in \mathbb{R}^n: \dist(x, A) < r \}$ and $\overline{B^n_r(A)} = \{x \in \mathbb{R}^n: \dist(x, A) \le r \}$.  For an open set $U\subset \mathbb{R}^n$, we following the notation for distributions, we shall let $\mathscr{D}(U, \mathbb{R}^k)$ be the collection of all $C^{\infty}(U, \mathbb{R}^k)$ functions with support compactly contained in $U$.

\begin{definition}\label{d: multiliear algebra}(Multilinear Algebra)
    Let $V, W$ be vector spaces.  Let $i \in \mathbb{N}$.  We shall write $\odot^i(V, W)$ for the set of all $i$-linear symmetric maps from $V^i$ to $W$. When $i=1$, we shall use the notation $\text{Hom}(V, W) = \odot^1(V, W)$.  

    When $V, W$ are inner product spaces and $\dim(V) <\infty$ we may define an inner product on $\odot^i(V, W)$ as in \cite{Federer69} 1.10.6. When $i =1$ and $V, W =\mathbb{R}^n$, we may realize $\phi \in \text{Hom}(\mathbb{R}^n, \mathbb{R}^n)$ as a matrix. In this case, the inner product is given by 
    \begin{align*}
    \phi \cdot \varphi := \text{trace}(\phi \circ \varphi^T).
    \end{align*}
    
    Assuming $V, W$ are inner product spaces and $\dim(V) = j \in \mathbb{N}$, the norm induced by the inner product on $\odot^i(V, W)$ satisfies the following formula.  For any orthonormal basis $\{e_1, ..., e_j\}$ of $V$
    \begin{align*}
        |\phi|^2 = \frac{1}{i!}\sum_{\sigma \in \mathscr{S}(i,j)} |\phi(e_{s(1)}, ..., e_{s(i)})|^2.
    \end{align*}
    where $\mathscr{S}(i,j)$ is the collection of maps from $\{1, ..., i \}$ into $\{1, ..., j \}$.
\end{definition}

\begin{definition}\label{d:projections}(Grassmanian and Projections)
    Let $1\le m \le n$ be integers.  We write $G(n,m)$ for the set of all $m$-dimensional subspaces of $\mathbb{R}^n$. If $T \in G(n,m)$, we let $T_{\#} \in \text{Hom}(\mathbb{R}^n, \mathbb{R}^n)$ be orthogonal projection onto $T$. We let $T_{\#}^*$ denote the adjoint map, $T^{\perp} = \text{ker}(T_{\#})$, and define $T_{\#}^{\perp} \in \text{Hom}(\mathbb{R}^n, \mathbb{R}^n)$ to be orthogonal projection onto $T^{\perp}$.  Note that $T_{\#}$ satisfies
\begin{align*}
    T_{\#} = T_{\#}^* \qquad T_{\#} \circ T_{\#} = T_{\#} \qquad \text{im}(T_{\#}) = T.
\end{align*}

We endow $G(n,m)$ with the topology induced by the standard metric
\begin{align*}
    \dist_{G}(T, S) = \sup_{v \in T \cap \mathbb{S}^{n-1}}|S_{\#}(v)|.
\end{align*}

We shall use the norm $\norm{T-S} = \dist_{G}(T, S)$.  Note that 
\begin{align*}
    |S_{\#} - T_{\#}| \le \sqrt{n}\norm{S_{\#} - T_{\#}}.
\end{align*}
\end{definition}

\begin{definition}\label{d:varifolds}(Varifolds) Let $1\le m \le n$ be integers, $U \subset \mathbb{R}^n$ be an open set.  An $m$-dimensional \textit{varifold} in $U$ is a Radon measure on $U \times G(n,m)$. The \textit{weight measure} $\norm{V}$ is defined by 
\begin{align*}
    \norm{V}(A) = V(A \times G(n, m))
\end{align*}
    for any $A \subset U$. For any $W \subset U$ we shall use the notation $\norm{V} \res W$ to denote the restriction of $\norm{V}$ to $W$.  We shall use the notation $\textbf{V}_m(U)$ for the set of all $m$-dimensional varifolds in $U$.

For $V \in \textbf{V}_m(U)$ and $\theta \in \mathscr{D}(U, \mathbb{R}^n)$ the \textit{first variation of $V$ with respect to area in the direction $\theta$} is given by
    \begin{align*}
        \delta V(\theta) = \int D\theta \cdot S_{\#} dV(x, S).
    \end{align*}
The first variation defines an associated variational measure
    \begin{align*}
        \norm{\delta V}(A) = \sup\{\delta V(\theta): \theta \in \mathscr{D}(A, \mathbb{R}^n), |\theta|\le 1 \}
    \end{align*}
    for $A \subset U$ an open set. Note that $\norm{\delta V}$ is Borel regular. If $\norm{\delta V}=0$ we say that $V$ is \textit{stationary}.

 We say that a varifold $V \in \textbf{V}_m(U)$ has \textit{locally bounded first variation} if $\norm{\delta V} \res W$ is a Radon measure for every $W \subset \subset U$, in which case the Riesz Representation Theorem gives that there is a vector-valued function $\hbar: \text{spt}\norm{V} \rightarrow \mathbb{S}^{n-1}$ such that 
 \begin{align*}
        \delta V(\theta) & = -\int \theta \cdot \hbar d\norm{\delta V}(x).
    \end{align*}
When $\norm{\delta V} << \norm{V}$, we may define the \textit{generalized mean curvature vector} $\boldsymbol{h}(V; x) \in L^1(\norm{V}, \mathbb{R}^n)$ by 
\begin{align*}
        \boldsymbol{h}(V; x) = \hbar(x) \left(\lim_{r \rightarrow 0} \frac{\norm{\delta V}(B_r(x))}{\norm{V}(B_r(x))}\right),
\end{align*}
    whenever the limit exists.
    
We say that $V \in \textbf{V}_m(U)$ is a \textit{rectifiable varifold} if there is a sequence of $m$-dimensional $C^1$ submanifolds $\Sigma_i$ and non-negative $\mathcal{H}^m$-measurable functions $c_i$ such that for any $f \in C_c(U \times G(n, m); \mathbb{R})$ 
\begin{align*}
V(f) = \sum_i\int_{\Sigma_i} f(x, T_x\Sigma_i)c_i(x)d\mathcal{H}^m(x).
\end{align*}
We shall use the notation $\textbf{RV}_m(U)$ to denote the set of all $m$-dimensional rectifiable varifolds in $U$. We say that $V \in \textbf{RV}_m(U)$ is an \textit{integral varifold} if the $c_i$ take only integer values. Let $\textbf{IV}_m(U)$ denote that set of all $m$-dimensional integral varifolds in $U$.

Note that if $V \in \textbf{RV}_m(U)$ then for $\norm{V}$-almost every $x \in \text{spt}\norm{V}$, $\text{ap Tan}^m(\norm{V}, x)$ (the cone of approximate tangent vectors) satisfies $\text{ap Tan}^m(\norm{V}, x) \in G(n, m)$. For such $x$, we define the function $S(x):= \text{Tan}^m(\norm{V}, x)$. Often we shall use the notation $S, S_{\#}$ to mean $S(x), S(x)_{\#}$.

For $V \in \textbf{RV}_m(U)$ we define the \textit{density function}
\begin{align*}
    \Theta^m(V, x) := \lim_{r \rightarrow 0} \frac{\norm{V}(B_r^n(x))}{\omega_m r^m}.
\end{align*}
Note that $\Theta^m(V, x)$ exists for $\mathcal{H}^m$-almost every $x \in \text{spt}\norm{V}$.

We topologize $\textbf{RV}_m(U)$ with the topology of weak convergence.  That is, we say that a sequence $V_i \rightarrow V$ \text{as varifolds} or \textit{weakly} if for every $f \in C_c(U \times G(n,m), \mathbb{R})$
\begin{align*}
    \int f(x, S)dV_i(x, S) \rightarrow \int f(x, S)dV(x, S).
\end{align*}
\end{definition}

\begin{definition}(Functions on Varifolds)
Let $1\le m < n$ be integers, $U \subset \mathbb{R}^n$, and $V \in \textbf{RV}_m(U)$.  Let $Y$ be a separable normed space and $f:U \rightarrow Y$ be a $\norm{V}$-measurable function with compact support in $U$, we shall sometimes use the following notation
    \begin{align*}
        \norm{V}_{(p)}(f) := \left(\int |f(x)|^p_{Y} d\norm{V}(x)\right)^{1/p}\\
        L^p(\norm{V}, Y) := \{f: U \rightarrow Y \text{  s.t.  } \norm{V}_{(p)}(f)<\infty \}.
    \end{align*}
When $Y$ is clear from context, we shall often write $L^p(\norm{V}) = L^p(\norm{V}, Y)$ for concision.  When $B_r^n(x) \subset \subset U$, we shall identify $\norm{V}_{(p)}(f \chi_{B_r^n(x)}) = \norm{V \res B_r^n(x)}_{(p)}(f)$.

Moreover, if $f \in C^1(U)$ and $x \in \text{spt}\norm{V}$, we use the notation 
\begin{align*}
    \nabla^M f (x):= (S(x)_{\#} \circ \nabla f)(x).
\end{align*}
Additionally, we use the notation 
\begin{align*}
    \nabla^\perp f (x):= \nabla f(x) - \nabla^{M} f(x).
\end{align*}
\end{definition}

\begin{definition}(Almgren's $Q$-valued functions) 
For any $Q \in \mathbb{N}$ we let $Q_Q(\mathbb{R}^n)$ be the collection of unordered $Q$-tuples of points in $\mathbb{R}^n$.  For $A \subset \mathbb{R}^m$, a relation $f:A \rightarrow Q_Q(\mathbb{R}^n)$ will be called a \textit{$Q$-valued function}.  For $a \in A$, we shall thing of $f(a)= \{b_i\}_{i=1}^Q$ as a sum of $Q$ Dirac masses, using the notation 
\begin{align*}
    f(a) = \sum_{i=1}^Q \llbracket b_i \rrbracket.
\end{align*}
Additionally, we shall write
\begin{align*}
    \text{graph}f = \{(a, b): a \in A, b \in \text{spt} f(a) \}.
\end{align*}

Given a $T \in G(n,m)$ and $s \in Q_Q(\mathbb{R}^{n-m})$ we define the corresponding \textit{$Q$-valued plane parallel to $T$}, denoted by $P_s$, by
\begin{align*}
    P_s := (\Theta^0(\norm{s}, \cdot ) \circ T_{\#}^{\perp})\mathcal{H}^m.
\end{align*}
Clearly, each $Q$-valued plane parallel to $T$ defines an $s \in Q_Q(\mathbb{R}^{n-m})$.
\end{definition}

\section{Preliminaries}\label{s:prelim results}

In this subsection, we record a few results that will be essential for the proof of Theorem \ref{t:BMO main theorem rect} and its corollaries.  The results of Section \ref{s:lipschitz approx} rely upon compactness and contradiction arguments.  Therefore, below we generalize Allard's compactness for integral varifolds (\cite{Allard72} Theorem 6.4) to rectifiable varifolds with almost-integral densities.  In fact, we shall see that the preliminary lemmata necessary for the Lipschitz approximation hold for a wider class of rectifiable varifolds (see Definition \ref{d:appropriate}).  Next, we state some fundamental results for rectifiable varifolds such as the monotonicity formula (Lemma \ref{l:monotonicity formula}) and the Caccioppoli estimate (Lemma \ref{l:Simon 22.2}).  We end this Section with some elementary results about harmonic functions which will be necessary for the tilt-excess decay lemma in Section \ref{s:tilt-excess decay}.

\subsection{Compactness}
\begin{definition}\label{d:appropriate}
For integers $1 \le m \le n$ and an open set $U \subset \mathbb{R}^n$, we borrow the terminology of \cite{Chou22_varifold_decomp} and say that a set $P \subset \textbf{RV}_m(U)$ is \textit{appropriate} if 
\begin{enumerate}
    \item For all $V, W \in P$, $V+W \in P$.
    \item For all $V \in P$, $\Theta^m(\norm{V}, x) \ge 1$ for $\norm{V}$-a.e. $x \in \text{spt}\norm{V}$.
    \item $P$ is closed with respect to the strong topology.
\end{enumerate}

In analogy, we shall say that a set $\mathbb{A} \subset [1, \infty)$ is \textit{appropriate} if and only if $\mathbb{A}$ is closed and closed under addition.  Following Definition \ref{d:almost integral densities}, we shall have the appropriate set $\mathbb{A} = \cup_{n \in \mathbb{N}}[n, n(1+\gamma_{\Theta})]$ in mind for applications.
\end{definition}

\begin{theorem}\label{t:rect varifolds compactness}\emph{(Compactness for Rectifiable Varifolds with Appropriate Densities)}[cf. Simons, Chapter 8, Theorem 5.8]
Let $\mathbb{A}\subset \mathbb{R}$ be an appropriate set, and let $\{V_i\}_{i=1}^{\infty} \subset \textbf{\emph{RV}}_m(U)$ be a sequence which satisfies
\begin{align*}
    \sup_i\{\norm{V_i}(W) + \norm{\delta V_i}(W)\}< \infty.
\qquad \text{for all  } W\subset \subset U\end{align*}
Suppose also that there exists a sequence of sets $A_i \subset U$ such that $\Theta^m(\norm{V_i}, x) \in \mathbb{A}$ for $\norm{V_i}$-a.e. $x \in \emph{spt}\norm{V_i} \setminus A_i$ for each $i \in \mathbb{N}$ and $\norm{V_i}(A_i) \rightarrow 0$.

Then, there is a subsequence $\{V_j \}_{j}$ and a $V \in \textbf{\emph{RV}}_m(U)$ such that 
\begin{enumerate}
    \item $V_i \rightarrow V$ in the sense of varifolds, i.e., in the weak topology.
    \item $\Theta^m(\norm{V}, x) \in \mathbb{A}$ for $\norm{V}$-a.e. $x \in \emph{spt}\norm{V}$.
    \item $\norm{\delta V}(W) \le \liminf_{i}\norm{\delta V_i}(W)$ for all $W\subset \subset U$.
\end{enumerate}.
\end{theorem}

\begin{remark}
    This compactness result implies that if $\mathbb{A}$ is an appropriate set and $P_\mathbb{A}(U) \subset \textbf{RV}_m(U)$ is the set of all rectifiable varifolds $V$ with $\Theta^m(\norm{V}, x) \in \mathbb{A}$ for $\norm{V}$-a.e. $x \in \text{spt}\norm{V}$ then $P_\mathbb{A}(U)$ is also appropriate in the sense of Definition \ref{d:appropriate}.
\end{remark}

Consulting \cite{Simon83_book} Chapter 8 Theorem 5.8, we see that the only new claim in Theorem \ref{t:rect varifolds compactness} is that the limit manifold $V$ satisfies $\Theta^m(\norm{V}, x) \in \mathbb{A}$ for $\norm{V}$-a.e. $x \in \text{spt}\norm{V}$. To see this claim, we need the following lemma.

\begin{lemma}\label{l:convergence to a flat plane}
    Let $\mathbb{A} \subset \mathbb{R}$ be an appropriate set.  Suppose that $W_i \in \textbf{\emph{RV}}_{m}(B_1^n(0))$ satisfy
    \begin{enumerate}
        \item $\Theta^m(\norm{W_i}, x) \in \mathbb{A}$ for $\norm{W_i}$-a.e. $x \in \emph{spt}\norm{W_i}$ for each $i \in \mathbb{N}$.
        \item There is a $V \in \textbf{\emph{RV}}_m(B_1^n(0))$ such that $V = \theta_0\mathcal{H}^m \res T$ for some $T \in G(n,m)$ and $\theta_0 \in [1, \infty)$, and 
        \begin{align*}
            W_i \rightarrow V.
        \end{align*}
        \item $\norm{\delta W_i}(B_1^n(0)) \rightarrow 0$.
    \end{enumerate}
    Then, $\theta_0 \in \mathbb{A}$.
\end{lemma}

\begin{proof}
Let $T_{\#}$ be orthogonal projection onto $T$.  Then for any $0<\epsilon$, we note that the pushforward varifold $T_{\#}\left(W_i \res B_{\epsilon}(P) \times G(n,m)\right)$ (see \cite{Simon83_book} Chapter 8, Section 2 for definition) satisfies
\begin{align*}
    T_{\#}\left(W_i \res B_{\epsilon}(T) \times G(n,m)\right) = \int J_{W_i}(T_{\#})(x) \Theta^m(\norm{W_i}, x)d\norm{W_i}(x),
\end{align*}    
where $J_{W_i}(T_{\#})(x)$ is the Jacobian of $T_{\#}$ restricted to $T_{x}W_i$.  That is, $J_{W_i}(T_{\#})(x) = (1+ \norm{T_xW_i - T}^2)^{\frac{-1}{2}}$.  Note that for $\{\vec{v}_k\}_{k = m+1}^{n}$ an orthonormal basis for $T^{\perp}$ the divergence theorem gives that for all $f \in C^{\infty}_c(B_1^n(0))$
\begin{align*}
    \int f\norm{T_xW_i - T}^2 d\norm{W_i} & = \sum_{k=m+1}^{n}\int f\norm{D \dist(\cdot, \vec{v}_k^{\perp})}^2d\norm{W_i}\\
    & \le \int \sum_{k=m+1}^{n} f\vec{v}_k \cdot \textbf{h}(W_i, x)d\norm{W_i}(x) + \int \dist(\cdot, \vec{v}_k^{\perp}) \nabla f \cdot \vec{v}_k d\norm{W_i}\\
    & \le |f|_{\infty}\norm{\delta W_i}(B_1^n(0)) + |\nabla f|_{\infty}(n-m)\norm{W_i}_{(1)}(\dist(\cdot, T)).
\end{align*}
Thus, the assumption that $\norm{\delta W_i}(B_1^n(0)) \rightarrow 0$ and $W_i \rightarrow V$ implies that $\norm{W_i}_{(1)}(\norm{T_xW_i - T}^2) \rightarrow 0$ as $i \rightarrow \infty$ and hence
\begin{align*}
    T_{\#}\left(W_i \res B_{\epsilon}(T) \times G(n,m)\right) \rightarrow \theta_0 V.
\end{align*}
Unwinding the definition of $T_{\#}\left(W_i \res B_{\epsilon}^n(T)\times G(n,m)\right)$ we see that this implies that for all $f \in C^\infty_{c}(B_1^n(0))$
\begin{align}\label{e:psi convergence to theta}
    \int_{T}f \psi_i d \mathcal{H}^m \rightarrow \int_T f \theta_0 d \mathcal{H}^m,
\end{align}
where $\psi_i(y) := \sum_{x \in T_{\#}^{-1}(y) \cap B_{\epsilon}^n(T)} \Theta^m(\norm{W_i}, x)$.  Note that since by assumption $\Theta^m(\norm{W_i}, x) \in \mathbb{A}$ for $\norm{W_i}$-a.e. $x \in \text{spt}\norm{W_i}$ and $\mathbb{A}$ is appropriate, we have $\psi_i(y) \in \mathbb{A}$ for $\norm{W_i}$-a.e. $x \in \text{spt}\norm{W_i}$.

Now, for any $0<\epsilon$ we define the set $A_i(\epsilon)\subset \text{spt}\norm{W_i}$ by 
\begin{align*}
    A_i(\epsilon):= \{x \in T: \psi_i(x) > \theta_0 + \epsilon \}.
\end{align*}
Note that by \eqref{e:psi convergence to theta} it must be that $\mathcal{H}^m(A_i(\epsilon)) \rightarrow 0$ as $i \rightarrow \infty$. Thus, the function $\min\{\theta_0 + \epsilon, \psi_i(x)\}$ converges in $L^1(B_1^n(0))$ to $\theta_0$. Since this holds for any $0< \epsilon$, if there did not exist a sequence of $a_i \in \mathbb{A}$ such that $a_i \le \theta_0$ and $a_i \rightarrow \theta_0^-$, \eqref{e:psi convergence to theta} would fail for non-negative $f$.  Thus, such a sequence $a_i$ must exist and $\theta_0 \in \mathbb{A}$ by the closedness of $\mathbb{A}$.
\end{proof}

\subsubsection{Proof of Theorem \ref{t:rect varifolds compactness}}
Follow the proof of \cite{Simon83_book} Chapter 8, ``Proof that $V$ is integer multiplicity if the $V_i$ are" after Theorem 5.8 verbatim, using Lemma \ref{l:convergence to a flat plane} to finish the last half of the argument. \qed

\subsection{Basic results for rectifiable varifolds.} 

\begin{lemma}\label{l:monotonicity formula}\emph{(Monotonicity Formula, \cite{Simon83_book}, (17.4) and p. 90, between (18.2) and (18.3))}
Let $U \subset \mathbb{R}^n$ be an open set, and let $V \in \textbf{\emph{RV}}_m(U)$ with generalized mean curvature vector $\boldsymbol{h}(V; x)$.  Let $B_R^n(x) \subset U$.  For any non-negative function $h \in C^1(U)$ and any $0<\sigma < \rho \le R$,
\begin{align*}
\frac{\int_{B_\sigma^n(x)}hd\norm{V}}{\omega_m \sigma^m} & \le \frac{\int_{B_{\rho}^n(x)}hd\norm{V}}{\omega_m \rho^m}\\
& \qquad + \frac{1}{\omega_m}\int_{\sigma}^{\rho}\tau^{-m}\int_{B_{\tau}^n(x)}\frac{|z-x|}{\tau}(|\nabla^M h(z)| + h |\boldsymbol{h}(V; z)|)d\norm{V}(z)d\tau.
\end{align*}
\end{lemma}

\begin{lemma}\label{l:Simon 22.2}\emph{(The Caccioppoli-type Estimate I, \cite{Simon83_book} Lemma 22.2)} Let $U \subset\mathbb{R}^n$ be an open set and let $V \in \textbf{\emph{RV}}_m(U)$ with generalized mean curvature $\boldsymbol{h}(V; \cdot) \in L^1_{loc}(\norm{V})$. Let $B_{\rho}^n(x) \subset U$.  Then for any $T \in G(n,m)$
\begin{align*}
    E(x, \rho/2, T) \le c \Big[\rho^{-m} \int_{B_{\rho}^n(x)}\left|\frac{\dist(z-x, T)}{\rho} \right|^2d\norm{V}(z) + \rho^{2-m}\int_{B_{\rho}^n(x)}|\boldsymbol{h}(V; z)|^2d\norm{V}(z)  \Big].
\end{align*}    
\end{lemma}

In the case $\boldsymbol{h}(V, \cdot) \not \in L^2(\norm{V})$ Lemma \ref{l:Simon 22.2} holds trivially, but is useless for obtaining tilt-excess estimates.  Thus, in the critical case with $m=1$ we need the following, related estimate.

\begin{lemma}\label{l:brakke 5.5 m=1}\emph{(The Caccioppoli-type Estimate II, cf. \cite{Brakke78} 5.5)}
    Let $V \in \emph{\textbf{RV}}_m(B_1^n(0))$, $T \in G(n,m)$, and $\phi \in C^1_0(\mathbb{R}^n, \mathbb{R}^+)$.  Then,
    \begin{align*}
        & \int |S(x)_{\#} - T_{\#}|^2 \phi^2(x)d\norm{V}(x) \le \\
        & \qquad 4m\left(\int \phi^2(x)|\boldsymbol{h}(V, x)|d\norm{V}(x)\right)^{1/3} \left(\int|T^{\perp}_{\#}(x)|^2 \phi^2(x) d\norm{V}(x)\right)^{1/3}\\
        & \qquad \qquad + 16 \int |T^{\perp}_{\#}(x)|^2 |D\phi(x)|^2 d\norm{V}(x).
    \end{align*}
\end{lemma}

\begin{proof} We note that if $V \in \textbf{IV}_m(\mathbb{R}^n)$, the lemma is proven in \cite{Brakke78} 5.5.  However, the proof relies entirely upon the manipulation of the integrand and goes through identically for rectifiable varifolds.
\end{proof}

\begin{theorem}\label{t:Allard Lemma 8.4}\emph{(\cite{Allard72} Lemma 8.4)}
Let $V \in \textbf{\emph{RV}}_m(B_1^n(0))$ which satisfies the hypotheses of the $(m, m, \tilde{\gamma})$-Allard regularity regime. For all $0<\delta<1$ there is an $0<\epsilon(\delta)$ such that if $0<\tilde{\gamma} \le \epsilon(\delta)$ then then there is a $T \in G(n,m)$ such that  
\begin{align*}
\dist_{\mathcal{H}}(\emph{spt}\norm{V} \cap B_{1-\delta}^n(0),  T \cap B_{1-\delta}^n(0)) \le \delta.
\end{align*}
\end{theorem}

\begin{remark}\label{r:reg regime implies tilt decay regime}
Together, Lemma \ref{l:Simon 22.2}, Lemma \ref{l:brakke 5.5 m=1}, and Theorem \ref{t:Allard Lemma 8.4} imply that for sufficently small $\tilde{\gamma}'$, varifolds which satisfy the $(m,m,\tilde{\gamma}')$-Allard regularity regime in $B_r^n(0)$ also satisfy the critical $(\delta_0, \tilde{\gamma})$-Allard tilt-excess decay regime in $B_{r/2}^n(0)$. Specifically, when $m \ge 2$ we may use Theorem \ref{t:Allard Lemma 8.4} and Lemma \ref{l:Simon 22.2}.  For $m=1$, we use Theorem \ref{t:Allard Lemma 8.4} and Lemma \ref{l:brakke 5.5 m=1} in $B_{r}^n(0)$ with $\phi = 1$ in $B_{r/2}^n(0)$, $\text{spt}\phi \subset B_r^n(0)$, and $|\nabla \phi| \le 4/r$.
\end{remark}

\subsection{Harmonic and almost-weakly harmonic functions}

\begin{lemma}\label{l:harmonic approximation}\emph{(\cite{Simon83_book} Lemma 21.3)}
For any $0<\epsilon$ there exists a $0< \delta(m, \epsilon)$ such that if $f \in W^{1, 2}(B^m_1(0))$ satisfies 
\begin{align*}
    \int_{B_1^m(0)}|\nabla f|^2 d\mathcal{L}^m \le 1\\ |\int_{B^m_1(0)}\nabla f \cdot \nabla \xi d\mathcal{L}^m| \le \delta \sup |\nabla \xi|
\end{align*}
for all $\xi \in C^\infty_c(B_1^m(0); \mathbb{R})$, then there exists a harmonic function $u: B_1^m(0) \rightarrow \mathbb{R}$ such that 
\begin{align*}
    \int_{B^m_1(0)}|\nabla u|^2 d\mathscr{L}^m \le 1  \qquad \int_{B_1^m(0)}|f-u|^2 d\mathscr{L}^m \le \epsilon.
\end{align*}
\end{lemma}

\begin{lemma}\label{l: harmonic D2 to D1 estimate}\emph{(cf. \cite{Simon83_book}, (21.3))}
    Let $B^m_r(0) \subset B_{R/4}^m(0) \subset \mathbb{R}^m$. For any harmonic function $u: B^m_{R}(0) \rightarrow \mathbb{R}$, if we let $\ell: B^m_{R}(0) \rightarrow \mathbb{R}$ be 
    \begin{align*}
         \ell(y) = u(0) + Du(0)\cdot y,
    \end{align*}
    then there is a dimensional constant $C(m)$ such that
    \begin{align*}
        \sup_{B^m_{\rho}(0)}|u-\ell| \le C(m)\rho \left(\frac{\rho}{R}\right)(\frac{R}{2}-\rho)^{-m/2} \norm{Du}_{L^2(B_{R}^m(0))}.
    \end{align*}
    Note that when $m=1$ we may take $C(m) =0$.
\end{lemma}

\begin{proof}
By Taylor's expansion, for any $y \in B_{\rho}^m(0)$ $|u(y)-\ell(y)| \le |y|^2\max_{B_{\rho}^m(0)}|D^2u|$.  For the $x_{\max} \in \overline{B_\rho^m(0)}$ such that 
    $|D^2u(x_{\max})| = \max_{B_{\rho}^m(0)}|D^2u|$, we note that by Jensen's inequality
    \begin{align*}
        D^2u(x_{\max}) & = \dashint_{B_{\frac{R}{2}-\rho}^m(x_{\max})}Du d\mathcal{L}^m\\
        & \le (\frac{R}{2}-\rho)^{-m/2}\norm{D^2u}_{L^2(B_{\frac{R}{2}-\rho}^m(x_{\max}))}\\
        & \le (\frac{R}{2}-\rho)^{-m/2}\norm{D^2u}_{L^2(B_{R/2}^m(0))}. 
    \end{align*}
    The proof is now finished by recalling the classical observation that by compactness there is a dimensional constant $C(m)$ such that $\norm{D^2u}_{L^2(B^m_{R/2}(0))}\le C(m)R^{-1}\norm{Du}_{L^2(B_{R}^m(0))}.$
\end{proof}

\section{Lipschitz Approximation}\label{s:lipschitz approx}

This Section is the technical heart of the paper.  In it we adapt the $Q$-valued Lipschitz approximation for integral varifolds from \cite{Menne10_Sobolev} Lemma 3.15 to the case of rectifiable varifolds with almost integral densities. More specifically, the Lipschitz approximation (Lemma \ref{l:menne 2.18}) relies upon a number of preliminary lemmata (Lemma \ref{l:menne 2.2} - Lemma \ref{l:Menne 2.13}) which are minor generalizations of corresponding lemmata in \cite{Menne10_Sobolev} Section 3.  It is interesting to observe that these lemmata actually hold for a rectifiable varifolds with densities in some appropriate set $\mathbb{A}$ (see Definition \ref{d:appropriate}).  In the interest of generality, we prove them in this case, but only use the appropriate set $\cup_{n \in \mathbb{N}}[n, n(1+\gamma_{\Theta})]$ in the proof of Lemma \ref{l:menne 2.18}.

Because the differences between Lemma \ref{l:menne 2.2} - Lemma \ref{l:Menne 2.13} and their corresponding integral varifold results come entirely from considering gaps in $\mathbb{A}$ rather than gaps in $\mathbb{N}$, we shall give full statements and indicate the changes required to the corresponding proofs in \cite{Menne10_Sobolev}.  However, the proof of Lemma \ref{l:menne 2.18} itself is more complicated.  For the reader's convenience we provide full details and keep careful track of constants.

\subsection{Results preliminary to the Lipschitz approximation}

\subsubsection{Density estimates} In the critical case lower density bounds are closely related to the isoperimetric inequality.  This observation goes back to \cite{Federer69} 4.1.25 for integral currents and was first used in the context of varifolds by Allard in \cite{Allard72} 8.3.

\begin{definition}
Let $\gamma_m<\infty$ be defined to be the smallest number such that the following property holds.

\textit{If $1 \le m \le n\in \mathbb{N}$ and $V \in \textbf{RV}_m(\mathbb{R}^n)$ with $\norm{V}(\mathbb{R}^{n})<\infty$ and $\norm{\delta V}(\mathbb{R}^n)<\infty$, then }
    \begin{align*}
        \norm{V}(\{x \in \mathbb{R}^n: \Theta^m(\norm{V}, x)\ge 1\}) \le \gamma_m \norm{V}(\mathbb{R}^n)^{1/n}\norm{\delta V}(\mathbb{R}^n).
    \end{align*}
By the isoperimetric inequality, $0<\gamma_m<\infty$.  See \cite{Menne09_isoperimetric} Section 2 for more details and discussion on $\gamma_m$.
\end{definition}

\begin{lemma}\label{l:menne sobolev 2.4}\emph{(Lower density bound at all scales, \cite{Menne10_Sobolev} Lemma 3.4)}
Let $n \in \mathbb{N}$, $0<r<\infty$, and $V \in \textbf{\emph{RV}}(B_r^n(0))$ with  $\Theta^m(\norm{V}, x ) \ge 1$ for $\norm{V}$-a.e. $x \in \emph{spt}\norm{V}$. For any $0<\kappa_1<1$ there exists a $0<\gamma_0= \gamma_0(m, n, \kappa_1)$ such that the following holds.

If $V$ satisfies 
\begin{align*}
    \norm{\delta V}(\overline{B_\rho^n(0)}) & \le (2\gamma_m)^{-1}\norm{V}(\overline{B_{\rho}^n(0)})^{1-1/m} \quad \text{for all } 0<\rho<r\\
    \norm{\delta V}(B_{1}^n(0)) & \le \epsilon \norm{V}(B_r^n(0))^{1-1/m},
\end{align*}
and $\epsilon \le \gamma_0(m, n, \kappa_1)$ then for all $0<\rho<r$
\begin{align*}
    1-\kappa_1 \le \frac{\norm{V}(B_{\rho}^n(0))}{\omega_m \rho^m}.
\end{align*}
\end{lemma}

\subsubsection{The Limit Case}

Lipschitz approximation theorems may be seen as a kind of continuity result relative to the following lemma.  That is, the key steps are based upon limit-compactness arguments where the limiting case satisfies the hypotheses of Lemma \ref{l:flat stationary RV}.  We state the integral varifold version first for use in the proof.

\begin{lemma}\label{l:flat stationary IV}\emph{(\cite{Menne10_Sobolev} Lemma 3.1)}
Let $1 \le m\le n$ be integers, and let $B_r^n(a) \subset \mathbb{R}^n$.  Suppose that $V \in \textbf{\emph{IV}}_m(B_r^n(a))$ is stationary and there is a $T \in G(n, m)$ such that $S(x) = T$ for $\norm{V}$-a.e. $x \in \emph{spt}\norm{V}$, then the following statements hold.
\begin{enumerate}
    \item $T^{\perp}(\emph{spt}\norm{V})$ is discrete and closed in $T^{\perp}(B_r(a))$.
    \item For every $x \in \emph{spt}\norm{V}$ and $y \in B_r^n(a)$, if $y-x \in T$ then $\Theta^m(\norm{V}, y) = \Theta^m(\norm{V}, x) \in \mathbb{N}$.
    \item For any $x \in B_r^n(a)$ if we let $S_x = \{y \in B_r^n(a): y - x \in T\}$, then $V \res S_x = \Theta^m(\norm{V}, x)\mathcal{H}^m \res S_x$.
\end{enumerate}
\end{lemma}

Therefore, extending these results to almost-integral varifolds is essential to establishing an analogous result. The argument comes from the monotonicity formula and in the proof of \cite{Simon83_book} Theorem 19.5 and Theorem 20.2.

\begin{lemma}\label{l:flat stationary RV}\emph{(The Limit Case)}
    Let $1 \le m\le n$ be integers, and let $B_r^n(a) \subset \mathbb{R}^n$.  Suppose that $V \in \textbf{\emph{RV}}_m(B_r^n(a))$ with $\Theta^m(\norm{V}, x) \ge 1$ for $\norm{V}$-a.e. $x \in \emph{spt}\norm{V}$. If $V$ is stationary and there is a $T \in G(n, m)$ such that $S(x) = T$ for $\norm{V}$-a.e. $x \in \emph{spt}\norm{V}$, then the following statements hold.
\begin{enumerate}
    \item $T^{\perp}(\emph{spt}\norm{V})$ is discrete and closed in $T^{\perp}(B_r(a))$.
    \item For every $x \in \emph{spt}\norm{V}$ and $y \in B_r^n(a)$, if $y-x \in T$ then $\Theta^m(\norm{V}, y) = \Theta^m(\norm{V}, x) \in \mathbb{N}$.
    \item For any $x \in B_r^n(a)$ if we let $S_x = \{y \in B_r^n(a): y - x \in T\}$, then $V \res S_x = \Theta^m(\norm{V}, x)\mathcal{H}^m \res S_x$.
\end{enumerate}
\end{lemma}

\begin{proof}
    Let $K \subset \subset B_r^n(a)$ be a compact subset, 
    \begin{align*}
        0<\delta < \dist(K, \mathbb{R}^n \setminus B_r^n(a)), \qquad \norm{V}(B_{\delta}^n(K)) \le M \omega_m r^m,\\
        0<\rho <\delta \text{ such that }(1-\rho/\delta)^{-m} \le 10,
    \end{align*}
    and let $x \in \text{spt}\norm{V} \cap K$. For any $0<\epsilon<\delta$ suppose that we can find $\{x_i\}_{i=1}^{N}$ such that
    \begin{align*}
        x_i \in B_{\rho}^n(x) \cap \text{spt}\norm{V}, \quad |T^{\perp}_{\#}(x_i - x_j)| \ge \epsilon.
    \end{align*}
In this case, we choose an appropriate function $h$ for use in Lemma \ref{l:monotonicity formula}. Let $g:\mathbb{R}^+ \rightarrow \mathbb{R}^+$ a smooth function such that 
\begin{align*}
    \begin{cases}
    g(t) \equiv 1 & 0 \le t < \epsilon/2\\
    g(t) \equiv 0 & t \ge \epsilon
\end{cases}
\end{align*}
and $|g'(t)| \le 4\epsilon^{-1}$.  For $\xi \in \{x_i\}_{i=1}^{M'}$, let $h_{\xi}(x):= g(|T_{\#}^{\perp}(x-\xi)|)$.  Note that for all $j = 1, ..., m, ..., n$,
\begin{align*}
    |\nabla^M_j T_{\#}^{\perp}(x-\xi)| \le |S(x)_{\#} \circ T_{\#}^{\perp}| & =|(S(x)_{\#} - T_{\#})\circ T_{\#}^{\perp}|\\
    & \le |S(x)_{\#} - T_{\#}|.
\end{align*}

Now, let $0<\sigma \le \epsilon$ and $\rho = (1-\rho/\delta)\delta$ and apply \ref{l:monotonicity formula} to $\xi \in \{x_i\}^{M'}_{i=1}$ for $h_{\xi}$ to obtain
\begin{align}\nonumber
    \frac{\norm{V}(B_{\sigma}^n(\xi))}{\omega_m \sigma^m} & \le \frac{\norm{V}(B_{\rho}^n(\xi) \cap \{z: |T_{\#}^{\perp}(z-\xi)| \le \epsilon \})}{\omega_m \rho^m}.
\end{align}
Summing these estimates up across $i=1,..., N$, we obtain from $\norm{\delta V}=0$ and $\Theta^m(\norm{V}, \cdot) \ge 1$ that
\begin{align*}
    N & \le  \sum_{i=1}^{N}\frac{\norm{V}(B_{\sigma}^n(x_i))}{\omega_m\sigma^m} \le \frac{\sum_{i=1}^{N} \norm{V}(B_{\rho}^n(x_i) \cap \{z: |T_{\#}^{\perp}(z-x_i)| \le \epsilon \})}{\omega_m \rho^m}\\
    & \le \frac{M}{(1-\rho/\delta)^m} \le 10M. 
\end{align*}
Since this holds for all sufficiently small $0<\epsilon$, we infer that $\mathcal{H}^0(T^{\perp}(\text{spt}\norm{V} \cap B_{\rho}^n(x))) \le 10M$.  Repeating this for each $x \in \text{spt}\norm{V}\cap K$ and choosing a finite subcover shows that $T^{\perp}_{\#}(W \cap K)$ is finite for all $K \subset \subset B^n_r(a)$. 

Therefore, taking a compact exhaustion $\{K_i\}_i$ of $B_r^n(a)$ we may apply the constancy theorem (see \cite{Duggan_varifold86} Corollary 4.5) to $V \res{K_i}$ to obtain $\Theta^m(\norm{V}, \cdot)$ is constant on each of the finite components of the support of $\norm{V} \cap K_i$. Taking an appropriate sequence of constants $\{c_i\}_i$ such that if $\{W_i\}$ is an enumeration of the connected components of $V$, $c_iW_i \in \textbf{IV}_m(B_r^n(a) \cap K_i)$, we may apply Lemma \ref{l:flat stationary IV}.  Letting $i \rightarrow \infty$ gives the desired conclusion. 
\end{proof}

\subsubsection{Lemmata preliminary to the Lipschitz approximation} 

For the purposes of the next several lemmata, we fix the following definitions.

\begin{definition}\label{d:C, Q, delta 12, gamma def}
    Let $1 \le m < n$ be integers. Fix $\mathbb{A} \subset \mathbb{R}$ an appropriate set.  For $Q \in \mathbb{A}$ we shall say that $0<\delta_1<\delta_2<1/2$ are \textit{lower admissible} for $Q$ if $[Q-2\delta_2, Q -\delta_1] \cap \mathbb{A} = \emptyset$. We shall say that $0<\delta_3<\delta_4<1$ are \textit{upper admissible} for $Q$ if $[Q+\delta_3, Q +\delta_4] \cap \mathbb{A} = \emptyset$.
\end{definition}

\begin{lemma}\label{l:menne 2.2}\emph{(Cone density with offset, cf. \cite{Menne10_Sobolev} Lemma 3.2)}
    Let $m, n, Q, \mathbb{A}$ be as in Defintion \ref{d:C, Q, delta 12, gamma def}. Let $0\le s<1$, $0<\delta<1$, and $0 \le M <\infty$. Then there is an $0<\epsilon$ such that the following holds.

Let $T \in G(n,m)$, $0 \le d \le \infty$, $0<t<\infty$, and $\zeta \in \mathbb{R}^n$ satisfy 
    \begin{align*}
        \max\{d, r \} \le M, \quad \zeta \in \overline{B_{d}^n(0)} \cap T, \quad d + t \le r.
    \end{align*}
    Suppose that $a \in \mathbb{R}^n$, $0<r<\infty$, and $V \in \textbf{\emph{RV}}_m(B_r^n(a))$ such that 
    \begin{enumerate}
        \item $V$ has locally bounded first variation.
        \item $\Theta^m(\norm{V}, x) \in \mathbb{A}$ for $\norm{V}$-a.e. $x \in \emph{spt}\norm{V}$.
        \item (Mean curvature control) $\norm{\delta V}(B_r(a)) \le \epsilon \norm{V}(B_r^n(a))^{1-1/n}$.
        \item (Tilt control) $\int_{B_r^n(a) \times G(n,m)}|S_{\#} - T_{\#}|dv(z, S) \le \epsilon \norm{V}(B_r^n(a))$.
        \item (Mass control) $\norm{V}(B_r^n(a)) \le M\omega_m r^m$, and for all $0<\rho<r$
        \begin{align*}
            \norm{V}(B_\rho^n(a)) \ge \delta \omega_m \rho^m.
        \end{align*}
Then, $\norm{V}(\{x \in B_t^n(a+\zeta): |T_{\#}(x-a)|> s|x-a|\}) \ge (1-\delta)\omega_m t^m$.
    \end{enumerate}
\end{lemma}

\begin{proof}
    The proof follows \cite{Menne10_Sobolev} Lemma 3.2, with Theorem \ref{t:rect varifolds compactness} in place of \cite{Allard72} 6.4.
\end{proof}

\begin{lemma}\label{l:Menne 2.6}\emph{(Multilayer monotonicity for rectifiable varifolds with appropriate density, cf. \cite{Menne10_Sobolev} Lemma 3.6)}
Let $m,n, \mathbb{A}, Q$ be as in Definition \ref{d:C, Q, delta 12, gamma def} and let $0<\delta_{\ref{l:Menne 2.6}, 1}< \delta_{\ref{l:Menne 2.6}, 2}<1$ be lower admissible for $Q$.  Let $0\le s<1$.  Then there is a positive, finite number $\epsilon$ with the following property.

Suppose $X \subset \mathbb{R}^n$, $0<r<\infty$, and $U = \cup_{x \in X}B_r^n(x)$.  Let $V \in \textbf{\emph{RV}}_m(\cup_{x \in X}B_r^n(x))$ be a varifold with locally bounded variation which satisfies the following properties.
\begin{enumerate}
    \item (Geometry) There exists a $T \in G(n,m)$ such that for all $x, y \in X$. 
\begin{align*}
    |T_{\#}(x-y)| \le s |x-y|,
\end{align*}
    \item (Density) For $\norm{V}$-a.e. $x \in \emph{spt}\norm{V}$, $\Theta^m(\norm{V}, x) \in \mathbb{A}$.  Moreover, assume that
    \begin{align*}
    \sum_{x \in X} \Theta^m_*(\norm{V}, x) \ge Q-2\delta_{\ref{l:Menne 2.6}, 2}.
\end{align*}
    \item (Mean Curvature and Tilt)  For all $x \in X \cap \emph{spt}\norm{V}$ and all $0<\rho \le r$
\begin{align*}
    \norm{\delta V}(\overline{B_{\rho}^n(x)})\le \epsilon \mu(\overline{B_{\rho}^n(x)})^{1-\frac{1}{n}}, \quad \int_{\overline{B_{\rho}^n(x)}}|S_{\#} - T_{\#}|dV(z, S) \le \epsilon \mu(\overline{B_{\rho}^n(x)}).
\end{align*}
\end{enumerate}

Then, for all $0<\rho \le r$
\begin{align*}
    \mu(\cup_{x \in X}B_\rho^n(x)) \ge (Q - \delta_{\ref{l:Menne 2.6}, 1})\omega_m \rho^m.
\end{align*}
\end{lemma}

\begin{proof}
The proof follows the proof of \cite{Menne10_Sobolev} Lemma 3.6 with Theorem \ref{t:rect varifolds compactness} in place of \cite{Allard72} 6.4. 

\end{proof}

\begin{lemma}\label{l:Menne 2.8}\emph{(cf. \cite{Menne10_Sobolev} Lemma 3.7)}
Suppose that $\mathbb{A}\subset \mathbb{R}$ is an appropriate set.  Let $M \not \in \mathbb{A}$, $0< \lambda_1 < \lambda_2<1$, and $T \in G(n, m)$.  Suppose that $F_\mathbb{A}$ is the collection of all $V \in \textbf{\emph{RV}}_m(B_1^n(0))$ such that 
\begin{enumerate}
    \item $S(x) = T$ for $\norm{V}$-a.e. $x \in \emph{spt}\norm{V}$.
    \item $\norm{\delta V} = 0$.
    \item $\Theta^m(\norm{V}, x) \in \mathbb{A}$ for $\norm{V}$-a.e. $x \in \emph{spt}\norm{V}$.
    \item $\norm{V}(B_1^n(0)) \le M \omega_m$.
\end{enumerate}
If we let $N = \sup_{V \in F_\mathbb{A}}\{\frac{\norm{V}(B_r^n(0))}{\omega_m r^m}: r \in [\lambda_1, \lambda_2]\}$, then we claim that $N<M$ and there exists $(V, r) \in F_\mathbb{A} \times [\lambda_1, \lambda_2]$ such that $N = \frac{\norm{V}(B_r^n(0))}{\omega_m r^m}$. 
\end{lemma}

\begin{proof}
The proof follows \cite{Menne10_Sobolev} Lemma 3.7 with $\mathbb{A}$ in place of $\mathbb{N}$ and Theorem \ref{t:rect varifolds compactness} in place of \cite{Allard72} 6.4.
\end{proof}

\begin{lemma}\label{l:Menne 2.10}\emph{(Quasi monotonicity, cf. \cite{Menne10_Sobolev} Lemma 3.9)} 
Suppose $\mathbb{A}$ is an appropriate set and $0<M<\infty$ is such that $M \not \in \mathbb{A}$.  Then, there exists a number $0<\epsilon< \infty$ such that the following holds.

If $a \in \mathbb{R}^n$, $0<r<\infty$, and $V \in \textbf{\emph{RV}}_m(B_r^n(a))$ such that 
\begin{enumerate}
    \item $\Theta^m(\norm{V}, x) \in \mathbb{A}$ for $\norm{V}$-a.e. $x \in \emph{spt}\norm{V}$.
    \item $\norm{V}(B_r^n(0)) \le M \omega_mr^m$.
    \item There exists a $T \in G(n,m)$ such that for all $0<\rho<r$ 
    \begin{align*}
        \norm{\delta V}(\overline{B_{\rho}^n(a)}) & \le \epsilon \norm{V}(\overline{B_{\rho}^n(a)})^{1-\frac{1}{n}},\\
        \int_{\overline{B_{\rho}^n(a)}} |S_{\#} - T_{\#}|dV(x, S) & \le \epsilon \norm{V}(\overline{B_{\rho}^n(a)}).
    \end{align*}
    Then $\norm{V}(\overline{B_{\rho}^n(a)}) \le M \omega_m \rho^m$ for all $0<\rho\le \lambda r$.
\end{enumerate}
\end{lemma}

\begin{proof}
    The proof follows \cite{Menne10_Sobolev} Lemma 3.9 verbatim with Lemma \ref{l:Menne 2.8} and Theorem \ref{t:rect varifolds compactness} in place of \cite{Menne10_Sobolev} Lemma 3.7 and \cite{Allard72} 6.4, respectively.
\end{proof}

\begin{lemma}\label{l:Menne 2.12}\emph{(Multilayer monotonicity with variable offset, cf. \cite{Menne10_Sobolev} Lemma 3.11)}
    Let $m,n, \mathbb{A}, Q$ be as in Definition \ref{d:C, Q, delta 12, gamma def} and let $0<\delta_{\ref{l:Menne 2.12}, 1}< \delta_{\ref{l:Menne 2.12}, 2}<1$ be lower admissible for $Q$. Let $0\le M<\infty$ and $0\le s<1$.  Then, there exists an $0<\epsilon$ such that the following holds.

    Let $X \subset \mathbb{R}^{n}$ and $T \in G(n,m)$ satisfy that for all $x, y \in X$
    \begin{align*}
        |T_{\#}(x-y)| \le s |x-y|.
    \end{align*}
    Let $0 < r< \infty$ and $0\le d, t<\infty$ such that 
    \begin{align*}
        d+t \le r, \qquad d \le Mt.
    \end{align*}
    Let $f: X \rightarrow \mathbb{R}^n$ by such that for all $x, y \in X$ 
    \begin{align*}
        |T_{\#}(f(x)-f(y))| \le s |f(x) - f(y)|\\
        f(x)-x \in B_{d}^n(0) \cap T.
    \end{align*}

    Let $V \in \textbf{\emph{RV}}_m(\cup_{x \in X}B_r^n(x))$ with locally bounded first variation such that $\Theta^m(\norm{V}, x) \in \mathbb{A}$ for $\norm{V}$-a.e. $x \in \emph{spt}\norm{V}$. Moreover, assume that
    \begin{align*}
        \sum_{x \in X} \Theta^m_{*}(\norm{V}, x) & \ge Q-\delta_{\ref{l:Menne 2.12}, 1}\\
        \norm{V}(B_r^n(x)) & \le M \omega_m r^m \quad \text{ for all } x \in X \cap \emph{spt}\norm{V},
    \end{align*}
    and for all $0<\rho<r$ and for all $x \in X \cap \emph{spt}\norm{V}$
    \begin{align*}
        \norm{\delta V}(\overline{B_\rho^n(x)}) \le \epsilon \norm{V}(\overline{B_\rho^n(x)})^{1-\frac{1}{n}}\\
        \int_{\overline{B_{\rho}^n(x)} \times G(n, m)}|S_{\#} - T_{\#}|dV(z, S) \le \epsilon \norm{V}(\overline{B_{\rho}^n(x)}).
    \end{align*}
    Then, $\norm{V}(\cup_{x \in X}\{y\in B_t^n(f(x)): |T_{\#}(x-y)|> s|x-y|\}) \ge (Q-\delta_{\ref{l:Menne 2.12}, 2})\omega_mt^m$.
\end{lemma}

\begin{proof}
    The proof follows \cite{Menne10_Sobolev} Lemma 3.11 with Lemma \ref{l:Menne 2.6}, Lemma \ref{l:Menne 2.10}, and Lemma \ref{t:rect varifolds compactness} in place of \cite{Menne10_Sobolev} Lemma 3.6, Lemma 3.9, and \cite{Allard72} 6.4, respectively.
\end{proof}

\begin{lemma}\label{l:Menne 2.13}\emph{(cf. \cite{Menne10_Sobolev} Lemma 3.12)}
    Let $m, n, \mathbb{A}, Q$ be as in Definition \ref{d:C, Q, delta 12, gamma def}.  Let $0<\delta_{ \ref{l:Menne 2.13}, 3}< \delta_{\ref{l:Menne 2.13}, 4} < 1$ be upper admissible for $Q$. We let $\mathbb{A}_Q\subset \mathbb{A}$ be the connected component of $\mathbb{A}$ which contains $Q$, and let $Q_{\max}:= \max \{x \in \mathbb{A}_Q \}$.  Let $0<M<\infty$, $0\le s, s_0<1$, $0<\delta_{\ref{l:Menne 2.13},1}<1$, and $0<\lambda< \infty$ be the unique number which satisfies
    \begin{align*}
        (1-\lambda^2)^{m/2} = (\delta_{\ref{l:Menne 2.13}, 4}- Q_{\max}) + \left(\frac{s_0^2}{1-s_0^2} \right)^{m/2}\lambda^m.
    \end{align*}
    Then there exists an $0<\epsilon$ such that the following holds. 

    Let $X \subset \mathbb{R}^n$, $T \in G(n,m)$, and suppose $\mathcal{H}^0(T(X))=1$.  Let $0<r<\infty$, $0 \le d \le \infty$, and $0<t<\infty$ be such that $d \le Mt$ and $d + t \le r$. Suppose that $\zeta \in \overline{B_{d}^n(0)} \cap T$.  

    Let $V \in \textbf{\emph{RV}}_m(\cup_{x \in X}B_r^n(x))$ with locally bounded first variation and $\Theta^m(\norm{V}, z) \in \mathbb{A}$ for $\norm{V}$-a.e. $z \in \text{spt}\norm{V}$ such that for all $x \in X$ 
    \begin{align*}
        \norm{V}(B_r^n(x)) \le M \omega_m r^m, \qquad \sum_{x \in X}\Theta^m(\norm{V}, x) \in \mathbb{A}_Q.
    \end{align*}
Moreover, assume that for all $x \in X$ and all $0<\rho \le r$
\begin{align*}
    \norm{\delta V}(\overline{B_\rho^n(x)}) \le \epsilon \norm{V}(\overline{B_\rho^n(x)})^{1-1/m}\\
    \int_{\overline{B_\rho^n(x)} \times G(n, m)}|S_{\#} - T_{\#}|^2dV(z, S) \le \epsilon \norm{V}(\overline{B_\rho^n(x)}).
\end{align*}
Finally, assume that
\begin{align*}
    \norm{V}(\cup_{x \in X}\{y \in B_t^n(x+\zeta): |T_{\#}(x-y)| > |x-y| \}) \le (Q+\delta_{\ref{l:Menne 2.13}, 4})\omega_n t^n.
\end{align*}
Then, the following statements hold.
\begin{enumerate}
    \item If $0< \tau \le \lambda t$ then $\norm{V}(\cup_{x \in X}\overline{B_{\tau}^n(x)}) \le (Q + \delta_{\ref{l:Menne 2.13}, 3})\omega_m\tau^m$.
    \item If $\xi \in \cup_{x \in X}\overline{B_{\lambda t/2}^n(x)}$ and for all $0<\rho < \delta_{\ref{l:Menne 2.13}, 1}\dist(\xi, X)$
    \begin{align*}
        \norm{V}(\overline{B_{\rho}^n(\xi)}) \ge \delta_{\ref{l:Menne 2.13}, 1} \omega_m \rho^m,
    \end{align*}
    then for some $x \in X$, $|T_{\#}(\xi - x)| \le s |\xi - x|$.
\end{enumerate}
\end{lemma}

\begin{proof}
    The proof follows \cite{Menne10_Sobolev} Lemma 3.12, with Lemma \ref{l:Menne 2.6}, Lemma \ref{l:Menne 2.10}, Lemma \ref{l:Menne 2.12}, and  Theorem \ref{t:rect varifolds compactness} in place of \cite{Menne10_Sobolev} Lemma 3.6, Lemma 3.9, Lemma 3.11, and \cite{Allard72} 6.4, respectively.
\end{proof}

\subsection{The Lipschitz map}
In order to state the Lipschitz approximation, we set the following definitions.

\begin{definition}\label{d:prep for lipschitz map}(Preliminaries for Lemma \ref{l:menne 2.18})
Let $1 \le m < n$ be integers. Given $T \in \textbf{G}(n, m)$, $a \in \mathbb{R}^{n}$, $0< r < \infty$, and $0< h \le \infty$ we define the cylindrical neighborhood $\textbf{C}(T, a, r, h)$ by
\begin{align*}
    \textbf{C}(T, a, r, h)= \mathbb{R}^n \cap \{z: |T_{\#}(z-a)|\le r \text{ and } |T_{\#}^\perp(z-a)|\le h \}.
\end{align*}
We shall use the notation $\textbf{C}(T, a, r) = \textbf{C}(T, a, r, \infty)$.

For $0<\gamma_{\Theta}< 1$ we fix $\mathbb{A}$ to be the appropriate set $\mathbb{A}= \cup_{n \in \mathbb{N}}[n, n(1+\gamma_{\Theta})]$.  Moreover, we let $Q \in \mathbb{N}$ be such that $[Q, Q(1+\gamma_{\Theta})]$ is a maximal connected component of $\mathbb{A}$.

For an open set $U$, let $V \in \textbf{RV}_m(U)$ with $\gamma_{\Theta}$-almost integral density such that $0 \in \spt \norm{V}$.  For such a $V$ and any $0< \epsilon_1$ and $\textbf{C}(T, 0, r, h)\subset U$, we define the following sets.
\begin{enumerate}
    \item Let $B(\epsilon_1)$ denote the set of all $x \in \textbf{C}(T, 0, r, h)$ with $\Theta^{*, m}(\norm{V}, z)>0$ such that
\begin{align*}
    \text{either  } \norm{\delta V} B_{\rho}^n(z)> \epsilon_1 \norm{V}(B_{\rho}^n(z))^{1-1/m} \quad \text{for some } 0<\rho<2r\\
    \text{or   } \int_{B_{\rho}^n(z)\times \textbf{G}(n, m)}|S_{\#} - T_{\#}|dV(\xi, S)> \epsilon_1 \norm{V}(B_{\rho}^n(z)) \text{ for some } 0<\rho \le 2r.
\end{align*}
\item Let $A(\epsilon_1):= \textbf{C}(T, 0, r, h) \setminus B(\epsilon_1)$.  Moreover, let $A(x) = A(\epsilon_1) \cap \{z:T_{\#}(z)=x \}$ for $x \in \mathbb{R}^m$.
\item Let $X_1$ be the set of all $x \in T \cap B_{r}^n(0)$ such that 
\begin{align*}
    \sum_{z \in A(x)} \Theta^m(\norm{V}, z) \in [Q, Q(1+\gamma_{\Theta})] \text{ and } \Theta^m(\norm{V}, z) \in \mathbb{A} \text{ for } z \in A(x),
\end{align*}
\item Let $X_2$ denote the set of all $x \in T \cap B_{r}^n(0)$ such that
\begin{align*}
    \sum_{z \in A(x)}\Theta^m(\norm{V}, z) \le (Q-1)(1+\gamma_{\Theta}) \text{ and } \Theta^m(\norm{V}, z) \in \mathbb{A} \text{  for } z \in A(x),
\end{align*}
\item Define $N = T \cap B_{r}^n(0) \setminus (X_1 \cup X_2)$,
\end{enumerate}
\end{definition}

\begin{lemma}\label{l:menne 2.18}\emph{(The Lipschitz Approximation)}
    Let $m,n, \gamma_{\Theta},\mathbb{A}$ be as in Definition \ref{d:prep for lipschitz map}.  Let $\{\delta_i\}_{i=1}^5$ and $Q \in \mathbb{N}$ be such that the following conditions hold.
    \begin{itemize}
        \item[a.]  Let $0<\delta_i <1 \text{ for } i \in \{1,2\}$ such that 
        \begin{align*}
            [Q-3\delta_1,Q - \delta_1] \cap \mathbb{A} = \emptyset, \quad Q + \delta_2 \not \in \mathbb{A}.
        \end{align*}
         \item[b.] Assume further that 
        \begin{align*}
        1-\delta_2 -3\delta_1>0.
    \end{align*}
        \item[c.] Let $0<\delta_i \le 1$ for $i \in \{3, 4, 5\}$ with $\delta_5 \le \min\{ 1/4, (2m\gamma_m)^{-m}/\omega_m\}$.
       \end{itemize}
    Let $0<L<\infty$ and $1\le M<\infty$. Then, there exists a positive, finite number $\epsilon$ such that the following holds.

    For $T \in G(n,m)$, $a \in \mathbb{R}^n$, $0<r<\infty$, and $2\delta_4r<h\le \infty$, let $U = B_{r}^n(\boldsymbol{C}(T, a, r, h))$. Suppose $V \in \textbf{\emph{RV}}_m(U)$ has locally bounded first variation and $\gamma_{\Theta}$-almost integral density.  Suppose further that $V$ satisfies
    \begin{align*}
        (Q-2\delta_1)\omega_m r^m \le \norm{V}(\boldsymbol{C}(T, a, r, h)) \le (Q+\delta_2)\omega_mr^m,\\
        \norm{V}(\boldsymbol{C}(T, a, r, h+ \delta_4 r) \setminus \boldsymbol{C}(T, a, r, h -2\delta_4r)) \le (1-\delta_3) \omega_m r^m,\\
        \norm{V}(U) \le M\omega_m r^m.
    \end{align*}
For $0<\epsilon_1 \le \epsilon$ we let $A(\epsilon_1), B(\epsilon_1)$ be as in Definition \ref{d:prep for lipschitz map} and define $H \subset A(\epsilon_1) \subset \boldsymbol{C}(T, a, r, h)$ be the set such that 
\begin{align*}
\norm{V}(\overline{B_{\rho}^n(x)}) \ge \delta_5 \omega_m\rho^m \qquad \text{for all }0<\rho<2r.
\end{align*}

Then $X_1$ is an $\mathcal{H}^m$-measurable subset of $T$ and there exists a function $f_Q:X_1 \rightarrow Q_{Q}(\mathbb{R}^m)$ with the following properties.
\begin{enumerate}
    \item $X_1 \subset \overline{B_{\rho}^m(T_{\#}(a))}$ and $\emph{Lip}(f_Q) \le L$.
    \item Then, $A(\epsilon_1), B(\epsilon_1)$ are Borel sets and for all $y \in X_1$ 
    \begin{align*}
        T^{\perp}_{\#}(A(\epsilon_1) \cap \emph{spt}\norm{V}) \subset \overline{B_{h-\delta_4r}^{n-m}(T^{\perp}_{\#}(a))}, \quad \emph{spt}f_Q(y) \subset T^{\perp}_{\#}(A(y)),\\
        \norm{f_Q(y)} = T^{\perp}_{\#}(\lfloor\Theta^m(\norm{V}, \cdot)\rfloor \mathcal{H}^{0} \res A(y)).
    \end{align*}
    \item Defining the sets $C(\epsilon_1) = \overline{B_{r}^{m}(T_{\#}(a))} \setminus (X_1 \setminus T_{\#}(B(\epsilon_1))) \subset \mathbb{R}^m$ and $D(\epsilon_1) = \boldsymbol{C}(T, a, r, h) \cap T^{-1}_{\#}(C(\epsilon_1)) \subset \mathbb{R}^n$. Then
    \begin{align*}
        \mathcal{L}^m(C(\epsilon_1)) + \norm{V}(D(\epsilon_1)) \le \Gamma_{(3)}\norm{V}(B(\epsilon_1))
    \end{align*}
    where $\Gamma_{(3)} = \max\{1+Q(1+\gamma_{\Theta}) +\frac{5^m}{(1-\alpha)\big(1-(1+\alpha)[Q\gamma_{\Theta} + 2\delta_1] + \gamma_{\Theta}\big)}, 4(Q+2)/\delta_1\}$ where $\alpha$ is as in \eqref{e:alpha}.
    \item For $x_1 \in H$, then $T^{\perp}_{\#}(x_1-a) \le h-\delta_4r$.  Furthermore, there is a $\lambda_{(4)}= \lambda_{\ref{l:Menne 2.13}}(m, \delta_2, s, Q, \gamma_{\Theta})/4$ such that for all $y \in X_1 \cap \overline{B_{\lambda_{(4)}r}^m(T_{\#}(x_1))}$ there exists an $x_2 \in A(y)$ with $\Theta^m(\norm{V}, x_2) \in \mathbb{A}$ and 
    \begin{align*}
        |T^{\perp}_{\#}(x_2  - x_1)| \le L|T_{\#}(x_2 - x_1)|.
    \end{align*}
    Moreover, $A(\epsilon_1) \cap \emph{spt}\norm{V} \subset H$ and 
    \begin{align*}
        H \cap T^{-1}_{\#}(Y) = \emph{graph} f_Q.
    \end{align*}
    \item The set $\overline{X_1} \setminus X_1$ has measure zero with respect to $\mathcal{H}^m$ and $T_{\#}(\norm{V} \res H)$.
    \item For $\lambda_{(4)}$ as above, if $\mathcal{H}^{m}(\overline{B_{r}^m(T_{\#}(a))} \setminus X_1) \le \frac{1}{2}\omega_m (\lambda_{(4)}r/6)^{m}$, $1\le q<\infty$, and $P = (\Theta^0(\norm{S}, \cdot) \circ T^{\perp}_{\#})\mathcal{H}^m$ is the $Q$-valued plane associated to $S \in Q_Q(\mathbb{R}^{n-m})$ via $T^{\perp}_{\#}$ then for $g:Y \rightarrow \mathbb{R}$ defined by
    \begin{align*}
        g(y) = \mathcal{G}(f(y), S),
    \end{align*}
    the following estimates hold.
    \begin{align*}
        \norm{\dist(\cdot , \emph{spt}P)}_{L^q(\norm{V}\res H)} & \le (12^{m+1})(Q+1) \left( \norm{g}_{L^q(\mathcal{H}^m\res X_1)} + \Gamma_{(6)} \mathcal{H}^m(\overline{B_r^m(T_{\#}(a))} \setminus X_1)^{\frac{1}{q} + \frac{1}{m}} \right),\\
        \norm{\dist(\cdot, \emph{spt}P) \res H}_{\infty} & \le \norm{g}_{L^{\infty}(\mathcal{H}^m \res X_1)} + 2(\mathcal{H}^m(\overline{B_r^m(T_{\#}(a))} \setminus X_1)/\omega_m)^{\frac{1}{m}},
    \end{align*}
    where $\Gamma_{(6)} = 2^3\omega_m^{-1/m}$.
    \item For $\mathcal{H}^m$-a.e. $y \in X_1$ the following statements hold.
    \begin{enumerate}
        \item[i.] $f_Q$ is approximately strongly affinely approximable at $y$.
        \item[ii.] Whenever $x \in H$ with $T_{\#}(x) = y$
        \begin{align*}
            T_x \emph{\spt}\norm{V} = \emph{Tan}(\emph{graph} \emph{ ap} Af_Q(y), (y, T^{\perp}_{\#}(x)))
        \end{align*}
        where $\emph{Tan}(S, a)$ denotes the classical tangent cone of $S$ at $a$ in the sense of \cite{Federer69} 3.1.21.
        \item[iii.] $\norm{S(x) - T} \le \norm{\emph{ap}Af_Q(y)}$ for $x \in H$ with $T_{\#}(x) = y$.
        \item[iv.] $\norm{\emph{ap}Af_Q(y)}^2 \le Q (1+\emph{Lip}(f_Q)^2)\max\{\norm{S(x) - T}^2: x \in T^{-1}_{\#}(\{y\}) \cap H\}$.
    \end{enumerate}
\end{enumerate}
\end{lemma}
For definitions of approximately strongly affinely approximability and $\text{ap }Af_Q$, see \cite{Menne10_Sobolev} Definition 2.2.  These concepts will not be used in the remainder of this work.

\begin{remark}\label{r:on the constants}(On the constants)
The assumptions on $\delta_1, \delta_2, Q$, and $0\le \gamma_{\Theta}<1$ imply the following statements.
\begin{enumerate}
    \item The assumption $a.$ implies $1-(Q-1)\gamma_{\Theta} - 9/4\delta_1> \frac{3}{4}\delta_1$.  The assumption $b.$ implies $\delta_1, 3\delta_1$ are lower admissible for $Q+1$.
    \item We define $0<\alpha(Q, \gamma_{\Theta}, \delta_1)$ to be the constant
\begin{align}\label{e:alpha}
    \alpha := \frac{\delta_1}{Q\gamma_{\Theta} + 2\delta_1}.
\end{align}
This condition and $a.$ ensures that
\begin{align}
    1-(1+\alpha)[Q\gamma_{\Theta} + 2\delta_1] + \gamma_{\Theta} > \gamma_{\Theta}.
\end{align}
\item We may assume $3L \le \delta_4$. 
\item We choose $0<s_0<1$ sufficiently close to $1$ such that $2(s_0^{-1} - 1)^{1/2} \le \delta_4$. Moreover, we choose $s_0\le s<1$ such that $(s^{-2}-1)^{1/2} \le \lambda_{(4)}/4$ and $Q^{1/2}(s^{-2} -1)^{1/2} \le L$.
\item Define $\epsilon>0$ to satisfy the following conditions.  First, $\epsilon$ must be less than or equal to the minimum of the following numbers
    \begin{align*}
        \epsilon_{\ref{l:menne sobolev 2.4}}(m,n, 1-\delta_3/2),\quad \epsilon_{\ref{l:menne 2.2}}(m,n, Q, \gamma_{\Theta}, \delta_2 -2\delta_1, \max\{M, 2\}, s),\\
        \epsilon_{\ref{l:Menne 2.6}}(m,n,Q+1, \delta_1, 2\delta_1, s), \quad \epsilon_{\ref{l:Menne 2.6}}(m,n,Q, \delta_1, 2\delta_1, s), \epsilon_{\ref{l:Menne 2.12}}(m,n,\delta_1/2, \delta_1, M, s),\\ \epsilon_{\ref{l:Menne 2.12}}(m,n, \delta_1, 2\delta_1, M, s), \quad \epsilon_{\ref{l:Menne 2.12}}(m,n, Q+1, M, \delta_2, s), \quad \epsilon_{\ref{l:Menne 2.12}}(m,n, Q, M, 1/4, s),\\
        \epsilon_{\ref{l:Menne 2.12}}(m,n, Q, M, \delta_2/3, s), \quad 
        \epsilon_{\ref{l:Menne 2.13}}(m,n, Q, \delta_5, \delta_2, s, s_0, M), \quad
        \epsilon_{\ref{l:Menne 2.13}}(m,n, Q, \gamma_{\Theta}, 1/2, \delta_2, s, s_0, M).
    \end{align*}
    Second, $\epsilon$ must satisfy
    \begin{align*}
        \epsilon \le (2\gamma_m)^{-1}, \quad 1-m\epsilon^2 \ge 1/2\\
        Q-2\delta_1 \le (1-m\epsilon^2)(Q-\delta_1).
    \end{align*}
\item Note that by inclusion, if $0 \le \gamma_{\Theta}' <\gamma_{\Theta}$, then $\gamma_{\Theta}'$-almost integral varifolds are also $\gamma_{\Theta}$-almost integral varifolds. This implies $$\epsilon(n,m,\gamma_{\Theta}, \delta_1, \delta_2, \delta_3, \delta_4, \delta_5, L, M) \le \epsilon(n,m,\gamma_{\Theta}', \delta_1, \delta_2, \delta_3, \delta_4, \delta_5, L, M).$$
\end{enumerate}
\end{remark}

 %   The proof of Lemma \ref{l:menne 2.18} follows the proof of \cite{Menne10_Sobolev} Lemma 2.18 with Lemma \ref{l:menne 2.2}, Lemma \ref{l:menne sobolev 2.4}, Lemma \ref{l:Menne 2.6}, Lemma \ref{l:Menne 2.8}, Lemma \ref{l:Menne 2.10}, Lemma \ref{l:Menne 2.12}, and Lemma \ref{l:Menne 2.13} in place of \cite{Menne10_Sobolev} Lemma 2.2, Lemma 2.4, Lemma 2.6, Lemma 2.8, Lemma 2.10, Lemma 2.12, and Lemma 2.13, respectively.  In the proof of (4), the allusion to \cite{Menne09_isoperimetric} 2.5 may be replaced by Corollary \ref{c: lower density rect}, above, by ensuring that $\epsilon \le \gamma_0(m,n,1/2)$.

%    To account for non-integral densities, we let $Y \subset \mathbb{R}^m$ be $X_1$ and $Z = X_2$ as in Definition \ref{d:prep for lipschitz map} and$f: Y \rightarrow Q_Q(\mathbb{R}^{n-m})$ defined so that
%    \begin{align}\label{e:f def}
%        f_Q(y) = T^{\perp}(\sum_{x \in A(y)} \lfloor\Theta^m(\norm{V}, x)\rfloor, \llbracket x \rrbracket) \quad \text{whenever } y \in Y.
%    \end{align}
%Therefore, observing that for $C= \cup_{n \in \mathbb{N}}[n, n(1+\gamma)]$ and any finite subset $\{\theta_i\}_i\subset C$, $\sum_i \lfloor \theta_i \rfloor \in C$, the argument from \cite{Menne10_Sobolev} Lemma 2.18 goes through with the above modifications.

Throughout the subsequent argument, we assume by translation and rescaling that $a = 0$, $r = 1$, and $T= \mathbb{R}^m \hookrightarrow \mathbb{R}^n$.

\textit{Proof of (1), (2)}
The measurability of $A(\epsilon_1), B(\epsilon_1)$, as well as the following basic properties of $A(\epsilon_1)$: \textit{For all $x \in A(\epsilon_1) \cap \emph{spt}\norm{V}$, $\Theta^m_*(\norm{V}, x) \ge \delta_3/2$,}
\begin{align*}
    \{\xi \in T^{-1}_{\#}(\overline{B_1^m(0)}): |T_{\#}(\xi - x)|>s|\xi-x| \} \subset (T^{\perp}_{\#})^{-1}(\overline{B_{\min\{\lambda_{(4)}/2, \delta_4\}}^{n-m}(T^{\perp}_{\#}(x))}),
\end{align*}
and $A(\epsilon_1) \cap \text{spt}\norm{V} \subset \overline{B_{h-\delta_4}^n(T)}$ follow identically as in \cite{Menne10_Sobolev} Lemma 3.15.

Next, we claim the following. \textit{If $X \subset A(\epsilon_1) \cap \emph{spt}\norm{V}$, $\Theta^m(\norm{V}, x) \in \mathbb{A}$ for all $x \in X$, and }
\begin{align*}
    |T_{\#}(x_1-x_2)| \le s|x_1 -x_2|
\end{align*}
\textit{for all $x_1, x_2 \in X$, then $\sum_{x \in X}\Theta^m(\norm{V}, x) \le Q(1+\gamma_{\Theta})$.}  To see the claim, we note that by our choice of $s$, for any $x \in A(\epsilon_1) \cap \text{spt}\norm{V}$
\begin{align*}
    \{ \xi \in B_1(T^{\perp}_{\#}(x)): |T_{\#}(\xi-x)|> s|\xi - x|\} \subset T^{-1}_{\#}(\overline{B_1^m(0)}) \cap (T^{\perp}_{\#})^{-1}(\overline{B_{\delta_4}^{n-m}(T^{\perp}_{\#}(x))}) \subset \boldsymbol{C}(T, 0,1, h).
\end{align*}
This implies that  
\begin{align*}
    \norm{V}(\cup_{x \in X}\{ \xi \in B_1^{n-m}(T^{\perp}_{\#}(x)): |T_{\#}(\xi-x)|> s|\xi - x|\}) \le \norm{V}(\boldsymbol{C}(T, 0,1, h)) \le (Q+\delta_2) \omega_m.
\end{align*}
If $\sum_{x \in X}\Theta^m(\norm{V}, x) \ge Q+1$, then we could apply Lemma \ref{l:Menne 2.6} in $B_1^n(0)$ with $Q+1$ and $\delta_{\ref{l:Menne 2.6}, 1}= \delta_1, \delta_{\ref{l:Menne 2.6}, 2}=2\delta_1$. Since $(Q+\delta_2)< Q+1- \delta_1$ we would obtain a contradiction.  Therefore, we have $\sum_{x \in X}\Theta^m(\norm{V}, x) \le Q(1+\gamma_{\Theta})$, as claimed.

Recall $X_1, X_2$ as in Definition \ref{d:prep for lipschitz map}.  The measurability of $X_1, X_2$ follows identically as in \cite{Menne10_Sobolev} Lemma 3.15.  Moreover, $\mathcal{H}^m(\overline{B_1^m(0)} \setminus (X_1 \cup X_2)) = 0$ and $X_1 \cap X_2 = \emptyset$. We let $f_Q: X_1\rightarrow Q_Q(\mathbb{R}^{n-m})$ be defined by 
\begin{align}\label{e:f def}
        f_Q(y) = \sum_{x \in A(y)} \lfloor\Theta^m(\norm{V}, x)\rfloor \llbracket T^{\perp}_{\#}(x)\rrbracket \quad \text{whenever } y \in X_1.
\end{align}
Following \cite{Menne10_Sobolev} Remark 2.14 applied to the graph of $f_Q$, one infers that for all $y_1, y_2 \in X_1$
\begin{align*}
    \mathcal{G}(f_Q(y_2), f_Q(y_1)) \le Q^{1/2}(s^{-2}-1)^{1/2}|y_1 - y_2|.
\end{align*}
This proves (1), (2). \qed

\textit{Proof of (3)}
Following the argument of \cite{Menne10_Sobolev} Lemma 3.15 verbatim, we obtain by our choice of $\epsilon$ the following \textit{coarea estimate}
\begin{align*}
    (1-m\epsilon^2)\norm{V}(\boldsymbol{C}(T, 0, 1, h) \cap T^{-1}_{\#}(W)) & \le \norm{V}(B(\epsilon_1) \cap T^{-1}_{\#}(W)) + Q(1+\gamma_{\Theta})\mathcal{L}^m(X_1 \cap W)\\
    & \qquad + (Q-1)(1+\gamma_{\Theta})\mathcal{L}^m(X_2 \cap W)
\end{align*}
for all measurable $W \subset \mathbb{R}^m$.

Moreover, by our choice of $\Gamma_{\ref{l:menne 2.18}(3)} \ge 4(Q+2)/\delta_1$ we may assume that $\norm{V}(B) \le (\delta_1/4)\omega_m$.  The mass bounds therefore imply that $\norm{V}(A(\epsilon_1))>0$ and hence $\mathcal{L}^m(X_1)>0$.  Now, by our choice of $\epsilon$ we obtain from the coarea estimate with $W = \overline{B_1^m(0)}$
\begin{align*}
    Q-2\delta_1 & \le (1-m\epsilon^2)\norm{V}(\boldsymbol{C}(T, 0, 1, h))\\
    & \le \norm{V}(B(\epsilon_1)) + Q(1+\gamma_{\Theta})\mathcal{L}^m(X_1) + (Q-1)(1+\gamma_{\Theta})\mathcal{L}^m(X_2)\\
    & \le \frac{\delta_1}{4}\omega_m + (Q-\frac{9}{4}\delta_1)\omega_m + (Q\gamma_{\Theta} + \frac{9}{4}\delta_1)\mathcal{L}^m(X_1) - (1-(Q-1)\gamma_{\Theta} - \frac{9}{4}\delta_1)\mathcal{L}^m(X_2).
\end{align*}
Hence, $\mathcal{L}^m(X_2) \le \frac{4}{3\delta_1}\mathcal{L}^m(X_1).$ 

Now, we make the following claim.  Recall $\alpha$ as in \eqref{e:alpha}. \textit{For any $z \in X_2$ with $\Theta^m(\mathcal{L}^m \res (\mathbb{R}^m \setminus X_2), z) = 0$ there exists a ball $B_t^m(\zeta) \subset T$ such that} $z \in \overline{B_t^m(\zeta)} \cap \overline{B_1^m(0)}$ \textit{and}
\begin{align*}
     \mathcal{L}^m(\overline{B_{5t}^m(\zeta)}) \le \frac{1}{(1-\alpha)\big(1 -(1+\alpha)[Q\gamma_{\Theta} + 2\delta_1] + \gamma_{\Theta}\big)}\norm{V}(B(\epsilon_1) \cap T^{-1}_{\#}(\overline{B_t^m(\zeta)})). 
\end{align*}
To see the claim, we note that since $\mathcal{L}^m(X_1)>0$ for any $z$ as above there must exist some ball $\overline{B_{t}^m(\zeta)}$ in the collection $\{\overline{B_{\theta}^m((1-\theta)z)}\}_{\theta \in (0,1]}$ such that 
\begin{align}\label{e:offset measure estimates}
    z \in \overline{B_t^m(\zeta)} \cap \overline{B_1^m(0)}, \quad 0< \mathcal{L}^m(X_1 \cap \overline{B_t^m(\zeta)}) \le \alpha \mathcal{L}^m(X_2 \cap \overline{B_t^m(\zeta)}).
\end{align}
In particular then, there exists $y \in X_1 \cap \overline{B_t^m(\zeta)}$. We consider the translation $f: \mathbb{R}^n \rightarrow \mathbb{R}^n$ defined by $f(x) = x + T_{\#}^{*}(\zeta - y)$.  Note that, viewing $y \in T \subset \mathbb{R}^n$ we have $f(y) = \zeta$.  Hence, recalling the height bounds from the basic properties of $A(\epsilon_1)$ we have that for any $\xi \in A(y)$
\begin{align*}
    \{w \in B_t(f(\xi)): |T_{\#}(y-w)|>s|y - w|\} \subset \textbf{C}(T, 0, 1, h) \cap T_{\#}^{-1}(\overline{B_t^m(\zeta)}).
\end{align*}
Therefore, applying Lemma \ref{l:Menne 2.12} with 
\begin{align*}
\delta_{\ref{l:Menne 2.12}, 2} = \delta_1, \quad \delta_{\ref{l:Menne 2.12}, 1} = \delta_{1}/2, \quad
X = A(y), \quad d=t, \quad r=1, \text{ and $f$ as above,}
\end{align*}
we obtain
\begin{align*}
    (Q-\delta_1)\omega_m t^m \le \norm{V}(\textbf{C}(T, 0, 1, h) \cap T_{\#}^{-1}(\overline{B_t^m(\zeta)})). 
\end{align*}
Applying the coarea estimate with $W = \overline{B_t^m(\zeta)}$ and recalling \eqref{e:offset measure estimates} gives further that
\begin{align*}
    (Q-2\delta_1)\omega_m t^m & \le \norm{V}(B(\epsilon_1) \cap T_{\#}^{-1}(\overline{B_t^m(\zeta)})) + Q(1+\gamma_{\Theta})\mathcal{L}^m(X_1) + (Q-1)(1+\gamma_{\Theta})\mathcal{L}^m(X_2)\\
    & = \norm{V}(B(\epsilon_1) \cap T_{\#}^{-1}(\overline{B_t^m(\zeta)})) + (Q-2\delta_1)\omega_m t^m + [Q\gamma_{\Theta} + 2\delta_1]\mathcal{L}^m(X_1 \cap \overline{B_t^m(\zeta)})\\
    & \qquad + [(Q-1)\gamma_{\Theta} + 2\delta_1 - 1]\mathcal{L}^m(X_2 \cap \overline{B_t^m(\zeta)})\\
    & \le \norm{V}(B(\epsilon_1) \cap T_{\#}^{-1}(\overline{B_t^m(\zeta)})) + (Q-2\delta_1)\omega_m t^m\\
    & \qquad + \big([Q\gamma_{\Theta} + 2\delta_1]\alpha + [(Q-1)\gamma_{\Theta} + 2\delta_1 - 1]\big)\mathcal{L}^m(X_2 \cap \overline{B_t^m(\zeta)}).
\end{align*}
Hence, recalling \eqref{e:offset measure estimates} the following inequalities hold.
\begin{align*}
    (1-\alpha) \mathcal{L}^m(\overline{B_t^m(\zeta)}) \le \mathcal{L}^m(X_2 \cap \overline{B_t^m(\zeta)}) \le \frac{ \norm{V}(B(\epsilon_1) \cap T_{\#}^{-1}(\overline{B_t^m(\zeta)}))}{\big(1-(1+\alpha)[Q\gamma_{\Theta} + 2\delta_1] + \gamma_{\Theta}\big)}.
\end{align*}
The claim follows.

The estimate $(3)$ now follows exactly as it does in \cite{Menne10_Sobolev} Lemma 3.15.  Since $\mathcal{H}^m$-a.e. $z \in X_2$ satisfies the hypotheses of the claim by taking a Vitali covering we obtain
\begin{align}
    \mathcal{L}^m(X_2) \le \frac{5^m}{(1-\alpha)\big(1-(1+\alpha)[Q\gamma_{\Theta} + 2\delta_1] + \gamma_{\Theta}\big)}\norm{V}(B(\epsilon_1)).
\end{align}
Thus, since $C(\epsilon_1) \setminus N \subset X_2 \cup T_{\#}(B(\epsilon_1))$ we estimate
\begin{align*}
    \mathcal{L}^m(C(\epsilon_1)) \le \Big(1+\frac{5^n}{(1-\alpha)\big(1-(1+\alpha)[Q\gamma_{\Theta} + 2\delta_1] + \gamma_{\Theta}\big)}\Big)\norm{V}(B(\epsilon_1)).
\end{align*}
Now, we estimate $D(\epsilon_1)$.  By the coarea estimate with $W = C(\epsilon_1)$, we have
\begin{align*}
    (1-m\epsilon^2)\norm{V}(D(\epsilon_1)) & \le \norm{V}(B(\epsilon_1)) + Q(1+\gamma_{\Theta})\mathcal{L}^m(C(\epsilon_1))\\
    & \le  \big(1+Q(1+\gamma_{\Theta}) +\frac{5^m}{(1-\alpha)\big(1-(1+\alpha)[Q\gamma_{\Theta} + 2\delta_1] + \gamma_{\Theta}\big)}\big)\norm{V}(B(\epsilon_1)).
\end{align*}
\qed

\textit{Proof of (4)} Note that since $H \subset A(\epsilon_1)$, the basic properties of $A(\epsilon_1) \cap \text{spt}\norm{V}$ give the first claim.  To prove the later claim, recall that $\lambda_{(4)} = \lambda_{\ref{l:Menne 2.13}}(m, \delta_2, s, Q, \gamma_{\Theta})/4$ for $0<\delta_2/2<\delta_2<1$ lower admissible for $Q$, and let $x_1 \in H$ and $y \in X_1 \cap \overline{B_{\lambda_{(4)}}^m(T_{\#}(x_1))}.$  Let $X = \{\xi \in A(y): \Theta^m(\norm{V}, \xi) \in \mathbb{N}\}$ and note that by Lemma \ref{l:menne 2.2} with
\begin{align*}
 \delta, M, a,  r, d, t \text{ and } \zeta \text{  replaced by }\\
1-\delta_2-2\delta_1 + \kappa, \max\{M, 2\}, x_1, 2, 1, 1, \text{ and } -T_{\#}(x_1)
\end{align*}
for $0< \kappa$ sufficiently small to obtain the estimate $$\norm{V}(T_{\#}^{-1}(\overline{B_1^m(0)}) \cap (T_{\#}^{\perp})^{-1}(\overline{B_{\min\{\lambda_{(4)}/2,\delta_4\}}^{n-m}(T_{\#}^{\perp}(x_1))})) \ge (\delta_2+2\delta_1 + \kappa)\omega_m.$$ We claim that this estimate implies that there exists some $x \in X$ such that $|T_{\#}^{\perp}(x-x_1)| \le \lambda_{(4)}$. To see the claim, we note that by the basic properties of $A(\epsilon_1) \cap \text{spt}\norm{V}$, our choice of $s_0$, and the hypotheses of the theorem
\begin{align*}
T_{\#}^{\perp}(\{\xi \in B_1^n(x): |T_{\#}(\xi-x)|> s|\xi-x|\}) \subset T_{\#}^{\perp} (\overline{B_{\min\{\lambda_{(4)}/2, \delta_4 \}}^n(x)})\\
\cup_{x \in X}\{\xi \in B_1^n(x): |T_{\#}(\xi-x)|> s|\xi-x|\} \subset \boldsymbol{C}(T, 0, 1, h)\\
\norm{V}(\boldsymbol{C}(T, 0, 1, h)) \le (Q+\delta_2)\omega_m.
\end{align*}
On the other hand, applying Lemma \ref{l:Menne 2.12} to $X$ and $Q$ with the translation $f(x) = x-T_{\#}^*(y)$, $\delta_{\ref{l:Menne 2.12}, 1} = \delta_1$, $\delta_{\ref{l:Menne 2.12}, 2} = 2\delta_1$, $t=1$, $r = 2$, $d=1$, and $s=s_0$ we see that 
\begin{align*}
   \norm{V}(\cup_{x \in X}\{\xi \in B_1^n(x): |T_{\#}(\xi-x)|> s|\xi-x|\}) \ge (Q-2\delta_1)\omega_m.
\end{align*}
Thus, if  $|T_{\#}^{\perp}(x-x_1)| > \lambda_{(4)}$ for all $x \in X$ we would have 
\begin{align*}
    (Q-2\delta_1)\omega_m + (\delta_2+2\delta_1 +\kappa)\omega_m \le \norm{V}(\boldsymbol{C}(T, 0, 1, h)) \le (Q+\delta_2)\omega_m,
\end{align*} a contradiction. Thus, there must exist a $x \in X$ such that $|T_{\#}^{\perp}(x-x_1)| \le \lambda_{(4)}$. Note that this implies $|x-x_1| \le 2\lambda_{(4)} \le \lambda_{\ref{l:Menne 2.13}}(m, \delta_2, s, Q, \gamma_{\Theta})/2 \le 1/2$.

Now, we apply Lemma \ref{l:Menne 2.13}(2) to $X$ with $s=s_0$ and
\begin{align*}
    \delta_{\ref{l:Menne 2.13}, 1}, \delta_{\ref{l:Menne 2.13}, 3}, \delta_{\ref{l:Menne 2.13}, 4}, \lambda, d, r, t, \zeta, \xi \text{ replaced by }\\
    1/2, \delta_2 + (1-\delta_2)/2, \delta_2, \lambda_{\ref{l:Menne 2.13}}(m, \delta_2, s, Q, \gamma_{\Theta}), 1, 2, 1, -T_{\#}^*(y), \text{ and } x_1.
\end{align*}
Thus, there exists some $x_2 \in X$ such that $|T_{\#}(x_1 -x_2)| \ge s |x_1 - x_2|$, as claimed.

The third claim comes from the observation that $\epsilon \le \frac{1}{2\gamma_m}$ and using Lemma \ref{l:menne sobolev 2.4} for $\kappa_1$ such that $1-\kappa_1 > 1/4 \ge \delta_5$. The claim of graphicality follows immediately from the definition of $A(\epsilon_1), Y,$ and $f_Q$.
\qed

\textit{Proof of (5)}
To prove (5), we will show that $\overline{X_1}\setminus X_1 \subset N$. Since $\mathcal{L}^m(N)=0$ and $T^{-1}(\overline{X_1}) \cap H \subset A(\epsilon_1)$.  The latter implies by the coarea formula that \begin{align*}
    (1-m\epsilon^2)T_{\#}(\norm{V} \res H) \le Q(1+\gamma_{\Theta}) \mathcal{L}^m \res X_1.
\end{align*}  The statement of (5) then follows immediately. 

To prove the claim, we let $y \in \overline{X_1}$, $y_i \in Y$ such that $y_i \rightarrow y$.  Let $s \in Q_Q(\mathbb{R}^{n})$ be such that $s = f_Q(y)$.  Let $R = T_{\#}^{-1}(y) \cap \text{spt}(s)$.  Following \cite{Federer69} 3.2.22(3), we note that $A(\epsilon_1) \cap \text{spt}\norm{V}$ is closed.  Hence, $R \subset A(\epsilon_1) \cap \text{spt}\norm{V}$.  On the other hand, Lemma \ref{l:menne 2.18}(4) implies $H \cap T_{\#}^{-1}(y) \subset A(\epsilon)$.  This gives the claim $T_{\#}^{-1}(\overline{X_1}) \cap H \subset A(\epsilon_1)$.

Now, we let $X_i = \{\xi \in A(y_i)\}$.  Applying Lemma \ref{l:Menne 2.6} with
\begin{align*}
    X, \delta_{\ref{l:Menne 2.6}, 1}, \delta_{\ref{l:Menne 2.6}, 2}, r \text{ replaced by}\\
    X_i, \delta_1, 2\delta_1, 2
\end{align*}
we see that $\norm{V}(\cup_{x \in A(y_i)}B_{\rho}^n(x)) \ge (Q-2\delta_1)\omega_m \rho^m$ for all $0<\rho \le r$.  Thus, letting $i \rightarrow \infty$ we see that since $\text{spt}f_Q(y_i) \rightarrow s$ in the Hausdorff metric on compact subsets, we see that
\begin{align*}
    (Q-2\delta_1) \le \liminf_{\rho \rightarrow 0}\frac{\norm{V}(\cup_{x \in R}B_{\rho}^n(x))}{\omega_m \rho^m} = \sum_{x \in R} \Theta^m_{*}(\norm{V}, x).
\end{align*}
Thus, $y \not \in X_2$.  Hence, if $y \not \in X_1$ then $y \in N$. \qed

\textit{Proof of (6)} Let $a, S, P, g$ be given as in the hypotheses of (6).
Note that by (4), for all $0<\kappa$
\begin{align*}
    \{x \in H \cap T_{\#}^{-1}(X_1): \dist(x, \text{spt}P) >\kappa \} \subset H \cap T_{\#}^{-1}(\{y \in X_1: g(y)>\kappa \}).
\end{align*}
Hence, trivially $\norm{\dist(\cdot, \text{spt}P)}_{L^q(\norm{V}\res H \cap T_{\#}^{-1}(X_1))} \le (1+L^2)^{1/2}Q(1+\tilde{\gamma}^2)\norm{g}_{L^q(\mathcal{L}^m\res X_1)}$.  Hence, to obtain the first estimate it remains only to estimate $\norm{\dist(\cdot, \text{spt}P)}_{L^q(\norm{V}\res H \setminus T_{\#}^{-1}(X_1))}$.

Let $z \in \overline{B_1^m(0)} \setminus \overline{X_1}$. We claim that there exists a ball $\overline{B_t^m(\zeta)}$ with the following properties.
\begin{align*}
    z \in \overline{B_t^m(\zeta)} \subset \overline{B_1^m(0)}, \quad 0<t\le \lambda_{(4)}/6\\
    \mathcal{L}^m(\overline{B_t^m(\zeta)} \cap X_1) = \mathcal{L}^m(\overline{B_t^m(\zeta)} \setminus X_1).
\end{align*}
To verify the claim, we consider the family $\{\overline{B_{\theta}^m((1-\theta)z)}: 0<\theta \le \lambda_{(4)}/6\}$ and recall the hypotheses of (6).

Taking a Vitali covering, we obtain a countable, disjoint collection $I = \{\overline{B_{t_i}^m(\zeta_i)}\}_i$ such that if we label $E_i = \overline{B_{t_i}^m(5\zeta_i)} \cap \overline{B_1^m(0)}$ then $\overline{B_1^m(0)} \setminus \overline{X_1} \subset \cup_i E_i$.  
For each ball $\overline{B_{t_i}^m(\zeta_i)}\in I$, we pick a point $y_i \in X_i \cap \overline{B_{t_i}^m(\zeta_i)}$ and define 
\begin{align*}
    h_i = \mathcal{G}(f(y_i), S), \quad X_i = A(y_i).
\end{align*} We now partition $I= J \cup K$ for
\begin{align*}
    J =\{\overline{B_{t_i}^m(\zeta_i)} \in I: h_i \ge 18t_i\}, \quad K = I\setminus J. 
\end{align*}
By (5), then we have
\begin{align*}
    \norm{\dist(\cdot, \text{spt}P)}_{L^q(\norm{V}\res H \setminus T_{\#}^{-1}(\overline{X_1}))} & \le 
    \norm{\dist(\cdot, \text{spt}P)}_{L^q(\norm{V}\res H \cap T_{\#}^{-1}(\cup_{i\in J} E_i))}\\
    & \qquad + \norm{\dist(\cdot, \text{spt}P)}_{L^q(\norm{V}\res H \cap T_{\#}^{-1}(\cup_{i\in K} E_i))}.
\end{align*}

Now, we claim the following: \textit{For all $i \in I$, if $x_1 \in H \cap T_{\#}^{-1}(E_i)$ then $\dist(x_1, \emph{spt} P) \le 6t_i + h_i$.} To verify the claim, we note that because $P$ is parallel to $T$ for any $x_2 \in  X_i$, $\dist(x_1, \text{spt}P) \le |T_{\#}^{\perp}(x_1 -x_2)| + h_i$.  Noting that $|T_{\#}(x_1) - y_i| \le 6t_i \le \lambda_{(4)}$ allows us to invoke (4), which gives an $x_2 \in A(y_i)$ such that $|T_{\#}^{\perp}(x_1-x_2)| \le L |T_{\#}(x_1-x_2)| = L|T_{\#}(x_1)-y| \le 6t_i$.  This gives the claim.  

For the $x_2 \in X_i$ constructed, above, we note that $|x_1-x_2| = 12t_i$.  Thus, we have $H \cap T_{\#}^{-1}(E_i) \subset \cup_{x \in X_i}\overline{B_{12t_i}^m(x)}$.  Therefore, recalling the basic properties of $A(\epsilon) \cap \text{spt}\norm{V}$ we may apply Lemma \ref{l:Menne 2.13} with
\begin{align*}
    \delta_{\ref{l:Menne 2.13}, 1}, \delta_{\ref{l:Menne 2.13}, 3}, \delta_{\ref{l:Menne 2.13}, 4}, s, \lambda, X, d, r, t, \zeta, \tau \text{ replaced by} \\
    1/2, \delta_2-(\delta_2 - Q\gamma_{\Theta} )/2, \delta_2, s_0, \lambda_{\ref{l:Menne 2.13}}(m, \delta_2, s_0, Q, \gamma_{\Theta}), 1,2,2, -T_{\#}^{*}(y_i), 12t_i.
\end{align*}
This gives $\norm{V}\res H(T_{\#}^{-1}(E_i)) \le (Q+ \delta_2)\omega_m(12t_i)^m$.  

Now, we estimate $\norm{\dist(\cdot, \text{spt}P)}_{L^q(\norm{V}\res H \cap T_{\#}^{-1}(\cup_{i\in J} E_i))}$.  For $i \in J$, we note that $h_i \ge 18t_i$ implies $\dist(x, \text{spt}P) \le (4/3)h_i$ for all $x \in H \cap T_{\#}^{-1}(E_i)$. Hence, if we label $J(\gamma):= \{i \in J: (4/3)h_i >\gamma\}$ we see that 
\begin{align*}
    \norm{V} \res H(T_{\#}^{-1}(\cup_{i\in J}E_i) \cap \{x \in \mathbb{R}^n: \dist(x, \text{spt}P)>\gamma \}) & \le \sum_{i \in J(\gamma)}\norm{V} \res H (T_{\#}^{-1}(E_i))\\
    & \le \sum_{i \in J(\gamma)}(Q+ \delta_2)\omega_m(12t_i)^m\\
    & \le \sum_{i \in J(\gamma)}(Q+ \delta_2)12^m\mathcal{L}^m(\overline{B_{t_i}^m(\zeta_i)})\\
    & \le 2(Q+ \delta_2)12^m\mathcal{L}^m(\cup_{i \in J(\gamma)}\overline{B_{t_i}^m(\zeta_i)} \cap X_1)
\end{align*}
But, for all $y \in X_1 \cap \overline{B_{t_i}^m(\zeta_i)}$ 
\begin{align*}
    \mathcal{G}(f(y), S) \ge \mathcal{G}(f(y_i), S) - L|y-y_i| \ge h_j - 2Lt_i \ge (2/3)h_i.
\end{align*}
Whence, we obtain 
\begin{align*}
    &\norm{V} \res H(T_{\#}^{-1}(\cup_{i\in J}E_i) \cap \{x \in \mathbb{R}^n: \dist(x, \text{spt}P)>\gamma \})\\
    & \qquad \le 2(Q+ \delta_2)12^m\mathcal{L}^m(\{y \in X_1: \mathcal{G}(f(y), S)>\gamma/2 \}).
\end{align*}
Thus, integrating one obtains
\begin{align*}
    \norm{\dist(\cdot, \text{spt}P)}_{L^q(\norm{V}\res H \cap T_{\#}^{-1}(\cup_{i\in J}E_j))} \le 4(Q+\delta_2)12^m\norm{g}_{L^q(\mathcal{L}^m \res X_1)}.
\end{align*}

To obtain an estimate on $\norm{\dist(\cdot, \text{spt}P)}_{L^q(\norm{V}\res H \cap T_{\#}^{-1}(\cup_{i\in K}E_j))}$ we note that for $i \in K$ and $x \in H \cap T^{-1}(E_i)$ then
\begin{align*}
    \dist(x, \text{spt} P) \le 24t_i. 
\end{align*}
Thus, for $0<\gamma<\infty$ we define $K(\gamma):= \{i \in K: 24t_i>\gamma \}$ and $u:\mathbb{R}^m \rightarrow \mathbb{R}$ by $u:\sum_{i\in I}2t_i \chi_{\overline{B_{t_i}^m(\zeta_i)}}$.  Thus, we estimate as follows.
\begin{align*}
    & \norm{V}\res H(T_{\#}^{-1}(\cup_{i \in K}E_i) \cap \{x \in \mathbb{R}^n: \dist(x, \text{spt}P)>\gamma \})\\
    & \qquad \le \sum_{i \in K(\gamma)}\norm{V}\res H(T_{\#}^{-1}(E_i))\\
    & \qquad \le \sum_{i \in K(\gamma)}(Q+\delta_2)\omega_m(12t_i)^m \le (Q+\delta_2)12^m\mathcal{L}^m(\cup_{i \in K(\gamma)}\overline{B_{t_i}^m(\zeta_i)})\\
    & \qquad \le (Q+\delta_2)12^m \mathcal{L}^m(\{y \in \mathbb{R}^m: u(y)>\gamma/12\}). 
\end{align*}
Hence, integrating we see that 
\begin{align*}
    \norm{\dist(\cdot, \text{spt}P)}_{L^q(\norm{V}
    \res H \cap T_{\#}^{-1}(\cup_{i\in K} E_i))} \le (Q+\delta_2)12^{m+1}\norm{u}_{L^q(\mathcal{L}^m)}.
\end{align*}
Combining this with the estimates
\begin{align*}
    \mathcal{L}^m(\cup_{i \in I}\overline{B_{t_i}^m(\zeta_i)}) \le 2\mathcal{L}^m(\overline{B_{1}^m(0)} \setminus X_1),\\
    \int |u|^qd\mathcal{L}^m = \sum_{i \in I}(2t_i)^q\omega_m t_i^m \le 2^q\omega_m^{-q/m}(\sum_{i \in I}\mathcal{L}^m(\overline{B_{t_i}^m(\zeta_i)}))^{1+q/m},\\
    \norm{u}_{L^q(\mathcal{L}^m)} \le 2^3\omega_m^{-1/m}\mathcal{L}^m(\overline{B_1^m(0)} \setminus Y)^{\frac{1}{q} +\frac{1}{n}},
\end{align*}
the first conclusion of (6) follows.

To prove the $q=\infty$ case, let $x_1 \in H$.  For all $(\mathcal{L}^m(\overline{B_1^m(0)} \setminus X_1)/\omega_m)^{1/m}< \theta <1$
\begin{align*}
T(x_1) \in \overline{B_{\theta}^m((1-\theta)T(x_1))}, \quad \mathcal{L}^m(\overline{B_{\theta}^m((1-\theta)T(x_1))})>0.    
\end{align*}
Thus, for all $0<\delta$ there exists a choice of $\theta$ and a $y \in X_1$ such that
\begin{align*}
    \mathcal{G}(f_Q(y), S) \le \norm{g}_{L^{\infty}(\mathcal{L}^m)}, \quad |T(x_1) - y| \le 2(\mathcal{L}^m(\overline{B_1^m(0)} \setminus X_1)/\omega_m)^{1/m}+\delta.
\end{align*}
Applying (4), we obtain $x_2 \in A(y)$ such that $|T_{\#}^{\perp}(x_2 - x_1)| \le L|T_{\#}(x_1 - x_2)| \le |T_{\#}(x_1) - y|$.  Therefore, because $P$ is parallel to $T$ we see 
\begin{align*}
    \dist(x_1, \text{spt}P) & \le \dist(x_2, \text{spt}P) + |T_{\#}^{\perp}(x_1- x_2)|\\
    & \le \mathcal{G}(f(y), S) + 2(\mathcal{L}^m(\overline{B_1^m(0)} \setminus X_1)/\omega_m)^{1/m} + \delta.
\end{align*}
Letting $\delta \rightarrow 0$, we obtain the second claim. \qed

\textit{Proof of (7)}
The proof of (7) follows as in the proof of \cite{Menne10_Sobolev} Lemma 3.15(7) verbatim. \qed

\section{A Sobolev-Poincar\'e Inequality}\label{s:sobolev poincare}

As an immediate consequence of Lemma \ref{l:menne 2.18}, we follow \cite{Menne10_Sobolev} and prove a Sobolev-Poincar\'e inequality that relates the $Q$-height with the tilt-excess, appropriately generalized for $L^p$ spaces.  In order to state the inequality, we need the following definitions to clarify what we mean by the $Q$-height of a varifold.

\begin{definition}
Let $1 \le m \le n$ be integers, $Q \in \mathbb{N}$, and $V\in \textbf{RV}(B_1^n(0))$.  Let $T \in G(n,m)$, $a \in B_1^n(0)$, and $0<r,h < 1$ be such that $\boldsymbol{C}(T, a, r, h) \subset B_1^n(0)$, where $\boldsymbol{C}(T, a, r, h)$ is as in Definition \ref{d:prep for lipschitz map}.

We define the \textit{$q$ tilt of $V$ with respect to $T$} as follows.
    \begin{align*}
        T_q(V, a, r, h, T): = r^{-m/q}\norm{S(\cdot) - T}_{L^q(\norm{V} \cap \boldsymbol{C}(T, a, r, h))}.
    \end{align*}
    
To define the $q$ height, we need the following preliminary definitions.  \begin{itemize}
\item Let $A \subset T \cap B_r^m(T_{\#}(a))$ be the set of points $x$ such that there exist $r(x) \in Q_Q(\mathbb{R}^{n-m})$ such that 
\begin{align*}
    \norm{r(x)} = \lfloor \Theta^m(\norm{V} \res \boldsymbol{C}(T, a, r, h), \cdot ) \rfloor \mathcal{H}^0 \res T^{-1}_{\#}(x).
\end{align*}
Moreover, we note that $A$ is $\mathcal{H}^m$-measurable by \cite{Menne10_Sobolev} Definition 4.1.
\item Let $P$ be a $Q$-valued plane parallel to $T$, and let $s: T \rightarrow Q_Q(\mathbb{R}^{n-m})$ be its corresponding defining function. For any $x \in A$, above, we define 
\begin{align*}
    g_P(x) := \mathcal{G}(r(x), s(x)).
\end{align*}
Note that $g:A \rightarrow \mathbb{R}$ is $\mathcal{H}^m$-measurable.
\item Let $P$ be a $Q$-valued plane parallel to $T$.  We define 
\begin{align*}
    H_q(V, a, r, h, P):= r^{-1-n/q}\norm{\text{dist}(\cdot, \text{spt}P)}_{L^q(\norm{V} \res \boldsymbol{C}(T, a, r, h))}.
\end{align*}
Letting $s: T \rightarrow Q_Q(\mathbb{R}^{n-m})$ be the corresponding defining function for $P$, we define 
\begin{align*}
    \tilde{H}_q(V, a, r, h, P) := \inf_{Y \subset A} \big\{r^{-1-n/q} \norm{g_P}_{L^q(\mathcal{H}^m \res Y)} + r^{-1-n/q}\mathcal{H}^m(T \cap \overline{B_r(T(a))} \setminus Y)^{\frac{1}{q}+\frac{1}{n}}\big\},
\end{align*}
where the infimum is over all $\mathcal{H}^m$-measurable subsets $Y \subset A$.
\end{itemize}
With these definitions, we define the $q$ \textit{height of} $V$ as follows.
\begin{align}
    H_q(V, a,r, h, Q, T) := \inf_{P}(H_q(V, a, r, h, P) + \tilde{H}_q(V, a, r, h, P)),
\end{align}
where the infimum is taken over all $Q$-valued planes parallel to $T.$
\end{definition}

\begin{lemma}\label{l:menne 3.4}\emph{(Sobolev-Poincar\'e Estimates, cf. \cite{Menne10_Sobolev} Theorem 4.4)}
    Let $1 \le m\le n$ be integers. Let $1\le M < \infty$, $0<\gamma_{\Theta}<1$, and $\mathbb{A}= \cup_{n \in \mathbb{N}}[n, n(1+\gamma_{\Theta})]$. Let $Q \in \mathbb{N}$ and $0<\delta_1, \delta_2, \delta_3, \delta_4 \le 1$ satisfy 
    \begin{align*}
        0<\delta_i <1 \text{ for } i \in \{1,2\} \text{ such that } [Q-3\delta_1,Q - \delta_1] \cap \mathbb{A} = \emptyset, \quad Q + \delta_2 \not \in \mathbb{A} \\
        0<\delta_i \le 1 \text{ for } i \in \{3, 4\}\\
        1-\delta_2 - 3\delta_1>0
    \end{align*}
    Then, there exists a $0<\epsilon< \infty$ such that the following holds.

    Let $a \in \mathbb{R}^n$, and $0<r<\infty$, $0<h\le \infty$ satisfy $\delta_4r >h$.  Let $T \in G(n,m)$, and suppose that $V \in \textbf{\emph{RV}}_m(\textbf{C}(T, a, 3r, h+2r))$ be a varifold with $\gamma_{\Theta}$-almost integral density and locally bounded first variation.  Suppose further that $V$ satisfies
    \begin{align*}
        (Q-2\delta_1)\omega_m r^m \le \norm{V}(\boldsymbol{C}(T, a,r, h)) \le (Q+\delta_2)\omega_m r^m,\\
        \norm{V}(\textbf{C}(T, a,r, h+\delta_4r) \setminus\textbf{C}(T, a,r, h-2\delta_4r)) \le (1-\delta_3)\omega_m r^m\\
        \norm{V}(B_{r}(\boldsymbol{C}(T, a, r, h))) \le M \omega_mr^m.
    \end{align*}
    Moreover, suppose that 
    \begin{align*}
        \norm{\delta V}(\boldsymbol{C}(T, a, 3r, h+2r)) \le \epsilon r^{n-1}, \qquad T_1(V, a, 3r, h+2r, T) \le \epsilon.
    \end{align*}
    Then, if we define for $0<\tau$
    \begin{align*}
        G(\tau):= \{x \in \textbf{C}(T, a, r, h) \cap \emph{spt}\norm{V}: \norm{\delta V}(\overline{B_\rho^n(x)}) \le \tau \norm{V}(\overline{B_{\rho}^n(x)})^{1-1/m} \text{ for all } 0<\rho<2r\},
    \end{align*}
    the following statements hold.
    \begin{enumerate}
        \item If $1\le q< m$ and $q^* = mq/(n-q)$, then
        \begin{align*}
            H_{q^*}(V \res G(\frac{1}{2\gamma_m}), a, r,h, Q, T) \le \Gamma_{\ref{l:menne 3.4}}\left( T_q(V, a, 3r, h+2r, T) + \left(\frac{\norm{V}(\boldsymbol{C}(T, a, r, h) \setminus G(\epsilon))}{r^m} \right)^{1/q} \right),
        \end{align*} where $\Gamma_{\ref{l:menne 3.4}} = \Gamma(m,n,Q, M, \delta_1, q)$
        \item If $m<q\le \infty$, then
        \begin{align*}
            H_{\infty}(V \res G(\frac{1}{2\gamma_m}), a, r,h, Q, T) \le \Gamma_{\ref{l:menne 3.4}}\left( T_q(V, a, 3r, h+2r, T) + \left(\frac{\norm{V}(\boldsymbol{C}(T, a, r, h) \setminus G(\epsilon))}{r^m} \right)^{1/q} \right),
        \end{align*} where $\Gamma_{\ref{l:menne 3.4}} = \Gamma(m,n,Q, M, \delta_1, q)$.
    \end{enumerate}
\end{lemma}

\begin{proof}
    The proof follows \cite{Menne10_Sobolev} Theorem 4.4 with Lemma \ref{l:menne 2.18} in place of \cite{Menne10_Sobolev} Lemma 3.15.
\end{proof}

\begin{remark}(On the constants in Lemma \ref{l:menne 3.4})
The $\Gamma_{\ref{l:menne 3.4}}$ may be calculated as in the proof of \cite{Menne10_Sobolev} Lemma 4.4.  Letting $0<\epsilon_0=\epsilon_{\ref{l:menne 2.18}}$ with
\begin{align*}
    m,n,Q, L=1, M, \delta_1, \delta_2, \delta_3, \delta_4, \delta_5= \min\{1/4, (2m\gamma_m)^{-m}/\omega_m\}.  
\end{align*}
We let $\epsilon_1:=\epsilon_0$ and $\lambda = \lambda_{\ref{l:menne 2.18}(4)}(m, \delta_2, s, Q, \gamma_{\Theta})$.  Then there is a dimensional constant $N(n)$ such that the $0<\epsilon$ of Lemma \ref{l:menne 3.4} must satisfy
\begin{align*}
    \epsilon \le \epsilon_0(n\gamma_m)^{1-m}, \qquad 3^m\epsilon \le \epsilon_0(n\gamma_m)^{-m}\\
    \Gamma_1N(n)3^m\epsilon \le \omega_1(\lambda/12)\epsilon_0 \quad \text{if $m=1$}\\
    \Gamma_1N(n)(3^m\frac{\epsilon}{\epsilon_0} + (\frac{\epsilon}{\epsilon_0})^{\frac{m}{m-1}})\le \frac{1}{2}\omega_m(\lambda/6)^m \quad \text{if $m>1$}.
\end{align*}
\end{remark}

\section{The Height Bound}\label{s:height bound}

The height bound in Theorem \ref{t:Allard Lemma 8.4} is sufficient for the proof of Theorem \ref{t:BMO main theorem rect} and Theorem \ref{t: VMO main theorem rect}.  However, we wish to obtain more quantitative control of the height in terms of the $\tilde{\gamma}$ as in the $(\delta_0, \tilde{\gamma})$-Allard tilt-excess decay regime. Such quantitative control is easily obtained by a number of means.  For example, one may combine the arguments of \cite{Simon83_book} Theorem 19.5 and Theorem 20.2 to see that under the $(\delta_0, \tilde{\gamma})$-Allard tilt-excess decay regime for $0<\eta$ as in \eqref{e:eta scale requirements} an estimate of the form 
\begin{align}\label{e:height bound intro}
\sup_{x \in \text{spt}\norm{V} \cap B_{\eta r}^n(0)} \text{dist}(x, T) \le C(m,n,\delta_0)\tilde{\gamma}^{\frac{1}{m+1}}r
\end{align}
holds for $0< \tilde{\gamma}$ sufficiently small. Below, we obtain \eqref{e:height bound intro} as a consequence of the Sobolev-Poincar\'e estimates of Lemma \ref{l:menne 3.4}.

In order to apply Lemma \ref{l:menne 2.18} and Lemma \ref{l:menne 3.4}, we need some preliminary density estimates.

\begin{remark}\label{r: extend 1-lower density to balls}
    If $V$ is as in \ref{s: the set up} and satisfies the critical $(\delta_0, \tilde{\gamma})$-Allard tilt-excess decay regime in $B_r(a)$ then for $0<\tilde{\gamma} \le \gamma_0(m,n,\kappa_1)$ we may apply Lemma \ref{l:menne sobolev 2.4} at all $x \in \text{spt}\norm{V} \cap B_{\frac{r}{4}}(a)$ to obtain
    \begin{align}\label{e:apply 1-lower density}
    1-\kappa_1 \le \frac{\norm{V}(B_{\rho}^n(x))}{\omega_m \rho^m}
\end{align}
for all $0<\rho <r/2$.
\end{remark}

\begin{definition}\label{d:eta and eta 0}
    For $0<\delta_0 \le 1/4$ we consider $0< \eta$ such that 
\begin{align}\label{e:eta scale requirements}
    \frac{(2-\delta_0)}{(1-\eta)^m} \le 2-\delta_0/2.
\end{align}    
Moreover, we define $0<\eta_0(m,\delta_0)$ such that $\frac{(2-\delta_0)}{(1-\eta_0)^m} = 2-\delta_0/2$. Note that $\delta_0 \le 1/4$ implies that $\eta_0 \le 1/15$.
\end{definition}

\begin{lemma}\label{l:upper density bound 1 scale}\emph{(Upper Density Bound)} Let $0<\delta_0\le 1/4$ be given.  Let $\kappa_1 = \delta_0/10$ in Lemma \ref{l:menne sobolev 2.4}.  For any $0<\eta$ such that \eqref{e:eta scale requirements} holds there is a $0< \gamma_1(m, \delta_0, \eta)$ such that if $V$ satisfies hypotheses of the critical $(\delta_0,\tilde{\gamma})$-Allard tilt-excess decay regime in $B_{r}^n(0)$ with $\tilde{\gamma} \le \min\{\gamma_0(m,n, \delta_0), \gamma_1(m, \delta_0, \eta)\}$ then for all $x \in B_{\eta_0 r}^n(0) \cap \emph{spt}\norm{V}$ and $\sigma \in [\eta r, (1-\eta_0)r]$
    \begin{align*}
        \frac{\norm{V}(B_{\sigma}^n(x))}{\omega_n \sigma^n} \le (2-\delta_0/4).
    \end{align*}
\end{lemma}

\begin{proof}
    Let $x \in B_{\eta_0 r}^n(0)$.  By the monotonicity formula (Lemma \ref{l:monotonicity formula}), with the choice $h \equiv 1$, $\sigma$ as in the hypotheses, and $\rho = (1-\eta_0)r$ we obtain
\begin{align*}
       \frac{\norm{V}(B_{\sigma}^n(x))}{\omega_m \sigma^m} & \le \frac{\norm{V}(B_{(1-\eta_0)r}^n(x))}{\omega_m (1-\eta_0)^mr^m} + \frac{1}{\omega_n}\int_{\sigma}^{\rho}\tau^{-m}\int_{B_{\tau}^n(x)}|\boldsymbol{h}(V; z)|d\norm{V}(z)d\tau\\
       & \le \frac{(2-\delta_0)}{(1-\eta_0)^m} + \frac{1}{\omega_n}r \sigma^{-m}\int_{B_{r}^n(0)}|\boldsymbol{h}(V; z)|d\norm{V}(z)\\
       & \le \frac{(2-\delta_0)}{(1-\eta_0)^m} + \frac{1}{\omega_n}r \sigma^{-m}\gamma_{\hbar,m}(0,r) \norm{V}(B_r^n(0))^{1-1/m} \\
       & \le (2-\delta_0/2) + \eta^{-m}(2-\delta_0)^{1-1/m}\omega_m^{-1/m}\gamma_{\hbar,m}(0, r).
\end{align*}
Thus, if $\gamma_{\hbar,m}(0,r) \le \tilde{\gamma}^2 \le \gamma_2(m, \delta_0, \eta)^2$ is such that $\eta^{-m}(2-\delta_0)^{1-1/m}\omega_m^{-1/m}\gamma_{2}^2 \le \delta_0/8$, the claim of the lemma follows.
\end{proof}

\begin{remark}
    We note that similar upper density bounds may be obtained in balls larger than $B^n_{\eta_0 r}(0)$ by limit-compactness arguments.
\end{remark}

\subsection{Applying Lemma \ref{l:menne 2.18} and Lemma \ref{l:menne 3.4} in the Allard tilt-excess decay regime}

In this subsection, we verify that for $Q=1$, $0<\gamma_{\Theta}<1$, and sufficiently small $\tilde{\gamma}$ varifolds which satisfy the $(\delta_0, \tilde{\gamma})$-Allard tilt-excess decay regime satisfy the hypotheses of Lemma \ref{l:menne 2.18} and Lemma \ref{l:menne 3.4}.

\begin{lemma}\label{c:apply Lip approx rect}\emph{(Applying Lemma \ref{l:menne 2.18})} Let $0<\delta_0\le 1/4$ and $0<\gamma_{\Theta}<1$ be given. Let $\kappa_1 = \delta_0/10$ and $\gamma_0(m, n, \delta_0)$ be as in Lemma \ref{l:menne sobolev 2.4}.  Let $0<\eta_0$ as in \eqref{e:eta scale requirements} and $\gamma_1(m, \delta_0, 1/4)$ by as in Lemma \ref{l:upper density bound 1 scale}.

Then, there is a $0< \gamma_{2}(m, n, \delta_0, \gamma_{\Theta})$ such the following holds.  For any $V \in \textbf{\emph{RV}}_m(B_1^n(0))$ with $\gamma_{\Theta}$-almost integral density which satisfies the hypotheses of the critical $(\delta_0, \tilde{\gamma})$-Allard tilt-excess decay regime in $B_1^n(0)$ with $\tilde{\gamma} \le \min\{\gamma_0, \gamma_1, \gamma_2\}$, we may apply Lemma \ref{l:menne 2.18} in $\boldsymbol{C}(T,x,\rho, \rho)$ for any $x \in B_{\eta_0 r}^n(0)$ and any $\rho \in [r/4, (1-\eta_0)r/3]$ with constants
    \begin{align*}
        M = (2-\delta_0/4)(1+\sqrt{2})^m, &  \qquad L = 1/24\\
       \delta_1 = (1-\max\{1-\delta_0/8, (1+\gamma_{\Theta})/2\})/4, &  \qquad  \delta_2 = \max\{1-\delta_0/8, (1+\gamma_{\Theta})/2\},\\
       \delta_3=1, & \qquad \delta_4= 1/8,\\ & \qquad \delta_5 = \min\{ 1/4, (2m\gamma_m)^{-m}/\omega_m\}.
    \end{align*} 
Note that for these choices
\begin{enumerate}
    \item $1-\delta_2-3\delta_1 \ge \delta_1>0$.
    \item We may choose the constant $\alpha(1, \gamma_{\Theta}, \delta_1)$ to satisfy $\alpha \in [\delta_1, \frac{1}{2}]$. 
\item We have $\Gamma_{(3)} \le \max\{3 + 4\frac{5^m}{\delta_1 + \gamma_{\Theta}}, 12/\delta_1\}$.
\end{enumerate}
\end{lemma}

\begin{proof}
Let $x \in B_{\eta_0 r}^n(0) \cap \text{spt}\norm{V}$ and $\rho \in [r/4, (1-\eta_0)r/3]$. Let $U_\rho(x)= B_{\rho}^n(\boldsymbol{C}(T, x, \rho, \rho))$.  Note that $U_{\rho}(x) \subset B_{(1+\sqrt{2})\rho}^n(x) \subset B_1^n(0)$.

By Lemma \ref{l:upper density bound 1 scale} for $\tilde{\gamma} \le \gamma_1(m,\delta_0, 1/4)$ we have $\norm{V}(U_{\rho}(x)) \le (2 - \delta_0/4)\omega_m (1+\sqrt{2})^m\rho^m$.  Thus, we may take $M = (2 - \delta_0/4) (1+\sqrt{2})^m$. On the other hand, since $B_{\rho}^n(x) \subset U_{\rho}(x)$ we may take $\kappa_1 = \delta_0/10$ and by Lemma \ref{l:menne sobolev 2.4}, $\tilde{\gamma} \le \gamma_0(m,n, \delta_0)$ then
\begin{align*}
    (1- \delta_0/10)\omega_m \rho^m \le \norm{V}(B_{\rho}^n(x)).
\end{align*}

In order to determine $\delta_1, \delta_2, \delta_3, \delta_4$ we let $0< \kappa_2$ be a small constant to be determined later and $\delta = \kappa_2 r/4$. Let $\epsilon(m, \kappa_2)$ be as in Lemma \ref{t:Allard Lemma 8.4} so that if $\tilde{\gamma} \le \epsilon(m, \kappa_2)$ then 
\begin{align}\label{e:small scale height control}
    \text{spt}\norm{V} \cap \boldsymbol{C}(T, x, (1-\eta_0)r/3, (1-\eta_0)r/3) \subset B_{\kappa_2 r/4}^n(T).
\end{align}
Thus, in particular, $\text{spt}\norm{V} \cap \boldsymbol{C}(T, x, \rho, \rho) \subset B_{\sqrt{\rho^2+\kappa_2^2 (r/4)^2}}^n(x)$.  Applying Lemma \ref{l:upper density bound 1 scale}, we see that  
\begin{align*}
    \norm{V}(\boldsymbol{C}(T, x, \rho, \rho)) \le (2 - \delta_0/4)\omega_m (\sqrt{1+\kappa_2^2(\frac{r}{4\rho})^2})^m\rho^m.
\end{align*}
Thus, if $\kappa_2 \le \kappa(m, \delta_0)$ is small enough to ensure $(2 - \delta_0/4) (\sqrt{1+\kappa_2^2(\frac{r}{4\rho})^2})^m 
\le 2-\delta_0/8$. Thus, we may take $\delta_2 \ge 1- \delta_0/8$. In the case that $1+\gamma_{\Theta} \ge 1-\delta_0/8$, we take $\delta_2 = \max\{1-\delta_0/8, (1+\gamma_{\Theta})/2\}$. In order to satisfy the condition that $1-\delta_2 -3\delta_1>0$, we take $\delta_1 = (1-\delta_2)/4>0$.

By taking $\kappa_2$ smaller, if necessary, \eqref{e:small scale height control} implies we may take $\delta_4 = 1/8$ and $\delta_3 = 1$. Finally, we shall fix $L = \frac{1}{24} = \frac{1}{3}\delta_4$.  The remaining claims follow immediately from the definitions of $\delta_1, \delta_2$, and $\Gamma_{(3)}$.
\end{proof}

\begin{remark}\label{r:B(delta)}\label{c:apply lipschitz approx 2}
Assume that $Q=1$ and $\delta_0, \gamma_{\Theta}, V, \eta_0, \rho, \tilde{\gamma},$ and $x$ are as in Lemma \ref{c:apply Lip approx rect}. If we let $f_Q:T \rightarrow \mathbb{R}^{n-m}$ be the Lipschitz function guaranteed by Lemma \ref{c:apply Lip approx rect} and Lemma \ref{l:menne 2.18} then the following properties hold.
\begin{enumerate}
    \item Let $0<\epsilon(n,m,\delta_0, \gamma_{\Theta})$ be the $\epsilon$ guaranteed by Lemma \ref{l:menne 2.18} for the choices of $M, \delta_i, L$ as in Lemma \ref{c:apply Lip approx rect}.  For 
    \begin{align}\label{e: gamma less than delta}
        0<\gamma_3(m,n,\delta_0, \gamma_{\Theta})^2 \le \epsilon_1 =\epsilon(m,n,\delta_0, \gamma_{\Theta})
    \end{align} if $\tilde{\gamma} \le \gamma_3$, then for all $z \in B(\epsilon_1) \in \boldsymbol{C}(T, x, \rho, \rho)$ cannot satisfy 
    \begin{align*}
        \norm{\delta V} (B^n_{\tau}(z))> \epsilon_1 \norm{V}(B^n_\tau(z))^{1-1/m}
    \end{align*}
for any $0<\tau<2\rho$.  Indeed, we compute directly
    \begin{align*}
        \norm{\delta V}( B^n_{\tau}(z)) & = \int_{B^n_{\tau}(z)}|\boldsymbol{h}(V; x)|d\norm{V}(x)\\
        & \le \gamma_{\hbar,m}(0,r)\norm{V}(B^n_{\tau}(z))^{1-1/m} \le \tilde{\gamma}^2 \norm{V}(B^n_{\tau}(z))^{1-1/m}. 
    \end{align*}
Thus, the hypothesis \eqref{e: gamma less than delta} implies that for all $z \in B(\epsilon_1)$ there must exist a radius $0< r_z \le 2\rho$ such that
\begin{align*}
    \int_{B^n_{r_z}(z)\times G(n, m)}|S_\# - T_\#|dV(\xi, S) > \epsilon_1 \norm{V}(\overline{B_{r_z}(z)}).
\end{align*}
By H\"older's inequality, $\int_{B^n_{r_z}(z)\times G(n, m)}|S_{\#} - T_{\#}|^2dV(\xi, S) > \epsilon_1^2 \norm{V}(\overline{B_{r_z}(z)})$.
Considering $\{\overline{B^n_{r_z}(z)}\}_{z \in B(\epsilon_1)}$, by Besicovitch's covering theorem there exist at most $C(n)< \infty$ families $\mathscr{F}_i$ such that $\mathscr{F}_i = \{\overline{B^n_{r_j}(z_j^i)}\}_i$ is a pairwise disjoint collection and $B(\delta) \subset \cup_i \cup_{\mathscr{F}_i}\overline{B^n_{r_j}(z_j^i)}$.  Thus, we may bound 
\begin{align}\label{e:bad set bound} \nonumber
    \norm{V}(B(\epsilon_1)) & \le \cup_i \cup_{\mathscr{F}_i}\norm{V}(\overline{B^n_{r_j}(z_j^i)})\\ \nonumber
    & \le C(n) \epsilon_1^{-2} \int_{B^n_{2\rho}(x)}|S_{\#} - T_{\#}|^2dV(\xi, S)\\
    & \le C(n) \epsilon_1^{-2} E(0, r, T)^2 r^m.
\end{align}
    \item By truncation, we may assume that for all $x\in T \cap B^n_{\rho}(0)$, $\max\{|y|: y \in f_Q(x)\} \le \kappa_2 r/4$, where $0<\kappa_2(m, \delta_0)$ is as in the proof of Lemma \ref{c:apply Lip approx rect}.  Note that $\delta_0 \le 1/4$ implies that $\kappa_2<1/5$.
    \item By \eqref{e:f def} for any $x \in X_1$, $\mathcal{H}^0(A(x)) = 1$ and if $y \in A(x)$, $\Theta^m(\norm{V}, y) \in [1, 1+ \gamma_{\Theta}]$.
\end{enumerate}
\end{remark}

Following Lemma \ref{c:apply Lip approx rect}, we have the following statement allowing us to apply Lemma \ref{l:menne 3.4}.

\begin{lemma}\label{l:apply 3.4}\emph{(Applying Lemma \ref{l:menne 3.4})} Let $0<\delta_0\le 1/4$ and $0<\gamma_{\Theta}<1$ be given. Let $\kappa_1 = \delta_0/10$ and $\gamma_0(m, n, \kappa_1)$ be as in Lemma \ref{l:menne sobolev 2.4}.  Let $0<\eta_0$ be as in Definition \ref{d:eta and eta 0} and $\gamma_1(m, \delta_0, 1/4)$ by as in Lemma \ref{l:upper density bound 1 scale}. Let $\gamma_2(m, \delta_0, \gamma_{\Theta})$ be as in Lemma \ref{c:apply Lip approx rect}.

Then, there is a $0< \gamma_{3}(m, \delta_0, \gamma_{\Theta})$ such the following holds.  For any $V \in \textbf{\emph{RV}}_m(B_1^n(0))$ with $\gamma_{\Theta}$-almost integral density which satisfies the hypotheses of the critical $(\delta_0, \tilde{\gamma})$-Allard tilt-excess decay regime in $B_r^n(0)$ and $\tilde{\gamma} \le \min\{\gamma_0, \gamma_1, \gamma_2, \gamma_3\}$, we may apply Lemma \ref{l:menne 3.4} in $\boldsymbol{C}(T,x,\rho, \rho)$ for any $x \in B_{\eta_0 r}^n(0)$ and any $\rho \in [r/4, (1-\eta_0)r/3]$ with constants
    \begin{align*}
        M = (2-\delta_0/4)(1+\sqrt{2})^m, &  \qquad L = 1/24\\
       \delta_1 = (1-\max\{1-\delta_0/8, (1+\gamma_{\Theta})/2\})/4, &  \qquad  \delta_2 = \max\{1-\delta_0/8, (1+\gamma_{\Theta})/2\},\\
       \delta_3=1, & \qquad \delta_4= 1/8.
    \end{align*}

In particular, for any $x \in B_{\eta_0 r}^n(0)$ and any $\rho \in [r/4, (1-\eta_0)r/3]$
\begin{align}\nonumber
    \rho^{-1-m/2}\norm{\dist(\cdot, T)}_{L^2(\norm{V}, B_\rho^n(x))} & \le 2\Gamma_{\ref{l:menne 3.4}} \tilde{\gamma},\\ \label{e:height control on Lipschitz function}
    \rho^{-1}\norm{\dist(\cdot, T)}_{L^{\infty}(\norm{V} \res \textbf{C}(T, x, \rho/2, \rho/2))} & \le C_{\ref{l:apply 3.4}}\tilde{\gamma}^{\frac{1}{m+1}},
\end{align}
where $C_{\ref{l:apply 3.4}} = \frac{2^{\frac{m+3}{m+1}}\Gamma_{\ref{l:menne 3.4}}^{\frac{1}{m+1}}}{2}$.
\end{lemma}

\begin{proof}
The choice of constants in order to apply Lemma \ref{l:menne 3.4} follows exactly as in the choice of constants in Lemma \ref{c:apply Lip approx rect}.  It remains only to show the height bounds.  Using the notation of Lemma \ref{c:apply Lip approx rect}, for $\tilde{\gamma} \le \min\{\gamma_0, \gamma_1, \gamma_2\}$ $V$ satisfies the hypotheses of Lemma \ref{l:menne 3.4}.  Thus, we immediately obtain that $G(\epsilon) = \boldsymbol{C}(T, x, \rho, \rho) \cap \text{spt}\norm{V}$. Hence, the following estimates give the first claim.
\begin{align*}
    \rho^{-1-m/2}\norm{\dist(\cdot, T)}_{L^2(\norm{V}\res \textbf{C}(T, x, \rho, \rho))} & \le 2 \rho^{-1-m/2^*}\norm{\dist(\cdot, T)}_{L^{2^*}(\norm{V}\res \textbf{C}(T, x, \rho, \rho))}\\
    & \le 2 H_{2^*}(V, a, r,h, Q, T) \le 2\Gamma_{\ref{l:menne 3.4}} T_2(V, a, 3r, h+2r, T)\\
    & \le  2\Gamma_{\ref{l:menne 3.4}} \tilde{\gamma}.
\end{align*}

To obtain the $L^{\infty}$ estimate, let $0<\kappa$ and note that if $z \in \text{spt}\norm{V}$ and $B_{\kappa \rho}^n(z) \subset \boldsymbol{C}(T, x, \rho, \rho)$ and $\dist(x, T) = \kappa \rho$, then by Lemma \ref{l:menne sobolev 2.4} we have $\norm{V}(B_{\kappa \rho/2}^n(x)) \ge (1-\delta_0/10)(\kappa \rho/2)^m$.  Furthermore, for all $y \in \text{spt}\norm{V} \cap B_{\kappa \rho/2}^n(x)$, $\dist(y, T) \ge \frac{\kappa}{2} \rho$.  Thus, we see that
\begin{align*}
    \left(\frac{\kappa}{2}\right)^{m+1}\rho^m & \le \rho^{-1-m}\norm{\dist(\cdot, T)}_{L^1(\norm{V}\res \boldsymbol{C}(T, x, \rho, \rho))}\\ & \le 2\rho^{-1-m/2}\norm{\dist(\cdot, T)}_{L^2(\norm{V}\res \boldsymbol{C}(T, x, \rho, \rho))} \le 4\Gamma_{\ref{l:menne 3.4}} \tilde{\gamma}.
\end{align*}
Thus, if $\kappa > 2^{\frac{m+3}{m+1}}\Gamma_{\ref{l:menne 3.4}}\frac{1}{m+1} \tilde{\gamma}^{\frac{1}{m+1}}$, we obtain a contradiction.  Thus, by possibly taking $\tilde{\gamma}$ even smalled if necessary, we may assume that $1-\frac{2^{\frac{m+3}{m+1}}\Gamma_{\ref{l:menne 3.4}}^{\frac{1}{m+1}}}{2}\tilde{\gamma}^{\frac{1}{m+1}} \ge 1/2$ we may obtain these estimates on $\boldsymbol{C}(T, x, \rho/2, \rho/2)$ as claimed.
\end{proof}

\begin{remark}\label{r:on the constants application}
    If $(1+\gamma_{\Theta})/2 \le 1-\delta_0/8$ then the constants $\delta_1, \delta_2$ as in Lemma \ref{c:apply Lip approx rect} and Lemma \ref{l:apply 3.4} do not depend upon $\gamma_{\Theta}$. Moreover, for $Q=1$ the resulting constant $\Gamma_{(3)}$ depends only upon $m,\delta_0$.  Following Remark \ref{r:on the constants}(6), we observe that if $0<\gamma_{\Theta}' < \gamma_{\Theta} \le 1-\delta_0/4$ the the corresponding constants satisfy $\Gamma_{(3)}(\gamma_{\Theta}')\le \Gamma_{(3)}(\gamma_{\Theta})$ and $\epsilon(\gamma_{\Theta})\le \epsilon(\gamma_{\Theta}')$.
\end{remark}

\section{Estimates on the Lipschitz approximation.}

For the remainder of this paper, we follow Allard's tilt-excess decay argument as rewritten in \cite{Simon83_book} Theorem 22.5 closely.  In this Section we show that the Lipschitz function $f$ guaranteed by Lemma \ref{c:apply Lip approx rect} and Lemma \ref{l:menne 2.18} has small Dirichlet energy and is almost weakly harmonic.

\begin{definition}\label{d:prep for estimates on lip function}
We shall need the following definitions.
\begin{enumerate}
\item (The Lipschitz function) By rotation, we may assume $T = \mathbb{R}^m \hookrightarrow \mathbb{R}^{n-m}$, embedded as the first $m$ coordinate directions.  Let $0<\eta_0$ be as in Definition \ref{d:eta and eta 0}. For $x \in B^n_{\eta_0 r}(0) \cap \text{spt}\norm{V}$, we define $f = (f^{m+1},..f^{n})$ be the $1$-valued Lipschitz approximation obtained from Lemma \ref{l:menne 2.18} applied in $\textbf{C}(T, x, r/4, r/4)$ as guaranteed by Lemma \ref{c:apply Lip approx rect}. Note that each function $f^i: \mathbb{R}^m \rightarrow \mathbb{R}$ is itself Lipschitz.
    \item (Extension of functions) For any function $g \in C(\mathbb{R}^m; \mathbb{R}^{n-m})$ we may extend $g$ to a function $\tilde{g} \in C(\mathbb{R}^n;\mathbb{R}^{n-m})$ as follows.
\begin{align}\label{e:extension}
\tilde{g}(x^1, ..., x^m, ..., x^n) := g(x^1, ..., x^m).
\end{align}
Note that if $g \in C^1_c(B^m_{r/4}(T_{\#}(x));\mathbb{R}^{n-m})$ then the height bound in Remark \ref{r:B(delta)}(2) (and \textit{a fortiori} Lemma \ref{l:apply 3.4}) implies that $|T_{\#}^{\perp}(z-x)| \le \frac{2}{5}r/4<r/8$ for $\norm{V}$-a.e. $z \in \text{spt}\norm{V} \cap \textbf{C}(T, x, r/4, r/4)$.  Thus, we may modify $\tilde{g}$ such that \eqref{e:extension} holds on $\text{spt}\norm{V} \cap \boldsymbol{C}(T, x , r/4, r/4)$ and $\tilde{g} \in C^1_c(\boldsymbol{C}(T, x, r/4, r/4);\mathbb{R}^{n-m})$.
    \item (Large pieces of the Varifold) For $x, f$ as above, we define \begin{align*}
M_1 & := \text{spt}\norm{V} \cap \text{graph}f \cap \boldsymbol{C}(T,x, r/8, r/8).
\end{align*} 
\item (Extensions of the Lipschitz function) Let $\tilde{f}^j: \mathbb{R}^n\rightarrow \mathbb{R}$ the extension of $f^j$ given by \eqref{e:extension}. Observe that for $\norm{V}$-a.e. $x = (x^1, ..., x^n) \in M_1$, $$x^{m+j} = \tilde{f}^j(x), \qquad  \nabla^M x^{m+j} = \nabla^M \tilde{f}^j.$$ 
\end{enumerate}    
\end{definition}

\begin{corollary}\label{c:estimates on lipschitz f}
    Let $\rho = r/4$. Under the hypotheses of Lemma \ref{c:apply Lip approx rect}, with $f$ as in Definition \ref{d:prep for estimates on lip function}, $0< \eta_0(m, \delta_0)$ as in Definition \ref{d:eta and eta 0}, $x \in B_{\eta_0 r}^n(0) \cap \emph{spt}\norm{V}$, and $0<\tilde{\gamma}^2 \le \gamma_3^2 \le \epsilon_1 = \epsilon$ as in Remark \ref{c:apply lipschitz approx 2}(1), then 
    \begin{align}\label{grad f bounds}
    (\rho/2)^{-m}\int_{B^m_{\rho/2}(T(x))} |\nabla f^j|^2 d\mathscr{L}^m & \le \Gamma_{\ref{c:estimates on lipschitz f}, 1} \tilde{\gamma}^2\\
    \label{f zeta bound}
    (\rho/2)^{-m}\int_{B^m_{\rho/2}(T(x))} \nabla f^j \cdot \nabla \zeta_1 d\mathscr{L}^m & \le \Big[\Gamma_{\ref{c:estimates on lipschitz f}, 2}\tilde{\gamma}^2 + \Gamma_{\ref{c:estimates on lipschitz f}, 3}\gamma_{\Theta}\tilde{\gamma}\Big]\sup |\nabla \zeta_1|,
    \end{align}
    where 
    \begin{align*}
        \Gamma_{\ref{c:estimates on lipschitz f}, 1} & = \Big[8^m(1+ L^2 + C(n)\Gamma_{(3)}\epsilon_1^{-2}L^2)\Big]\\
        \Gamma_{\ref{c:estimates on lipschitz f}, 2} & := 8^{m}C(n)\Gamma_{(3)}\epsilon_1^{-2}L +L8^m\\
        & \quad+\frac{8^{m+1}(2-\delta_0/4)}{(1-\delta_0/10)^{1/m}}\cdot C_{\text{Poincar\'e}}(n) + 8^mC(n)\Gamma_{(3)}\epsilon_1^{-2}\\
        \Gamma_{\ref{c:estimates on lipschitz f}, 3} & = (1+L^2)\Gamma_{\ref{c:estimates on lipschitz f}, 1}^{1/2}.
    \end{align*}
\end{corollary}

\textbf{Proof of \eqref{grad f bounds}.} 
Using $X_1$ as in Definition \ref{d:prep for lipschitz map}, we see
\begin{align*}
    (\rho/2)^{-m}\int_{B^m_{\rho/2}(T_{\#}(x))} |\nabla f^j|^2 d\mathscr{L}^m & \le  (\rho/2)^{-m}\int_{X_1 \cap B^m_{\rho/2}(T_{\#}(x))} |\nabla f^j|^2 d\mathscr{L}^m\\
    & + (\rho/2)^{-m}\int_{B^m_{\rho/2}(T_{\#}(x)) \setminus X_1} |\nabla f^j|^2 d\mathscr{L}^m.
\end{align*}
By \eqref{e:bad set bound} and Lemma \ref{l:menne 2.18}(3)
\begin{align*}
(\rho/2)^{-m}\int_{B^m_{\rho/2}(T_{\#}(x)) \setminus X_1} |\nabla f^j|^2 d\mathscr{L}^m & \le (\rho/2)^{-m}L^2\mathscr{L}^{m}(B^m_{\rho/2}(T_{\#}(x)) \setminus X_1)\\ 
& \le 8^{m}L^2\Gamma_{(3)}C(n)\epsilon_1^{-2}\tilde{\gamma}^2,
\end{align*}
where $L = 1/24$ as in Lemma \ref{c:apply Lip approx rect}.

To estimate $|\nabla f^j|$ on $X_1$, we use the function $F(x) = (x, f^j(x))$ and its Jacobian deternimant $J(F) = \sqrt{\text{det}(DF \circ DF^{T})}$.  Noting that $|\nabla f|^2$ is non-negative, $\Theta^m(V; x) \ge 1$ for $\norm{V}$-a.e. $x \in M_1$, and $1 \le J(F) \le 1 + L^2/2$ the Area Formula and the Triangle inequality give 
\begin{align*} \nonumber
    (\rho/2)^{-m}\int_{X_1 \cap B^m_{\rho/2}(T_{\#}(x))} |\nabla f^j|^2 d\mathscr{L}^m & \le (\rho/2)^{-m}\int_{X_1 \cap B^m_{\rho/2}(T_{\#}(x))} |\nabla f^j|^2 (\Theta^m \circ F) J(F) d\mathscr{L}^m\\ \nonumber
    & = (\rho/2)^{-m}\int_{M_1} |\nabla \tilde{f}^j|^2 d\norm{V}\\ \nonumber
& \le 2(\rho/2)^{-m}\int_{M_1} |\nabla \tilde{f}^j -\nabla^M \tilde{f}^j|^2 d\norm{V}\\
& \qquad + 2(\rho/2)^{-m}\int_{M_1} |\nabla^M \tilde{f}^j|^2 d\norm{V}.
\end{align*}
Now, we estimate each term separately. 
\begin{align}\label{e:f grad error 1}    \nonumber
    (\rho/2)^{-m}\int_{M_1} |\nabla^M \tilde{f}^j|^2 d\norm{V} & = (\rho/2)^{-m}\int_{M_1} |\nabla^M x^{m+j}|^2 d\norm{V}\\ \nonumber
    & \le (\rho/2)^{-m}\int_{M_1 \times G(n, m)} |S_{\#} - T_{\#}|^2 dV(x, S)\\
    & \le 8^{m}E(0,r,T)^2.
\end{align}
Furthermore,
\begin{align}\label{e:f grad error 2}\nonumber
    \int_{M_1} |\nabla \tilde{f}^j -\nabla^M \tilde{f}^j|^2 d\norm{V}  & = \int_{M_1 \times G(n, m)} \left|S_{\#}^\perp (\nabla \tilde{f}^j) \cdot S_{\#}^\perp(\nabla \tilde{f}^j) \right| dV(x, S)\\ \nonumber
    & = \int_{M_1 \times G(n, m)} \left|S_{\#}^\perp \circ T_{\#}(\nabla \tilde{f}^j)\right|^2 dV(x, S)\\ \nonumber
    & \le \int_{M_1 \times G(n, m)} |S_{\#}^\perp \circ T_{\#}|^2 |\nabla \tilde{f}^j|^2 dV(x, S)\\ \nonumber
    & \le \int_{M_1 \times G(n, m)} |S_{\#} - T_{\#}|^2|\nabla \tilde{f}^j|^2 dV(x, S)\\
    & \le L^2 E(0,r,T)^2 r^m.
\end{align}
Combining \eqref{e:f grad error 1}, \eqref{e:f grad error 2} gives \eqref{grad f bounds}. \qed

\textbf{Proof of \eqref{f zeta bound}.}

For any $\zeta_1 \in C^1_c(B^m_{\rho/2}(T_{\#}(x)); \mathbb{R})$ we decompose
\begin{align*}
    \left|(\rho/2)^{-m}\int_{B^m_{\rho/2}(T_{\#}(x))} \nabla f^j \cdot \nabla \zeta_1 d\mathscr{L}^m \right| & \le \left|(\rho/2)^{-m}\int_{B^m_{\rho/2}(T_{\#}(x)) \cap X_1} \nabla f^j \cdot \nabla \zeta_1 d\mathscr{L}^m \right|\\
    & \qquad + \left|(\rho/2)^{-m}\int_{B^m_{\rho/2}(T_{\#}(x)) \setminus X_1} \nabla f^j \cdot \nabla \zeta_1 d\mathscr{L}^m \right|.
\end{align*}
By Remark \ref{r:B(delta)}(1), Lemma \ref{l:menne 2.18}(3), and Lemma \ref{c:apply Lip approx rect}, we may bound the second term by as follows.
\begin{align}\label{e: weak harm error 1} \nonumber
    \left|(\rho/2)^{-m}\int_{B^m_{\rho/2}(T_{\#}(x)) \setminus X_1} \nabla f^j \cdot \nabla \zeta_1 d\mathscr{L}^m \right| & \le (\rho/2)^{-m} \mathscr{L}^m(B^m_{\rho/2}(T_{\#}(x)) \setminus X_1) L \sup |\nabla \zeta_1|\\ 
    &\le  8^{m} L \Gamma_{(3)}C(n) \epsilon_1^{-2} \tilde{\gamma}^{2} \sup |\nabla \zeta_1|.
\end{align}
Now, we compare the integral over $B^m_{\rho/2}(T_{\#}(x)) \cap X_1$ with the integral over $M_1$.  By the triangle inequality, we estimate
\begin{align} \nonumber
    & \left|(\rho/2)^{-m}\int_{B^m_{\rho/2}(T_{\#}(x)) \cap X_1} \nabla f^j \cdot \nabla \zeta_1 d\mathscr{L}^m \right| \\ \label{e: weak harmonic term 1}
    & \qquad \le \left|(\rho/2)^{-m}\int_{B^m_{\rho/2}(T_{\#}(x)) \cap X_1} (\nabla f^j \cdot \nabla \zeta_1) \Theta^m(V; F(x)) J(F)d\mathscr{L}^m\right|\\ \label{e: weak harmonic term 2}
    & \qquad \qquad + \left|(\rho/2)^{-m}\int_{B^m_{\rho/2}(T_{\#}(x)) \cap X_1} (\nabla f^j \cdot \nabla \zeta_1) (1-\Theta^m(V; F(x)) J(F)) d\mathscr{L}^m\right|.
\end{align}
Noting that $|1-(\Theta^M \circ F) J(F)| \le \gamma_{\Theta}(1 + L^2)$ for every $x \in X_1$, we may use \eqref{grad f bounds} to estimate \eqref{e: weak harmonic term 2} as follows.
\begin{align} \label{e: weak harmonic error 2} \nonumber
    & \left|(\rho/2)^{-m}\int_{B^m_{\rho/2}(T_{\#}(x)) \cap X_1} \nabla f^j \cdot \nabla \zeta_1 (1- \Theta^m(V; F(x)) J(F))d\mathscr{L}^m\right| \\ \nonumber
    & \qquad \le (\rho/2)^{-m}\int_{B^m_{\rho/2}(T_{\#}(x))} |\nabla f^j|d\mathscr{L}^m \gamma_{\Theta}(1+L^2)\sup |\nabla \zeta_1|\\ & \qquad \le (1+L^2) \Gamma_{\ref{c:estimates on lipschitz f}, 1}^{1/2} \gamma_{\Theta}\tilde{\gamma} \sup |\nabla \zeta_1|. 
\end{align}

Therefore, it only remains to estimate the term in line \eqref{e: weak harmonic term 1}. To do so, we extend $f, \zeta_1$ to functions on $\boldsymbol{C}(T, x , \rho/2, \rho/2) \subset \mathbb{R}^n$. Recalling the discussion in Definition \ref{d:prep for estimates on lip function}(3) we let $\zeta \in C^1_c(\boldsymbol{C}(T, x , \rho/2, \rho/2); \mathbb{R})$ be an extension of $\zeta_1$ such that \eqref{e:extension} holds in $\text{spt}\norm{V} \cap \boldsymbol{C}(T, x , \rho/2, \rho/2)$. Note that on this set $\nabla \tilde{f}^j(x) \cdot \nabla \zeta(x) = \nabla f^j(T_{\#}(x)) \cdot \nabla \zeta_1(T_{\#}(x))$.  Therefore, we obtain
\begin{align}\label{e:energy estimate terms 1,2} \nonumber
& \left|(\rho/2)^{-m}\int_{B^m_{\rho/2}(T_{\#}(x)) \cap X_1} (\nabla f^j \cdot \nabla \zeta_1) \Theta^m(V; F(x)) J(F) d\mathscr{L}^m\right|\\ \nonumber
& \qquad = \left|(\rho/2)^{-m}\int_{M_1} \nabla \tilde{f}^j \cdot \nabla \zeta d\norm{V} \right|\\
& \qquad \le \left|(\rho/2)^{-m}\int_{M_1} (\nabla^M \tilde{f}^j \cdot \nabla^M \zeta) d\norm{V} \right| + \left|(\rho/2)^{-m}\int_{M_1} (\nabla \tilde{f}^j \cdot \nabla \zeta - \nabla^M \tilde{f}^j \cdot \nabla^M \zeta) d\norm{V} \right|.
\end{align}
We estimate the second term as follows.
\begin{align}\label{e: weak harmonic error 3} \nonumber
    & \left|(\rho/2)^{-m}\int_{M_1} (\nabla \tilde{f}^j \cdot \nabla \zeta - \nabla^M \tilde{f}^j \cdot \nabla^M \zeta) d\norm{V} \right|\\ \nonumber
    & \qquad \le (\rho/2)^{-m}\int_{M_1} \left|\nabla \tilde{f}^j \cdot \nabla \zeta  - \nabla^M \tilde{f}^j \cdot \nabla^M \zeta \right| d\norm{V}\\ \nonumber
    & \qquad = (\rho/2)^{-m}\int_{M_1 \times G(n, m)} \left|S^\perp_{\#} (\nabla \tilde{f}^j) \cdot S^\perp_{\#}(\nabla \tilde{\zeta}_1) \right| dV(x, S)\\ \nonumber
    & \qquad = (\rho/2)^{-m}\int_{M_1 \times G(n, m)} \left|S^\perp_{\#} \circ T_{\#}(\nabla \tilde{f}^j) \cdot S^\perp_{\#} \circ T_{\#}(\nabla \tilde{\zeta}_1) \right| dV(x, S)\\ \nonumber
    & \qquad \le (\rho/2)^{-m}\int_{M_1 \times G(n, m)} |S^\perp_{\#} \circ T_{\#}|^2 |\nabla \tilde{\zeta}_1| |\nabla \tilde{f}^j| dV(x, S)\\ \nonumber
    & \qquad \le (\rho/2)^{-m}L \int_{M_1 \times G(n, m)} |S_{\#} - T_{\#}|^2  dV(x, S) \sup |\nabla \zeta_1|\\
    & \qquad \le L 8^m E(0,r,T)^2 \sup |\nabla \zeta_1|.
\end{align}

To estimate the first term in \eqref{e:energy estimate terms 1,2} we recall that for all $X \in C^1_c(B^n_{r}(0); \mathbb{R}^n)$
\begin{align*}
    \int  DX \cdot S_{\#} dV(x, S) = -\int X\cdot \boldsymbol{h}(V; x) d\norm{V}(x).
\end{align*}
Thus, for any $\zeta \in C^1_c(B^n_{r}(0); \mathbb{R})$ and any $j=1, ..., n-m$ we may let $X= \zeta e_{m+j}$.  Observing that $D(\zeta e_{m+j}) \cdot S_{\#} = (\nabla^M x^{m+j})\cdot \nabla^M \zeta$ we obtain
\begin{align*}
    \int (\nabla^M x^{m+j})\cdot \nabla^M\zeta d\norm{V} = - \int e_{m+j}\cdot \boldsymbol{h}(V; x) \zeta d\norm{V}(x).
\end{align*}
Thus, for $\zeta \in C^1_c(\boldsymbol{C}(T, x, \rho/2, \rho/2))$ the extension of $\zeta_1$ and recalling $\nabla^M \tilde{f}^j = \nabla^M x^{m+j}$ on $M_1$ we estimate by H\"older's inequality and Poincar\'e's inequality
\begin{align}\label{e: weak harmonic error 4}\nonumber
    (\rho/2)^{-m}&\int_{M_1} (\nabla^M \tilde{f}^{j})\cdot \nabla^M\zeta d\norm{V}\\ \nonumber
    & = - (\rho/2)^{-m}\int e_{m+j}\cdot \boldsymbol{h}(V; x) \zeta d\norm{V}\\ \nonumber & \qquad \qquad - (\rho/2)^{-m}\int_{\text{\spt}\norm{V} \setminus M_1} (\nabla^M x^{m+j})\cdot \nabla^M\zeta d\norm{V}\\ \nonumber
    & \le (\rho/2)^{-m}\left(\int_{\text{spt}\zeta}|\boldsymbol{h}(V;x)|d\norm{V} \sup |\zeta| + \norm{V}(\boldsymbol{C}(T, x, \rho/2, \rho/2) \setminus M_1) \sup |\nabla \zeta| \right)\\ \nonumber
    & \le \gamma_{\hbar,m}(0,r)\norm{V}(B_{\rho}^n(x))^{1-1/m}(\rho/2)^{-m}\sup|\zeta_1|\\ \nonumber
    & \qquad \qquad + 8^m\Gamma_{(3)}C(n)\epsilon_1^{-2}\tilde{\gamma}^2 \sup|\nabla \zeta_1|\\
    & \le \left[\frac{8^{m+1}(2-\delta_0/4)}{(1-\delta_0/10)^{1/m}}\cdot C_{\text{Poincar\'e}}(n) \cdot \gamma_{\hbar,m}(0,r) + 8^mC(n)\Gamma_{(3)}\epsilon_1^{-2}\tilde{\gamma}^2\right]\sup|\nabla \zeta_1|.
\end{align}

Combining \eqref{e: weak harm error 1}, \eqref{e: weak harmonic error 2}, \eqref{e: weak harmonic error 3}, and \eqref{e: weak harmonic error 4} under the assumption that $E(0,r,T) \le \tilde{\gamma}$ and $\gamma_{\hbar,m}(0,r) \le \tilde{\gamma}^2$ we obtain \eqref{f zeta bound}. \qed

\begin{remark}\label{r:on the constants for lipschitz function}
    Following Remark \ref{r:on the constants} and Remark \ref{r:on the constants application}, if $\gamma_{\Theta} \le 1-\delta_0/4$ then $\Gamma_{(3)}$ may be taken to be independent of $\gamma_{\Theta}$.  Moreover, for $0<\gamma_{\Theta}' < \gamma_{\Theta} \le 1-\delta_0/4$ then the corresponding $\epsilon_1(\gamma_{\Theta})\le \epsilon_1(\gamma_{\Theta}')$, which implies $\Gamma_{\ref{c:estimates on lipschitz f}, i}(\gamma_{\Theta}') \le \Gamma_{\ref{c:estimates on lipschitz f}, i}(\gamma_{\Theta})$ for $i=1,2,3$.  This implies that for $\gamma_{\Theta} \le 1-\delta_0/4$ then the constants $\epsilon_1, \Gamma_{\ref{c:estimates on lipschitz f}, i}$ may be taken independent of $\gamma_{\Theta}$.
\end{remark}

\section{The tilt-excess decay lemma}\label{s:tilt-excess decay}

The previous lemmata allow us to follow Allard's tilt-excess decay lemma. The main novelty below is the observation that the density parameter $\gamma_{\Theta}$ is only used to obtain a Lipschitz approximation which is sufficiently almost-weakly harmonic in the sense of Lemma \ref{l:harmonic approximation}.  But, all the estimates on the harmonic function guaranteed by Lemma \ref{l:harmonic approximation} may be taken in terms of $\tilde{\gamma}$.  This is what allows us to iterate the argument without the density needing to be closer and closer to integral.

\begin{lemma}\label{l:7 rect}\emph{(The tilt-excess decay lemma, cf. \cite{Allard72} 8.16 or \cite{Simon83_book} 22.5)}
Let $0<\delta_0\le 1/4$ be given.  Then, there exists a scale $0< \eta_1(n, m, \delta_0) \le \eta_0(m,\delta_0)/4$ and constants
\begin{align*}
    0<\gamma_{\Theta}(m,\delta_0)\le 1/2, \quad 0< \gamma_4(n, m, \delta_0),
\end{align*}
such that the following statement holds. Let $\gamma_0(m,n,\delta_0/10)$ as in Lemma \ref{l:menne sobolev 2.4}.  Let $0< \gamma_1(m, \delta_0, \eta_1)$ as in Lemma \ref{l:upper density bound 1 scale}, $0<\gamma_{2}(m, n, \delta_0, \gamma_{\Theta})$ as in Lemma \ref{c:apply Lip approx rect}, and $0<\gamma_3(m,\delta_0, \gamma_{\Theta})$ be as in Lemma \ref{l:apply 3.4}.

If $V$ satisfies the hypotheses of the critical $(\delta_0, \tilde{\gamma})$-Allard tilt-excess decay regime in $B_r^n(0)$ with constant $0< \tilde{\gamma} \le \min\{\gamma_0, \gamma_1, \gamma_2, \gamma_3, \gamma_4 \}$ and $\gamma_{\Theta}$-almost integral density, then there exists an $R_x \in G(n, m)$ such that
\begin{align}
    \left((\eta_1 r)^{-m} \int_{B^n_{\eta_1 r}(x) \times G(n, m)}|S_\# - (R_x)_{\#}|^2dV(\xi, S)\right)^{\frac{1}{2}} & \le \frac{1}{2}\tilde{\gamma}.
\end{align}
\end{lemma}

\begin{proof}
First, by rotation we assume that $T = \mathbb{R}^m \hookrightarrow \mathbb{R}^n$ as the first $m$ coordinates. Note that by Remark \ref{r:on the constants for lipschitz function}, since $\gamma_{\Theta} \le 1/2 < 1-\delta_0/4$ the constants $\epsilon_1, \Gamma_{\ref{c:estimates on lipschitz f}, i}$ are independent of $\gamma_{\Theta}$.

To obtain $R \in G(n, m)$, we invoke Lemma \ref{c:apply Lip approx rect} to apply Lemma \ref{l:menne 2.18} in $\textbf{C}(T, x, r/4, r/4)$.  Thus, the estimates \eqref{f zeta bound} and \eqref{grad f bounds} from Corollary \ref{c:estimates on lipschitz f} apply and we may use Lemma \ref{l:harmonic approximation} to $[\Gamma_{\ref{c:estimates on lipschitz f},1}^{1/2}\tilde{\gamma}]^{-1}f^j$ in $B^m_{r/8}(x)$.  That is, for any $0<\kappa_3$ there exists an $0<\epsilon(n, m, \delta_0, \kappa_3)$ such that if \begin{align*}
\frac{\Gamma_{\ref{c:estimates on lipschitz f}, 2}}{\Gamma_{\ref{c:estimates on lipschitz f}, 1}^{1/2}}\tilde{\gamma} \le \epsilon(n, m, \delta_0, \kappa_3)/2, \qquad (1+L^2)\gamma_{\Theta} \le \epsilon(n, m, \delta_0, \kappa_3)/2
\end{align*}
then for $j = m+1, ...,n$ there exists harmonic functions $u^j: B^m_{r/8}(T_{\#}(x)) \rightarrow \mathbb{R}$ such that
\begin{align}\label{energy constraints on harmonic function}
    (r/8)^{-m}\int_{B^m_{r/8}(T_{\#}(x))}|\nabla u^j|^2d\mathscr{L}^m & \le \Gamma_{\ref{c:estimates on lipschitz f}, 1} \tilde{\gamma}^2\\ \label{u^j close to f^j}
    (r/8)^{-m-2}\int_{B^m_{r/8}(T_{\#}(x))}|f^j-u^j|^2d\mathscr{L}^m & \le \Gamma_{\ref{c:estimates on lipschitz f}, 1} \kappa_3 \tilde{\gamma}^2.
\end{align}

Recall that under the hypotheses of Lemma \ref{l:apply 3.4}, 
$\sup_{B^m_{r/8}(T_{\#}(x))}|f| \le C_{\ref{l:apply 3.4}}(n,m,\delta_0) \tilde{\gamma}^{\frac{1}{m+1}}(r/8)$.  We claim that 
\begin{align}\label{e:u small claim}|u^j(0)| \le 3 C_{\ref{l:apply 3.4}} \tilde{\gamma}^{\frac{1}{m+1}} r/8.
\end{align}
To see the claim, we use \eqref{energy constraints on harmonic function} to estimate
\begin{align*}
    \sup_{y \in B^m_{r/8}(T(x))}|u^j(y) -u^j(T_{\#}(x))| & \le (r/8)\sup_{y \in B^m_{r/8}(T_{\#}(x))}|\nabla u(y)| \\
    & \le (r/8) 2^m \Gamma_{\eqref{c:estimates on lipschitz f}, 1}^{1/2}\tilde{\gamma}.
\end{align*}
Thus, for $0<\tilde{\gamma}$ sufficiently small such that $$2^m \Gamma_{\ref{c:estimates on lipschitz f}, 1}^{1/2}\tilde{\gamma} < 3C_{\ref{l:apply 3.4}} \tilde{\gamma}^{\frac{1}{m+1}},$$ then $|u(y) - f(y)| \ge \frac{C_{\ref{l:apply 3.4}} \tilde{\gamma}^{\frac{1}{m+1}}}{2} r/8$ for all $y \in B^m_{r/8}(T_{\#}(x))$. However, by \eqref{u^j close to f^j} this would imply that for $\tilde{\gamma}$ sufficiently small
\begin{align}\label{e:u(0) small}
    2^{2-m}C_{\ref{l:apply 3.4}} \tilde{\gamma}^{\frac{1}{m+1}} \le \Gamma_{\ref{c:estimates on lipschitz f}, 1}^{1/2} \kappa_3 \tilde{\gamma}.
\end{align}
Taking $\tilde{\gamma}$ sufficiently small depending upon $n,m,\delta_0, \kappa_3$ this is a contradiction.  Thus, \eqref{e:u small claim} holds.

Therefore, if we let $\ell^j(z) = u^j(T_{\#}(x)) + \nabla u^j(T_{\#}(x))\cdot(z-T_{\#}(x))$ and $\ell(z)=(\ell^1(z), ..., \ell^{n-m}(z))$ then
\begin{align*}
    |\ell(T_{\#}(x))| \le \sum_{j=1}^{n-m}|\ell^j(T_{\#}(x))| \le 3(n-m) C_{\ref{l:apply 3.4}} \tilde{\gamma}^{\frac{1}{m+1}}r/8.
\end{align*}
Let $\tau= (T_{\#}(x), \ell(T_{\#}(x)))$.  Let $R_x \in G(n, m)$ be such that $R_x = \text{Image}(\ell)$.

Now, we focus on obtaining the scale $\eta_1(m,n,\delta_0)$.  For $0< \eta \le 1/2$ we let $\rho = \eta r/8$. Let $0< \kappa_4\le 1$. By taking $0< \gamma_4(n, m, \delta_0, \kappa_3, \kappa_4, \eta)$ sufficiently small, we may assume that the following conditions hold in $B^n_{\rho}(x).$
\begin{enumerate}
\item (Height Control) By Lemma \ref{l:apply 3.4} applied in $\boldsymbol{C}(T, x, r/4, r/4)$ we may assume that 
\begin{align*}
    \text{spt}\norm{V} \cap B^n_{\rho}(x) \subset B^n_{2C_{\ref{l:apply 3.4}} \tilde{\gamma}^{\frac{1}{m+1}} r/8}(T)
\end{align*}
and $|\tau| \le 12(n-m)C_{\ref{l:apply 3.4}} \tilde{\gamma}^{\frac{1}{m+1}} r/8$.  By taking $0<\gamma_4$ small we may assume that $$[12(n-m)+4]C_{\ref{l:apply 3.4}} \gamma_4^{\frac{1}{m+1}} r/8 \le \kappa_4 \rho.$$ 

\item (Density control) We may also assume that $\gamma_4(m, \delta_0, \kappa_3, \kappa_4, \eta) \le \gamma_2(m, \delta_0, \eta)$ as in Lemma \ref{l:upper density bound 1 scale}.
\item (Tilt Control) From \eqref{energy constraints on harmonic function} and Jensen's inequality, we may assume that 
\begin{align*}
    \norm{R_x-T} = |Du(T_{\#}(x))| \le \Gamma_{\ref{c:estimates on lipschitz f}, 1}^{1/2}\gamma_4 \le 2C_{\ref{l:apply 3.4}} \tilde{\gamma}^{\frac{1}{m+1}}.
\end{align*}
\end{enumerate}
Note that under these assumptions, for all $z \in \text{spt}\norm{V} \cap B^n_{\rho}(x)$ 
\begin{align}\label{e:crude height bound k4}\nonumber
    \dist(z-\tau, R_x) & \le |\tau| + \norm{R_x - T} \cdot \rho + 2C_{\ref{l:apply 3.4}}\tilde{\gamma}^{\frac{1}{m+1}} r/8\\
    & \le  (12(n-m)+4)C_{\ref{l:apply 3.4}} \tilde{\gamma}^{\frac{1}{m+1}} r/8 \le \kappa_4 \rho.
\end{align}

With these assumptions, for $m \ge 2$ we apply Lemma \ref{l:Simon 22.2} and we obtain
\begin{align*}
    &\rho^{-m}\int_{B^n_{\rho}(\tau) \times G(n,m)}|S_{\#}-(R_x)_{\#}|^2dV(x, S) \\
    & \quad \le C \left[(2\rho)^{-m}\int_{B^n_{2\rho}(\tau)}\left|\frac{\text{dist}(x-\tau, R_x)}{2\rho}\right|^2d\norm{V}(x) + (2\rho)^{2-m}\int_{B^n_{2\rho}(\tau)}|\boldsymbol{h}(V; x)|^2d\norm{V}(x)\right]\\
    & \quad \le C \left[(2\rho)^{-m}\int_{B^n_{2\rho}(\tau)}\left|\frac{\text{dist}(x-\tau, R_x)}{2\rho}\right|^2d\norm{V}(x) + (2\rho)^{2-m}\gamma_{\hbar,m}(0,r)^2\norm{V}(B^n_{2\rho}(\tau))^{1-2/m}\right].
    \end{align*}
Now, we decompose the first term and estimate by Lemma \ref{l:menne 2.18}, Remark \ref{r:B(delta)}, and the same argument by which we obtained \eqref{e: weak harmonic error 2}.
\begin{align*}
    & (2\rho)^{-m-2}\int_{B^n_{2\rho}(\tau)}|\text{dist}(z-\tau, R_x)|^2d\norm{V}(z)\\\nonumber
    & \qquad \le (2\rho)^{-m-2}\int_{B^n_{2\rho}(\tau) \cap M_1}|\text{dist}(z-\tau, R_x)|^2d\norm{V}(z)\\ \nonumber 
    & \qquad \qquad + (2\rho)^{-m-2}\int_{B^n_{2\rho}(\tau) \setminus M_1}|\text{dist}(z-\tau, R_x)|^2d\norm{V}(z)\\ \nonumber
    & \qquad \le (2\rho)^{-m-2}\int_{B^m_{2\rho}(T(x)) \cap X_1}|\text{dist}(F(z)-\tau, R_x)|^2 \Theta^m(V; z)J(F)d\mathscr{L}^m(z)\\ \nonumber 
    & \qquad \qquad + \kappa_4^2(2\rho)^{-m}\norm{V}(\boldsymbol{C}(T, 0, \eta_0 r, \eta_0 r) \setminus M_1).
\end{align*}
Note that by Remark \ref{r:B(delta)}(1), we have 
\begin{align}\label{e:tilt excess error 0}
    \kappa_4^2(2\rho)^{-m}\norm{V}(\boldsymbol{C}(T, x, \eta_0 r, \eta_0 r) \setminus M_1) & \le \kappa_4^2\eta^{-m}\Gamma_{(3)}C(n)\epsilon_1^{-2}\tilde{\gamma}^2.
\end{align}
Now, we estimate the remaining term. Note that for any $z \in X_1$, $\dist(F(z)-\tau, R) \le (1+L^2)^{1/2}|f(z)-\ell(z)|$, $|J(F(z))| \le (1+L^2)^{1/2}$, and  $(1+\gamma_\Theta)(1+L^2)<3$. Moreover, by the triangle inequality, Lemma \ref{l: harmonic D2 to D1 estimate} with $R = r/4$, \eqref{energy constraints on harmonic function}, and \eqref{u^j close to f^j} we estimate for $m>1$ any $j = 1, ..., n-m$
\begin{align}\label{e: tilt excess error 2} \nonumber
  &(2\rho)^{-m-2}\int_{B^m_{2\rho}(T_{\#}(x)) \cap X_1}|\text{dist}(F(z)-\tau, R)|^2 \Theta^m(V; z)J(F)d\mathscr{L}^m(z)\\ \nonumber
    & \qquad \le 3(2\rho)^{-m-2}\int_{B^m_{2\rho}(T_{\#}(x))}|f(z)-\ell(z)|^2 d\mathscr{L}^m(z)\\ \nonumber
    & \qquad \le 6(2\rho)^{-m-2}\int_{B^m_{2\rho}(T_{\#}(x))}|u-\ell|^2d\mathscr{L}^m + 6\sum_{j=m+1}^{n}(2\rho)^{-m-2}\int_{B^m_{2\rho}(T_{\#}(x))}|f^j-u^j|^2d\mathscr{L}^m\\ 
    & \qquad \le  3C(m)8^{m}\eta^2\Gamma_{\ref{c:estimates on lipschitz f}, 1}\tilde{\gamma}^2 + 6(n-m)\eta^{-m-2}\Gamma_{\ref{c:estimates on lipschitz f}, 1} \kappa_3 \tilde{\gamma}^2.
\end{align}

Now, we turn to estimate the term $(2\rho)^{2-m}\gamma_{\hbar,m}(0,r)^2\norm{V}(B^n_{2\rho}(\tau))^{1-2/m}$. For our choices of $\gamma_4(m, \delta_0, \kappa_3, \kappa_4, \eta)$, we see that by Lemma \ref{l:menne sobolev 2.4} and Lemma \ref{l:upper density bound 1 scale}
\begin{align*}
    \norm{V}(B^n_{2\rho}(\tau))^{1-2/m} & \le \frac{\norm{V}(B^n_{2\rho}(\tau))}{\norm{V}(B^n_{\rho}(x))^{2/m}}\\
    & \le \frac{(2-\delta_0/4)\omega_m (2\rho)^m}{[(1-\delta_0/10)\omega_m\rho^m]^{2/m}} = 2^m\frac{2-\delta_0/4}{1-\delta_0/10}\omega_m^{1-2/m}\rho^{m-2}.
\end{align*}
Hence, trivially we may estimate
\begin{align}\label{e: tilt excess error 3}
    (2\rho)^{2-m}\gamma_{\hbar,m}(0,r)^2\norm{V}(B^n_{2\rho}(\tau))^{1-2/m} & \le 4\frac{2-\delta_0/4}{1-\delta_0/10}\omega_m^{1-2/m}\tilde{\gamma}^4. 
\end{align}

Now, we may assume $0<\kappa_4<1/2$ so that our height bound implies $B^n_{\rho/2}(x) \subset B_{\rho}^n(\tau)$.  Thus, putting \eqref{e:tilt excess error 0}, \eqref{e: tilt excess error 2}, and \eqref{e: tilt excess error 3} together we obtain
\begin{align}\label{e: tilt excess 2}\nonumber
    &(\rho/2)^{-m}\int_{B^n_{\rho/2}(x) \times G(n, m)}|S_{\#}-(R_x)_{\#}|^2dV(x, S)\\ \nonumber
    & \qquad \qquad \le 2^m\rho^{-m}\int_{B^n_{\rho}(\tau) \times G(n, m)}|S_{\#}-(R_{x})_{\#}|^2dV(x, S)\\
    & \qquad \qquad \le C(n,m,\delta_0) [\kappa_4^2\eta^{-m} + \eta^2 + \eta^{-m-2} \kappa_3 + \tilde{\gamma}^2]\tilde{\gamma}^2,
\end{align}
for $C(n,m, \delta_0) = 2^m\max\{\Gamma_{(3)}C(n)\epsilon_1^{-2}, 3C(m)8^{m}\Gamma_{\ref{c:estimates on lipschitz f}, 1}, 6(n-m)\Gamma_{\ref{c:estimates on lipschitz f}, 1}, 4\frac{2-\delta_0/4}{1-\delta_0/10}\omega_m^{1-2/m}\}$. 

Now we choose our constants as follows.
\begin{enumerate}
    \item Choose $0<\eta$ such that $C(n, m, \delta_0) \eta^2 < \frac{1}{20}$.
    \item Choose $0< \kappa_3$ such that $C(n,m,\delta_0)\eta^{-m-2}\kappa_3<\frac{1}{20}.$
    \item Choose $0<\kappa_4\le 1$ such that $C(n,m,\delta_0)\kappa_4^2\eta^{-m}<1/20$.
    \item If necessary, we may additionally ensure that $C(n, m, \delta_0)\tilde{\gamma}^2 < \frac{1}{20}$.
\end{enumerate}
Let $0<\gamma_4(m, \delta_0, \kappa_3, \kappa_4, \eta) = \gamma_4(n,m,\delta_0)$  and $\gamma_{\Theta}(n,m, \delta_0) =\min\{\gamma_{\Theta}(n,m, \delta_0, \kappa_3), 1/2\} $ corresponding to these choices.  Then, for
\begin{align*}
    0< \tilde{\gamma} \le \min\{\gamma_0(n,m,\delta_0), \gamma_1(m, \delta_0), \gamma_2(n, m, \delta_0), \gamma_3(n, m,\delta_0), \gamma_4(n,m,\delta_0)\},
\end{align*} we have
\begin{align}\label{e: tilt excess 3}
    &(\frac{\rho}{2})^{-m}\int_{B^n_{\rho/2}(x) \times G(n, m)}|S_{\#}-(R_x)_{\#}|^2dV(x, S) \le \frac{1}{4}\tilde{\gamma}^2.
\end{align}
Therefore, choosing $\eta_1 := \eta(n, m, \delta_0)/16$ proves the lemma in the $m\ge 2$ case.

For $m=1$, recall that we are assuming $\gamma_{\hbar, 1} \le \tilde{\gamma}^4$.  Now, we invoke Lemma \ref{l:brakke 5.5 m=1} with $\phi \in C^1_0(B_{2\rho}^n(\tau))$ such that $\phi =1$ on $B_{\rho}^n(\tau)$ and $|\nabla \phi| \le 2\rho^{-1}$. This gives
\begin{align}\label{e:m=1 apply cacciopoli}\nonumber
        E(\tau, \rho, R_x)^2 & \le 4m\norm{\boldsymbol{h}(V, \cdot)}_{L^1(\norm{V}\res B_{2\rho}^n(\tau))}^{1/3} (2\rho)^{\frac{-m-2}{3}}\norm{\dist(\cdot, R_x+\tau)}_{L^2(\norm{V}\res B_{2\rho}^n(\tau))}^{2/3}\\
        & \qquad \qquad + 64 (2\rho)^{-m-2}\norm{\dist(\cdot, R_x+\tau)}_{L^2(\norm{V}\res B_{2\rho}^n(\tau))}.
    \end{align}
Because $m=1$, $u^j = \ell^j$ for each $j=1, ..., n-m$.  Therefore, splitting the height integral into the intersection of $M_1$ and the complement we obtain the following estimates.  First, we obtain by exactly the same argument as in \eqref{e: tilt excess error 2}
\begin{align*}
& \rho^{-m}\int_{B_{2\rho}^n(\tau) \cap M_1} \frac{|\dist(x-\tau,R_x)|^2}{\rho^2} d\norm{V}(x)\\
& \qquad \le 3\rho^{-m}\int_{B_{2\rho}^m(T_{\#}(\tau)) \cap X_1} \frac{|f(x)- u(x)|^2}{\rho^2} d\mathscr{L}^m(x)\\
& \qquad \le 6(n-m)\eta^{-m-2}\Gamma_{\ref{c:estimates on lipschitz f}, 1} \kappa_3 \tilde{\gamma}^2.
\end{align*}
Secondly, by exactly the same argument as in \eqref{e:tilt excess error 0}, we obtain
\begin{align*}
(2\rho)^{-m-2}\int_{B^n_{2\rho}(\tau) \setminus M_1}|\text{dist}(z-\tau, R_x)|^2d\norm{V}(z) \le \kappa_4^2\eta^{-m}\Gamma_{(3)}C(n)\epsilon_1^{-2}\tilde{\gamma}^2.
\end{align*}
Plugging these into \eqref{e:m=1 apply cacciopoli}, we obtain the following.
\begin{align*}
    E(\tau, \rho, R_x)^2 & \le C(n,m,\delta_0)\left(\tilde{\gamma}^{4/3}(\kappa_3 \eta^{-m-2} + \kappa_4^2 \eta^{-m})^{1/3}\tilde{\gamma}^{2/3} + (\kappa_3 \eta^{-m-2} + \kappa_4^2 \eta^{-m})\tilde{\gamma}^{2}\right).
\end{align*}
Therefore, choosing $\eta = 1/2$ and $0<\kappa_3, \kappa_4$ such that $C(n,m,\delta_0)(\kappa_3 2^{-3} + \kappa_4^2 2^{-1})^{1/3} \le 1/16$, we obtain thresholds $0<\gamma_4(n,m,\delta_0), \gamma_{\Theta}(n,m,\delta_0)$ as before such that \eqref{e: tilt excess 3} holds for $m=1$.

\end{proof}

\section{Proof of the Main Results}\label{s:proof of main results}

The proofs of Theorem \ref{t:BMO main theorem rect} and Corollary \ref{t: VMO main theorem rect}, follow from inductively applying Lemma \ref{l:7 rect}.  Before we can do so, we need to show that under the assumptions of Lemma \ref{l:7 rect}, $V$ satisfies the density hypotheses of the $(\delta_0, \tilde{\gamma})$-Allard tilt-excess decay regime in $B_{\eta_1 r}^n(x)$.

\begin{lemma}\label{l:improved density bounds}\emph{(Improved density bounds)}
    Let $V, \eta_1, x, \tilde{\gamma}, \gamma_{\Theta}$ be as in Lemma \ref{l:7 rect}.  Then for all $\rho \in [\eta_1 r, (1-\epsilon_0)r/6]$
    \begin{align*}
        & \max\{1-\delta_0/10, \Big[(1-C_{\ref{l:apply 3.4}}\tilde{\gamma}^{\frac{1}{m+1}})^m - \Gamma_{(3)}C(n)\epsilon_1^{-2}\eta_1^{-m}\tilde{\gamma}^2\Big]\}\\
        & \qquad \le \frac{\norm{V}(B_{\rho}^n(x))}{\omega_m\rho^m}\\
        & \qquad \le (1+\gamma_{\Theta})(1+L^2)^{1/2}+\Gamma_{(3)}C(n)\epsilon_1^{-2}2^m \tilde{\gamma}^2.
    \end{align*} 
    In particular, $\gamma_{\Theta}\le 1/2$ ensures that
    \begin{align*}
        \delta_0 \le \frac{\norm{V}(B^n_{\eta_1r}(x))}{\omega_mr^m}\le 2-\delta_0.
    \end{align*}
\end{lemma}
\begin{proof}
    We note that the $\gamma_4$ in Lemma \ref{l:7 rect} is chosen such that by Lemma \ref{c:apply Lip approx rect} and Lemma \ref{l:apply 3.4} we may apply Lemma \ref{l:menne 2.18} in $\textbf{C}(T, x, 2\rho, 2\rho)$.  By the height bound in \eqref{e:height control on Lipschitz function} and the estimates in Lemma \ref{l:menne 2.18}(3) and Remark \ref{r:B(delta)}(1) we estimate
    \begin{align*}
        \norm{V}(B_\rho^n(x)) & \ge \norm{V}(B_\rho^n(x) \cap \text{graph}f)\\
        & \ge \norm{V}(\text{graph}f|_{B^n_{(1-C\tilde{\gamma}^{\frac{1}{m+1}})\rho}(T_{\#}(x))})\\
        & \ge (1-C_{\ref{l:apply 3.4}}\tilde{\gamma}^{\frac{1}{m+1}})^m\omega_m\rho^m - \mathcal{L}^m(C(\epsilon_1))\\
        & \ge (1-C_{\ref{l:apply 3.4}}\tilde{\gamma}^{\frac{1}{m+1}})^m\omega_m\rho^m - \Gamma_{(3)}C(n)\epsilon_1^{-2}\tilde{\gamma}^2r^m\\
        & \ge \Big[(1-C_{\ref{l:apply 3.4}}\tilde{\gamma}^{\frac{1}{m+1}})^m - \Gamma_{(3)}C(n)\epsilon_1^{-2}\eta_1^{-m}\tilde{\gamma}^2\Big]\omega_m\rho^m.
    \end{align*}
The other lower bound comes from Lemma \ref{l:menne sobolev 2.4}.

To obtain the upper bound, we again use the Lipschitz approximation and the estimates in Lemma \ref{l:menne 2.18}(3) and Remark \ref{r:B(delta)}(1). 
\begin{align*}
    \norm{V}(B_{\rho}^n(x)) & \le \norm{V}(\textbf{C}(T, x, \rho, \rho) \cap \text{graph}f) + \norm{V}(\textbf{C}(T, x, \rho, \rho) \setminus \text{graph}f)\\
        & \le  (1+\gamma_{\Theta})(1+L^2)^{1/2}\omega_m\rho^m + \norm{V}(D(\epsilon_1))\\
        & \le [(1+\gamma_{\Theta})(1+L^2)^{1/2}+\Gamma_{(3)}C(n)\epsilon_1^{-2}2^m \tilde{\gamma}^2]\omega_m\rho^m.
\end{align*}

Now, we note that from the proof of Lemma \ref{l:7 rect} (the choice of constants (4) and the definition of $C(n,m,\delta_0)$) we see that
\begin{align*}
\Gamma_{(3)}C(n)\epsilon_1^{-2}2^m \tilde{\gamma}^2 < \frac{1}{20}.   
\end{align*}
Whence, we may calculate that for $\gamma_{\Theta} \le 1/2$
\begin{align*}
    (1+\gamma_{\Theta})(1+24^{-2})^{1/2}+\frac{1}{20} < 2-\frac{3}{4} \le 2-\delta_0.
\end{align*}
\end{proof}

\subsection{Proof of Theorem \ref{t:BMO main theorem rect}}\label{s:proof of BMO}
If $0 \le \tilde{\gamma} \le \{\gamma_0,\gamma_1,\gamma_2,\gamma_3,\gamma_4\}$ and $0<\gamma_{\Theta}(n,m,\delta_0)$ as in Lemma \ref{l:7 rect} and Lemma \ref{l:improved density bounds} then for all $x \in \text{spt}\norm{V} \cap B^n_{\eta_0 r}(0)$, the varifold $V$ had $\gamma_{\Theta}$-almost integral density and satisfies the hypotheses of the $(\delta_0, \tilde{\gamma})$-Allard tilt-excess decay regime in $B^n_{\eta_1 r}(x)$.

Therefore, we may reapply Lemma \ref{l:7 rect} and Lemma \ref{l:improved density bounds} in $B^n_{\eta_1 r}(x)$ relative to $R_x \in G(n,m)$.  More generally, for $i \in \mathbb{N} \cup \{0\}$ define scales $r_i = \eta_1^i r$ and $R_{x, 0}:=T, R_{x, 1} = R_x$, as above.  We inductively define $R_{x,i} \in G(n,m)$ to be the plane guaranteed by applying Lemma 
\ref{l:7 rect} in $B^n_{r_{i-1}}(x)$ relative to $R_{x,i} \in G(n,m)$.  Thus, for any $0< \rho< r$, if we let $j \in \mathbb{N} \cup \{0\}$ be such that $r_{j+1} < \rho \le r_{j}$ then
\begin{align}\label{e:tilt excess at all scales} \nonumber
    \inf_{R \in G(n,m)} E(x, \rho, R) & \le E(x, \rho, R_{x, j})\\
    & \le \eta_1^{-m/2} E(x, r_{j}, R_{x,j}) \le C(n,m,\delta_0) \tilde{\gamma}.
\end{align} 
\qed

\subsection{Proof of Corollary \ref{t: VMO main theorem rect}}
We use the notation of the proof of Theorem \ref{t:BMO main theorem rect}, above.  Let $x \in \text{spt}\norm{V} \cap B^n_{\eta_0 r}(0)$ and $r_i = \eta_1^i r$. Under the hypotheses of Theorem \ref{t:BMO main theorem rect}, as noted above, Lemma \ref{l:7 rect} and Lemma \ref{l:improved density bounds} imply that $V$ satisfies the hypotheses of the $(\delta_0, \tilde{\gamma})$-Allard tilt-excess decay regime in $B^n_{r_1}(x)$ with constant $\tilde{\gamma}(x, r_1) := \max\{\frac{1}{2}\tilde{\gamma}, \gamma_{\hbar,m}(x, r_1)^{1/2}\}$. By induction, we see that for $i \ge 2$
\begin{align*}    
\tilde{\gamma}(x, r_i) := \max\{\frac{1}{2}\tilde{\gamma}(r_{i-1}), \gamma_{\hbar,m}(x, r_i)^{1/2}\}.
\end{align*}
Since $\boldsymbol{h}(V; \cdot) \in L^1_{loc}(\norm{V}; \mathbb{R}^n)$ and $\lim_{i \rightarrow \infty} \sup_{x \in B_{\eta_0 r}^n(0)}\gamma_{\hbar,m}(x, r_i) = 0$ for all $x \in \text{spt}\norm{V} \cap B^n_{\eta_0 r}(0)$.  Whence, $\lim_{i \rightarrow \infty}\sup_{x \in B_{\eta_0 r}^n(0)}\tilde{\gamma}(x, r_i) = 0$.  Therefore the estimate \eqref{e:tilt excess at all scales} gives the conclusion of the corollary. \qed

\subsection{Proof of Porism \ref{t:decay rates}}
The proof of Porism \ref{t:decay rates} follows from considering the proof of Corollary \ref{c:estimates on lipschitz f} and Lemma \ref{l:7 rect}. First, we assume $m \ge 2$ and $0<\alpha<1$.  Let $A \subset B_{\eta_0 r}^n(0) \cap \text{spt}\norm{V}$ be the set of points which are both Lebesgue points of $\boldsymbol{h}(V; \mathbb{R}^n) \in L^1(\norm{V}; \mathbb{R}^n)$ and Lebesgue point of $\Theta^m(\norm{V}, \cdot)$.  Note that since $\norm{V}$ is Radon and $\boldsymbol{h}(V, \cdot),\Theta^m(\norm{V}, \cdot)$ are $\norm{V}$-measurable, $\norm{V}(B_{\eta_0 r}^n(0) \setminus A) = 0$.  

Let $x \in A$.  By Corollary \ref{t: VMO main theorem rect} we may choose a scale $0< \rho_0 \le r$ depending upon the point $x$ such that the following conditions hold. 
\begin{enumerate}
    \item (Density point of $\boldsymbol{h}(V, \cdot)$) There exists a $0<\Gamma< \infty$ for which 
\begin{align*}
\gamma_{\hbar,m}(x, \rho) \le \Gamma \rho^m \le \left(\min\{\gamma_0, \gamma_1, \gamma_2, \gamma_3, \gamma_4\}\right)^2
\end{align*}
for all $0<\rho \le \rho_0$.
\item For $C(n,m,\delta_0)$ as in \eqref{e: tilt excess 2}, we may choose $\eta_2(n,m,\delta_0, \alpha)$ such that $C(n,m,\delta_0)\eta_2^{2(1-\alpha)} \le 1/20$.
\item We let $\kappa_3(n,m,\delta_0)$ be such that
\begin{align*}
C(n,m,\delta_0)\eta_2^{-m-2}\kappa_3 \le \frac{\eta_2^{2\alpha}}{20}.
\end{align*}
\item Furthermore, we may assume that $2(12(n-m)+4)C_{\ref{l:apply 3.4}}\gamma(x, \rho)^{\frac{1}{m+1}}\eta_1^{-1} \le \kappa_4$ is sufficiently small so that 
\begin{align*}
C(n,m,\delta_0)\kappa_4^2\eta_2^{-m} \le \frac{\eta_2^{2\alpha}}{20}
\end{align*} for all $0<\rho \le \rho_0$.
\item Further, we assume that $0<\rho_0$ is sufficiently small such that $\tilde{\gamma}(x, \rho)^2$ satisfies
\begin{align*}
C(n,m,\delta_0)\tilde{\gamma}(x, \rho)^2 \le \frac{\eta_2^{2\alpha}}{20}.
\end{align*}
for all $0<\rho \le \rho_0$.
\item (Density of $\Theta^m(\norm{V}, \cdot)$) For $0<\kappa_5, \kappa_6$ small we may assume that 
\begin{align*}
    \norm{V}(B_\rho^n(x) \setminus \{z: |\Theta^m(\norm{V}, z) - \Theta^m(\norm{V}, x)| \ge \kappa_5\}) \le \kappa_6 \omega_m \rho^m
\end{align*}
for all $0<\rho\le\rho_0$.
\end{enumerate} 

With this choice of $0<\rho_0(m,n,\delta_0, \kappa_5, \kappa_6, x)$ we inspecting the proof of Corollary \ref{c:estimates on lipschitz f} applied in $B_{\rho}^n(x)$ for $0<\rho\le \rho_0$. We see that
we may obtain in the place of \eqref{e: weak harmonic error 2} the following
\begin{align*}
    &\left|(\rho/2)^{-m}\int_{B^m_{\rho/2}(T(x)) \cap X_1} \nabla f^j(z) \cdot \nabla \zeta_1(z) (\Theta^m(V, F(x)) - \Theta^m(V; F(z)) J(F(z)))d\mathscr{L}^m(z)\right| \\
    & \qquad \le \kappa_5(\rho/2)^{-m}\int_{B^m_{\rho/2}(T(x))} |\nabla f^j|d\mathscr{L}^m (1+L^2)\sup |\nabla \zeta_1| + 2^m\kappa_6\gamma_{\Theta}L\sup |\nabla \zeta_1|\\  
    & \qquad \le \Big[(1+L^2) \Gamma_{\ref{c:estimates on lipschitz f}, 1}^{1/2} \tilde{\gamma}\kappa_5 + 2^m\kappa_6\gamma_{\Theta}L\Big]
\sup |\nabla \zeta_1|.
\end{align*}
This gives in place of \eqref{f zeta bound}
\begin{align*}
(\rho/2)^{-m}\int_{B^m_{\rho/2}(T(x))} \nabla f^j \cdot \nabla \zeta_1 d\mathscr{L}^m & \le \Big[\Gamma_{\ref{c:estimates on lipschitz f}, 2}\tilde{\gamma}^2 + (1+L^2) \Gamma_{\ref{c:estimates on lipschitz f}, 1}^{1/2} \tilde{\gamma}\kappa_5 + 2^m\kappa_6\gamma_{\Theta}L\Big]\sup |\nabla \zeta_1|.
\end{align*}
Thus, for $\kappa_3$ as above and $\epsilon(n,m,\delta_0, \kappa_3)>0$ as in Lemma \ref{l:harmonic approximation}, we choose $\tilde{\gamma}, \kappa_5$, and $\kappa_6$ such that
\begin{align*}
    \frac{\Gamma_{\ref{c:estimates on lipschitz f}, 2}}{\Gamma_{\ref{c:estimates on lipschitz f}, 1}^{1/2}}\tilde{\gamma} \le \epsilon(n, m, \delta_0, \kappa_3)/3, \quad (1+L^2)\kappa_5 \le \epsilon(n, m, \delta_0, \kappa_3)/3,\quad 2^mL\kappa_6\frac{\gamma_{\Theta}}{\Gamma_{\ref{c:estimates on lipschitz f}, 1}^{1/2}\tilde{\gamma}} \le \epsilon(n, m, \delta_0, \kappa_3)/3.
\end{align*}

Whence, following the proof of Lemma \ref{l:7 rect} applied in $B_{\rho_0}^n(x)$ allows us to estimate as in \eqref{e: tilt excess 2} as follows.
\begin{align*}
    C(n,m,\delta_0)\eta_2^2 \le \frac{\eta_2^{2\alpha}}{20},\qquad C(n,m,\delta_0)\eta_2^{-m-2}\kappa_3 \le \frac{\eta_2^{2\alpha}}{20}\\
C(n,m,\delta_0)\kappa_4^2\eta_2^{-m} \le \frac{\eta_2^{2\alpha}}{20}, \qquad C(n,m,\delta_0)\tilde{\gamma}^2 \le \frac{\eta_2^{2\alpha}}{20}\\
    E(x, \eta_2\rho_0, T)^2 \le \frac{\eta_2^{2\alpha}}{4}\tilde{\gamma}^2.
\end{align*}

Therefore, iterating this estimate as in Proof of Corollary \ref{t: VMO main theorem rect}, above, we see that for $r_i = (\eta_2)^i\rho_0$ and $i \in \mathbb{N}$, $V$ satisfies the hypotheses of the $(\delta_0, \tilde{\gamma}(x, r_i))$-Allard tilt-excess decay regime in $B_{r_i}^n(x)$, where
\begin{align*}
    \tilde{\gamma}(x, r_i) = \max\{\frac{\eta_2^{i\alpha}}{2^i}\tilde{\gamma}(x, \rho_0), \Gamma^{1/2} r_i^{m/2} \}.
\end{align*}
Since as $i \rightarrow \infty$ the sequence $\Gamma^{1/2} r_i^{m/2}$ decreases faster than $(\eta_2^{\alpha}/2)^{i}\tilde{\gamma}(x, \rho_0)$ for $m\ge 2$ we obtain the estimate 
\begin{align*}
    \lim_{i \rightarrow \infty} \inf_{R \in G(n,m)} E(x, r_i, R) & \le \lim_{i \rightarrow \infty} \left(\frac{\eta_2^\alpha}{2}\right)^{i} \tilde{\gamma}(x, \rho_0)\\
   & = \lim_{i \rightarrow \infty}(\rho_0^{-\alpha}2^{-i})r_i^{\alpha}\tilde{\gamma}(x, \rho_0).
\end{align*}
Because $\rho_0^{-\alpha}2^{-i} \rightarrow 0$ as $i \rightarrow \infty$, \eqref{e:tilt excess at all scales} gives the conclusion of the Porism.

When $m=1$, we choose $\eta_2 = 1/8$ and $\kappa_3, \kappa_4, \kappa_5, \kappa_6$ and $\rho_0$ such that 
\begin{align*}
C(n,m,\delta_0)\eta_2^{-m-2}\kappa_3 \le \frac{\eta_2^{2}}{20}, \qquad
C(n,m,\delta_0)\kappa_4^2\eta_2^{-m} \le \frac{\eta_2^{2}}{20}, \qquad C(n,m,\delta_0)\tilde{\gamma}^2 \le \frac{\eta_2^{2}}{20}\\
    E(x, \eta_2\rho_0, T)^2 \le \frac{\eta_2^{2}}{4}\tilde{\gamma}^2.
\end{align*}
Following the rest of the argument above gives the porism. \qed

\subsection{Proof of Corollary \ref{c:chord arc}}

Let $0<\tau$ be given.  The tilt-excess estimate in the definition of an $m$-dimensional chord arc varifold with constant $\tau$ may be satisfied by \ref{e:tilt excess at all scales} for sufficiently small $\tilde{\gamma}(n,m,\delta_0, \tau)$.  Moreover, since this estimate is realized for $R_{x, j-1}$, by the proof of Lemma \ref{l:7 rect} we have the one-sided Reifenberg condition
\begin{align}\label{e:one-sided Reif}
   \text{spt}\norm{V} \cap B^n_{r_j}(x) \subset B^n_{C(n,m,\delta_0)\tilde{\gamma}^{\frac{1}{m+1}}\eta_1 r_{i-1}}(R_{x,j-1}) = B^n_{\kappa_4r_{i}}(R_{x,j-1}).
\end{align}
To obtain the other inclusion, we note that by Corollary \ref{c:apply Lip approx rect} we may obtain the estimate of Lemma \ref{l:menne 2.18}(3) in $B^n_{r_{i-1}}(x)$.  According to Remark \ref{r:B(delta)}(1), this gives
\begin{align*}
    \mathscr{L}^m(C(\epsilon_1)) \le \Gamma_{(3)}C(n)\epsilon_1^{-2}\tilde{\gamma}(x, r_{i-1})^2r_{i-1}^m,
\end{align*}
which, in combination with \eqref{e:one-sided Reif}, gives
\begin{align*}
    \sup_{x \in \text{spt}\norm{V} \cap B_{r_{i-1}}^n(0)} \dist(x, T \cap B_{r_{i-1}}^n(x)) \le \left[\kappa_4 + (\Gamma_{(3)}C(n)\epsilon_1^{-2}\tilde{\gamma}(x, r_{i-1})^2)^{1/m}\right]r_{i-1}.
\end{align*}
By choosing $\gamma(n,m,\delta_0, \tau)$ sufficiently small so that $\kappa_4 + (\Gamma_{(3)}C(n)\epsilon_1^{-2}\tilde{\gamma}^2)^{1/m} \le \tau$ the Reifenberg condition is satisfied.

To satisfy the Ahlfors condition, we note that we may take $\tilde{\gamma} \le \gamma_0(m,\tau)$ in Lemma \ref{l:menne sobolev 2.4} to obtain the lower bound.  For the upper bound, we see that if \begin{align*}
    (1+\gamma_{\Theta})(1+L^2)^{1/2}+\Gamma_{(3)}C(n)\epsilon_1^{-2}2^m \tilde{\gamma}^2 \le 1+\tau
\end{align*}
then Lemma \ref{l:improved density bounds} gives the upper bound. Note that this may require choosing $L$ small depending upon $\tau$, as well.  In this case, a new $L$ gives a new $\epsilon_1$.  The $\gamma(n,m,\delta_0, \tau)$ which satisfies the corollary must be chosen compared to this new $\epsilon_1$. \qed

\bibliography{references}
\bibliographystyle{amsalpha}

\end{document}